\documentclass[12pt]{article}
\usepackage{array}
\usepackage{amssymb}
\usepackage{amsmath}
\usepackage{ntheorem}
\usepackage{epsf}
\usepackage{color}
\usepackage{xcolor}
\usepackage{graphicx}
\usepackage{epstopdf}
\usepackage{subfigure}

\usepackage{multirow}
\usepackage{longtable}
\usepackage{threeparttable}
\usepackage{url}
\usepackage{makecell}
\usepackage{booktabs}

\usepackage{setspace}

\usepackage{listings}
\usepackage[ruled, linesnumbered]{algorithm2e}
\def\bfm#1{\mbox{\boldmath$#1$}}
\numberwithin{equation}{section}

\begin{document}

\voffset -0.5truecm \hoffset -1.5truecm


\theoremstyle{plain}
\newtheorem{corollary}{Corollary}[section]
\newtheorem{definition}{Definition}[section]
\newtheorem{example}{Example}[section]
\newtheorem{lemma}{Lemma}
\newtheorem{property}{Property}[section]
\newtheorem{remark}{Remark}[section]
\newtheorem{theorem}{Theorem}
\newtheorem{Table}{Table}
\newtheorem{proof}{Proof}


\newcommand{\bbB}{{\mathbb{B}}}
\newcommand{\bbC}{{\mathbb{C}}}
\newcommand{\bbD}{{\mathbb{D}}}
\newcommand{\bbE}{{\mathbb{E}}}
\newcommand{\bbI}{{\mathbb{I}}}
\newcommand{\bbJ}{{\mathbb{J}}}
\newcommand{\bbK}{{\mathbb{K}}}
\newcommand{\bbN}{{\mathbb{N}}}
\newcommand{\bbO}{{\mathbb{O}}}
\newcommand{\bbP}{{\mathbb{P}}}
\newcommand{\bbR}{{\mathbb{R}}}
\newcommand{\bbS}{{\mathbb{S}}}
\newcommand{\bbT}{{\mathbb{T}}}
\newcommand{\bbV}{{\mathbb{V}}}
\newcommand{\bbX}{{\mathbb{X}}}
\newcommand{\bbY}{{\mathbb{Y}}}


\newcommand{\cA}{{\cal A}}         \newcommand{\dA}[1]{{\cal A}_{#1}}
\newcommand{\cB}{{\cal B}}         \newcommand{\dB}[1]{{\cal B}_{#1}}
\newcommand{\cI}{{\cal I}}
\newcommand{\cJ}{{\cal J}}
\newcommand{\cK}{{\cal K}}
\newcommand{\cN}{{\cal N}}
\newcommand{\cO}{{\cal O}}
\newcommand{\cP}{{\cal P}}
\newcommand{\cS}{{\cal S}}         \newcommand{\dS}[1]{{\cal S}_{#1}}
\newcommand{\cR}{{\cal R}}
\newcommand{\cT}{{\cal T}}
\newcommand{\cU}{{\cal U}}
\newcommand{\cV}{{\cal V}}
\newcommand{\cX}{{\cal X}}
\newcommand{\cY}{{\cal Y}}
\newcommand{\cZ}{{\cal Z}}


\newcommand{\f}[1]{f_{_{\tiny\rm #1}}}  \newcommand{\hf}[1]{\hat{f}_{#1}}
\newcommand{\tf}[2]{f_{{#1}}^{{#2}}}

\newcommand{\h}[1]{h_{_{\tiny\rm #1}}}
\newcommand{\n}[1]{n_{_{\tiny\rm #1}}}
\newcommand{\q}[2]{q_{_{\tiny\rm #1\!, \, #2}}}
\newcommand{\del}[1]{\delta_{_{\tiny\rm #1}}}
\newcommand{\pii}[2]{\pi_{_{\tiny\rm #1\!, \, #2}}}

\newcommand{\y}[2]{y_{#1}^{\rm #2}}
\newcommand{\yb}[1]{{\bar{y}^{\rm #1}}}
\newcommand{\Y}[1]{Y_{\rm #1}}
\newcommand{\YY}[1]{Y_{\! \mbox{\scriptsize #1}}}

\newcommand{\no}[2] {   \| #1 \|_{_{#2}}   }
\newcommand{\nor}[3]{\|#1\|_{_{#2}}^{#3}}

\newcommand{\bismall}[2]{(\!\!\begin{array}{c} #1 \\ \vspace*{-0.6cm} \\ #2 \end{array} \!\! )}
\newcommand{\bi}[2]{\Big(\!\!\begin{array}{c} #1 \\ #2 \end{array} \!\!\Big)}
\newcommand{\twoc}[2]{\bigg(\!\!\begin{array}{c} #1 \\ #2  \end{array} \!\!\bigg)}
\newcommand{\twolr}[2]{\left(\!\!\begin{array}{c} #1 \\ #2  \end{array} \!\!\right)}

\newcommand{\twob}[2]{\lbrace {#1 \atop #2} \rbrace}
\newcommand{\twoB}[2]{\bigg\{\!\!\begin{array}{c} #1 \\ #2  \end{array} \!\!\bigg\}}

\newcommand{\fourc}[4]{\bigg(\!\!\begin{array}{cc} #1 & #2 \\ #3 & #4 \end{array} \!\!\bigg)}
\newcommand{\fourlr}[4]{\left(\!\!\begin{array}{cc} #1 & #2 \\ #3 & #4 \end{array} \!\!\right)}
\newcommand{\Cases}[4]{\left\{\!\!\begin{array}{ll} #1, & #2,  \\ [3mm]
                                                    #3, & #4  \end{array} \right. }
\newcommand{\Casesthree}[6]
{\left\{\!\!\begin{array}{ll} #1 & #2  \\ [2mm]
                              #3 & #4  \\ [2mm]
                              #5 & #6
\end{array} \right. }

\newcommand{\threec}[3]{\left(\!\!\begin{array}{c} #1 \\ #2 \\ #3  \end{array} \!\!\right)}
\newcommand{\threetwo}[6]{\left(\!\!\begin{array}{cc}
#1 & #2 \\
#3 & #4 \\
#5 & #6
\end{array} \!\!\right)}

\newcommand{\threethree}[9]{\left(\!\!\begin{array}{ccc}
#1 & #2 & #3 \\
#4 & #5 & #6 \\
#7 & #8 & #9
\end{array} \!\!\right)}

\newcommand{\sr}[1]{\stackrel{#1}{=}}
\newcommand{\srl}[1]{\stackrel{#1}{<}}
\newcommand{\srg}[1]{\stackrel{#1}{>}}
\newcommand{\sle}[1]{\stackrel{#1}{\le}}
\newcommand{\sge}[1]{\stackrel{#1}{\ge}}

\newcommand{\mc}[1]{\multicolumn{1}{c}{#1}}
\newcommand{\MC}[3]{\multicolumn{#1}{#2}{#3}}
\newcommand{\TC}[4]{\contentsline {#1}{\numberline {#2}#3}{#4}}
\newcommand{\LIST}[3]{\contentsline {section}{\numberline {\hspace*{-0.5cm}#1}#2}{#3}}

\newcommand{\Hpi}[1]{\hat{\pi}_{_{\tiny\rm #1}}}
\newcommand{\titem}[1]{\vspace*{-0.15cm}\item[{\rm #1} ]\hspace*{0.1cm}}
\newcommand{\ub}[1]{\underline{\bf #1}}
\newcommand{\RED}[1]{{\color{red} #1}}
\newcommand{\GREEN}[1]{{\color{green} #1}}
\newcommand{\BLUE}[1]{{\color{blue} #1}}

\newcommand{\two}[2]{{#1}_{_{\tiny\rm #2}}}
\newcommand{\three}[3]{{#1}_{_{#2, \tiny\rm #3}}}

\newcommand{\threev}[3]{\left(\!\!\begin{array}{c} #1 \\   #2 \\   #3 \end{array} \!\!\right)}

\newcommand{\ct}[1]{c^{(t)}_{#1}}
\newcommand{\dt}[1]{d^{(t)}_{#1}}
\newcommand{\rt}[1]{r^{(t)}_{#1}}


\newcommand{\0}{{\bf 0}\!\!\!{\bf 0}}
\newcommand{\1}{{\bf 1}\!\!\!{\bf 1}}        \def\b1{{\bf 1\!\!\!1}}  
\newcommand{\2}{{\bf 2}\!\!\!{\bf 2}}        \def\b2{{\bf 2\!\!\!2}}  


\newcommand{\ba}{{\bf a}}
\newcommand{\bb}{{\bf b}}
\newcommand{\bd}{{\bf d}}
\newcommand{\bn}{{\bf n}}
\newcommand{\bm}{{\bf m}}
\newcommand{\bp}{{\bf p}}
\newcommand{\bq}{{\bf q}}
\newcommand{\br}{{\bf r}}
\newcommand{\bs}{{\bf s}}
\newcommand{\bt}{{\bf t}}
\newcommand{\bu}{{\bf u}}
\newcommand{\bv}{{\bf v}}
\newcommand{\bw}{{\bf w}}
\newcommand{\bx}{{\bf x}}
\newcommand{\by}{{\bf y}}
\newcommand{\bz}{{\bf z}}

\newcommand{\bA}{{\bf A}}
\newcommand{\bB}{{\bf B}}
\newcommand{\bC}{{\bf C}}
\newcommand{\bD}{{\bf D}}
\newcommand{\bH}{{\bf H}}
\newcommand{\bI}{{\bf I}}
\newcommand{\bL}{{\bf L}}
\newcommand{\bO}{{\bf O}}
\newcommand{\bP}{{\bf P}}
\newcommand{\bQ}{{\bf Q}}
\newcommand{\bR}{{\bf R}}
\newcommand{\bS}{{\bf S}}
\newcommand{\bT}{{\bf T}}
\newcommand{\bU}{{\bf U}}
\newcommand{\bV}{{\bf V}}
\newcommand{\bW}{{\bf W}}
\newcommand{\bX}{{\bf X}}
\newcommand{\bY}{{\bf Y}}
\newcommand{\bZ}{{\bf Z}}


\newcommand{\iba}{{\boldsymbol{a}}}    \newcommand{\ibA}{{\boldsymbol{A}}}
\newcommand{\ibb}{{\boldsymbol{b}}}    \newcommand{\ibB}{{\boldsymbol{B}}}
                                       \newcommand{\ibC}{{\boldsymbol{C}}}

\newcommand{\ibh}{{\boldsymbol{h}}}    \newcommand{\ibH}{{\boldsymbol{H}}}

\newcommand{\ibI}{{\boldsymbol{I}}}
\newcommand{\ibla}{{\boldsymbol{\lambda}}}
\newcommand{\ibe}{{\bfm e}}
\newcommand{\ibm}{{\bfm m}}
\newcommand{\ibn}{{\bfm n}}             \newcommand{\ibN}{{\boldsymbol{N}}}

\newcommand{\ibp}{{\boldsymbol{p}}}     \newcommand{\ibP}{{\boldsymbol{P}}}
\newcommand{\ibs}{{\boldsymbol{s}}}
\newcommand{\ibt}{{\bfm t}}
\newcommand{\ibu}{{\boldsymbol{u}}}
\newcommand{\ibv}{{\boldsymbol{v}}}
\newcommand{\ibw}{{\boldsymbol{w}}}
\newcommand{\ibx}{{\boldsymbol{x}}}     \newcommand{\ibX}{{\boldsymbol{X}}}
\newcommand{\iby}{{\boldsymbol{y}}}
\newcommand{\ibz}{{\boldsymbol{z}}}


\newcommand{\hb}{\hat{b}}
\newcommand{\hm}{\hat{m}}
\newcommand{\hp}{\hat{p}}
\newcommand{\hq}{\hat{q}}
\newcommand{\hR}{\hat{R}}

\newcommand{\hbr}{\hat{\br}}

\newcommand{\hth}{\hat{\theta}}   \newcommand{\hTh}{\hat{\Theta}}
\newcommand{\hbth}{\boldsymbol{\hth}}
\newcommand{\hvth}{\hat{\vartheta}}
\newcommand{\hmu}{\hat{\mu}}
\newcommand{\hsi}{\hat{\sigma}}   \newcommand{\hSi}{\hat{\Sigma}}
\newcommand{\hal}{\hat{\alpha}}   \newcommand{\hbe}{\hat{\beta}}     \newcommand{\hbbe}{\bfm{\hat \beta}}
\newcommand{\hbal}{\boldsymbol{\hal}}
\newcommand{\hga}{\hat{\gamma}}   \newcommand{\hbga}{\bfm{\hat \gamma}}
\newcommand{\hpsi}{\hat{\psi}}
\newcommand{\hxi}{\hat{\xi}}
\newcommand{\hpi}{\hat{\pi}}      \newcommand{\hbpi}{\boldsymbol{\hpi}}
\newcommand{\hla}{\hat{\lambda}}  \newcommand{\hbla}{\bfm{{\hat \lambda}}}
\newcommand{\hde}{\hat{\delta}}
\newcommand{\hvp}{\hat{\varphi}}
\newcommand{\hrho}{\hat{\rho}}
\newcommand{\hbrho}{\boldsymbol{\hrho}}
\newcommand{\hbphi}{\bfm{\hat \phi}}
\newcommand{\hphi}{\hat \phi}
\newcommand{\hnu}{\hat{\nu}}

\newcommand{\whVar}{\widehat{\mbox{Var}}}
\newcommand{\whse}{\widehat{\mbox{se}}}
\newcommand{\whPr}{\widehat{\Pr}}


                                     \newcommand{\tD}{\tilde{D}}

\newcommand{\tp}{\tilde{p}}
                                     \newcommand{\tQ}{\tilde{Q}}
\newcommand{\tx}{\tilde{x}}
\newcommand{\ty}{\tilde{y}}

\newcommand{\tal}{\tilde{\alpha}}
\newcommand{\tbe}{\tilde{\beta}}
\newcommand{\tth}{\tilde{\theta}}     \newcommand{\tbth}{{\bfm \tth}}
\newcommand{\tpi}{\tilde{\pi}}
\newcommand{\tmu}{\tilde{\mu}}
\newcommand{\tga}{\tilde{\gamma}}
\newcommand{\tsi}{\tilde{\sigma}}     \newcommand{\tSi}{\tilde{\Sigma}}
\newcommand{\tpsi}{\tilde{\psi}}
\newcommand{\txi}{\tilde{\xi}}
\newcommand{\tvp}{\tilde{\varphi}}

\newcommand{\Bt}{\bar{t}}
\newcommand{\Bu}{\bar{u}}
\newcommand{\Bv}{\bar{v}}
\newcommand{\Bw}{\bar{w}}
\newcommand{\Bx}{\bar{x}}       \newcommand{\BX}{\bar{X}}
\newcommand{\By}{\bar{y}}       \newcommand{\BY}{\bar{Y}}
\newcommand{\Bz}{\bar{z}}

\newcommand{\Bxi}{\bar{\xi}}

\newcommand{\BVar}{\overline{\mbox{Var}}}


\newcommand{\al}{\alpha}            \newcommand{\alt}{\alpha^{(t)}}
                                    \newcommand{\altI}{\alpha^{(t+1)}}
\newcommand{\be}{\beta}             \newcommand{\Bet}{\beta^{(t)}}
                                    \newcommand{\betI}{\beta^{(t+1)}}
\newcommand{\ga}{\gamma}            \newcommand{\Ga}{\Gamma}
\newcommand{\gat}{\gamma^{(t)}}
\newcommand{\de}{\delta}            \newcommand{\De}{\Delta}
\newcommand{\la}{\lambda}           \newcommand{\La}{\Lambda}
\newcommand{\si}{\sigma}           \newcommand{\Sig}{\Sigma}     
\newcommand{\tha}{\theta}           \newcommand{\Tha}{\Theta}
\newcommand{\ths}{\theta^*}
\newcommand{\tht}{\theta^{(t)}}     \newcommand{\thtI}{\theta^{(t+1)}}

\newcommand{\thx}{\theta_x}         \newcommand{\Thx}{\Theta_x}
\newcommand{\thy}{\theta_y}
\newcommand{\piit}{\pi_i^{(t)}}
\newcommand{\xt}{x^{(t)}}

\newcommand{\ka}{\kappa}
\newcommand{\om}{\omega}            \newcommand{\Om}{\Omega}

\newcommand{\ve}{\varepsilon}
\newcommand{\vp}{\varphi}
\newcommand{\vr}{\varrho}
\newcommand{\vth}{\vartheta}


\newcommand{\bxi}{{\bfm \xi}}              \newcommand{\bet}{{\bfm \eta}}
\newcommand{\bphi}{{\boldsymbol{\phi}}}    \newcommand{\bphit}{\bphi^{(t)}}

\newcommand{\bmu}{{\boldsymbol{\mu}}}      \newcommand{\bnu}{{\boldsymbol{\nu}}}
\newcommand{\bla}{{\boldsymbol{\lambda}}}  \newcommand{\bLa}{{\bfm \Lambda}}
\newcommand{\bSi}{{\bfm \Sigma}}
\newcommand{\bom}{{\bfm \omega}}   \newcommand{\bOm}{{\bfm \Omega}}
\newcommand{\bde}{{\bfm \delta}}   \newcommand{\bDe}{{\bfm \Delta}}
\newcommand{\bth}{{\boldsymbol{\theta}}} \newcommand{\btht}{\bth^{(t)}}
\newcommand{\bTh}{{\boldsymbol{\Theta}}}
\newcommand{\bsi}{{\bfm \sigma}}
\newcommand{\bal}{{\boldsymbol{\alpha}}}
\newcommand{\bbe}{{\boldsymbol{\beta}}}
\newcommand{\brho}{{\boldsymbol{\rho}}}

\newcommand{\bpsi}{{\bfm \psi}}          \newcommand{\bPsi}{{\bfm \Psi}}
\newcommand{\bpi}{{\boldsymbol{\pi}}}    \newcommand{\bpit}{\bpi^{(t)}}
\newcommand{\bga}{{\boldsymbol{\gamma}}} \newcommand{\bgat}{\bga^{(t)}}
\newcommand{\bvp}{{\bfm \varphi}}


\newcommand{\Bernoulli}{\mbox{Bernoulli}}
\newcommand{\Binomial}{\mbox{Binomial}}
\newcommand{\BBB}{\mbox{Bernoulli}_2^{\rm bk}}
\newcommand{\BBinomial}{\mbox{BBinomial}}
\newcommand{\Categorical}{\mbox{Categorical}}
\newcommand{\CP}{\mbox{CP}}
\newcommand{\D}{\mbox{D}}
\newcommand{\Degenerate}{\mbox{Degenerate}}
\newcommand{\DExponential}{\mbox{DExponential}}
\newcommand{\Dirichlet}{\mbox{\rm Dirichlet}}
\newcommand{\DMultinomial}{\mbox{DMultinomial}}
\newcommand{\Exponential}{\mbox{Exponential}}
\newcommand{\FDiscrete}{\mbox{FDiscrete}}
\newcommand{\Finite}{\mbox{Finite}}
\newcommand{\GD}{\mbox{GD}}
\newcommand{\GDirichlet}{\mbox{\rm GDirichlet}}
\newcommand{\GLiouville}{\mbox{GLiouville}}
\newcommand{\GPoisson}{\mbox{GPoisson}}
\newcommand{\GP}{\mbox{GP}}  \newcommand{\GPI}{\GP^{\mbox{\scriptsize{(I)}}}} \newcommand{\GPII}{\GP^{\mbox{\scriptsize{(II)}}}}
\newcommand{\Hgeometric}{\mbox{Hgeometric}}
\newcommand{\HPP}{\mbox{HPP}}
\newcommand{\IBeta}{\mbox{IBeta}}
\newcommand{\Ichi}{\mbox{I}\raisebox{0.5ex}{$\chi$}}
\newcommand{\IGamma}{\mbox{IGamma}}
\newcommand{\IGaussian}{\mbox{IGaussian}}
\newcommand{\IWishart}{\mbox{IWishart}}
\newcommand{\Laplace}{\mbox{Laplace}}
\newcommand{\Liouville}{\mbox{Liouville}}
\newcommand{\Logistic}{\mbox{Logistic}}
\newcommand{\Lognormal}{\mbox{Lognormal}}
\newcommand{\mBeta}{\mbox{Beta}}
\newcommand{\mGamma}{\mbox{Gamma}}
\newcommand{\Multinomial}{\mbox{Multinomial}}
\newcommand{\NBinomial}{\mbox{NBinomial}}
\newcommand{\ND}{\mbox{ND}}
\newcommand{\NDirichlet}{\mbox{NDirichlet}}
\newcommand{\NHPP}{\mbox{NHPP}}
\newcommand{\Poisson}{\mbox{Poisson}}
\newcommand{\TBeta}{\mbox{TBeta}}
\newcommand{\VZIP}{\mbox{VZIP}}
\newcommand{\Wishart}{\mbox{Wishart}}
\newcommand{\ZAP}{\mbox{ZAP}}
\newcommand{\ZDP}{\mbox{ZDP}}
\newcommand{\Zeta}{\mbox{Zeta}}
\newcommand{\ZIP}{\mbox{ZIP}}
\newcommand{\ZTP}{\mbox{ZTP}}



\newcommand{\rhor}{\rho^{^{\rm real}}}
\newcommand{\rhof}{\rho^{^{\tiny \rm fit}}}

\newcommand{\AIC}{\mbox{AIC}}
\newcommand{\BIC}{\mbox{BIC}}
\newcommand{\Corr}{\mbox{Corr}} \newcommand{\CorrR}{\Corr^{\small \rm{R}}} \newcommand{\CorrF}{\Corr^{\small \rm{fit}}}
\newcommand{\Cov}{\mbox{Cov}}
\newcommand{\CV}{\mbox{CV}}
\newcommand{\data}{\mbox{data}}
\newcommand{\diag}{\mbox{diag}}
\newcommand{\I}{\mbox{I}}
\newcommand{\IG}{\mbox{IG}}
\newcommand{\KL}{\mbox{KL}}
\newcommand{\logit}{\mbox{logit}}
\newcommand{\median}{\mbox{median}}
\newcommand{\MSE}{\mbox{MSE}}
\newcommand{\rank}{\mbox{rank}\,}
\newcommand{\sol}{\mbox{sol}}
\newcommand{\se}{\mbox{se}}
\newcommand{\SE}{\mbox{SE}}
\newcommand{\sgn}{{\rm sgn}}
\newcommand{\tr}{\mbox{$\,$tr$\,$}}
\newcommand{\Var}{\mbox{Var}}

\newcommand{\qand}{\quad \mbox{and} \quad}
\newcommand{\qqand}{\quad \qand \quad}
\newcommand{\qag}{\quad \mbox{against} \quad}
\newcommand{\qas}{\quad \mbox{as} \quad}
\newcommand{\qor}{\quad \mbox{or} \quad}
\newcommand{\qve}{\quad \mbox{versus} \quad}
\newcommand{\qwh}{\quad \mbox{where} \quad}
\newcommand{\qwith}{\quad \mbox{with} \quad}
\newcommand{\col}{\mbox{: }}
\newcommand{\RE}{\mbox{RE}}
\newcommand{\yes}{\mbox{yes}}
\newcommand{\No}{\mbox{no}}
\newcommand{\DPP}{\mbox{DPP}}

\newcommand{\rd}{\,\mbox{d}}
\newcommand{\e}{\mbox{e}}
\newcommand{\w}{\mbox{w}}
\renewcommand{\ge}{\geqslant}
\renewcommand{\le}{\leqslant}
\newcommand{\Zp}{Z\hspace{0.02cm}'}
\newcommand{\mif}{\mbox{if } \ }


\newcommand{\II}{I\hspace*{-0.4mm}I}
\newcommand{\III}{I\!I\!I}
\newcommand{\IR}{{I\!\! R}}
\newcommand{\IV}{I$\!$V}
\newcommand{\et}{{\it et al}.}
\newcommand{\Et}{{\it et al}.\,}

\newcommand{\yikH}{y_{ik}^{\rm H}}
\newcommand{\ybH}{\bar{y}^{\rm H}}
\renewcommand{\1}{\uppercase\expandafter{\romannumeral1}}
\renewcommand{\2}{\uppercase\expandafter{\romannumeral2}}


\newcommand{\sd}{\stackrel{{\rm d}}{=}}
\newcommand{\sdt}{$ {\small $\sd$} $}
\newcommand{\heq}{\;\hat{=}\;}
\newcommand{\teq}{\triangleq}
\newcommand{\iid}{\stackrel{{\rm iid}}{\sim}}
\newcommand{\tiid}{i.i.d.$\hspace*{0.08cm}$}
\newcommand{\ind}{\stackrel{{\rm ind}}{\sim}}
\newcommand{\dsim}{\stackrel{.}{\sim}}
\newcommand{\dis}{\displaystyle}
\newcommand{\tex}{\textstyle}
\newcommand{\cf}{cf.$\hspace*{0.1cm}$}

\newcommand{\T}{\!\top\!}
\newcommand{\na}{\nabla}
\newcommand{\noi}{\noindent}
\newcommand{\ra}{\rightarrow}
\newcommand{\pr}{\propto}
\newcommand{\eq}{\equiv}
\newcommand{\pa}{\partial}
\newcommand{\ol}{\overline}
\newcommand{\non}{\nonumber}
\newcommand{\ap}{\approx}
\newcommand{\Bot}{\;\bot\;}
\newcommand{\inde}{{\Bot\!\!\!\!\!\!\!\Bot}}
\newcommand{\btu}{\bigtriangleup}


\newcommand{\vs}{\vspace*{-0.25cm}}
\newcommand{\vkl}{\vskip 0.10in}
\newcommand{\vkL}{\vskip 0.15in}
\newcommand{\vkU}{\vskip 0.30in}


\newcommand{\namelistlabel}[1]{\mbox{#1}\hfil}
\newenvironment{namelist}[1]{%
\begin{list}{}
       {\let \makelabel \namelistlabel
        \settowidth{\labelwidth}{#1}
        \setlength{\leftmargin}{1.1\labelwidth}   }
        }{%
\end{list} }

\def\bds{\begin{description} \itemsep=-\parsep \itemindent=-0.7 cm}
\def\eds{\end{description}}
\def\i{\item}

\newcommand{\tphi}{\tilde{\phi}}
\newcommand{\tla}{\tilde{\lambda}}
\newcommand{\MP}{\mbox{MP}}
\newcommand{\twocc}[2]{\bigg\{\!\!\begin{array}{c} #1 \\ #2  \end{array} \!\!\bigg\}}
\newcommand{\bibx}{{\boldsymbol{\bar{x}}}}
\newcommand{\tbmu}{{\boldsymbol{\tilde{\mu}}}}
\newcommand{\gD}{\mbox{\rm gD}}
\newcommand{\HDirichlet}{\mbox{\rm HDirichlet}}
\newcommand{\hbp}{\bfm {\hat p}}


\newcommand{\Bpi}{\bar{\pi}}
\newcommand{\Bal}{\bar{\al}}
\newcommand{\Bbe}{\bar{\be}}
\newcommand{\p}[1]{p_{\rm #1}}
\newcommand{\s}[1]{s_{\rm #1}}
\newcommand{\GA}{\mbox{GA}}
\newcommand{\NR}{\mbox{NR}}

\doublespacing

\begin{flushleft}
{\Large \bf  The upper--crossing/solution (US) algorithm with strongly stable convergence and an acceleration technique for root--finding problems}
\end{flushleft}

\smallskip
\begin{flushleft}
{\bf Xun-Jian LI$^{{\rm a}}$},   \ \
{\bf Hua ZHOU$^{{\rm a,b}}$},             \ \
{\bf Kenneth LANGE$^{{\rm b,c}}$}         \ \  and \ \
{\bf Guo-Liang TIAN$^{{\rm d},*}$}

\vkl
{\small \it $^{\rm a}$Department of Biostatistics,
            University of California, Los Angeles, CA 90095, USA} \\
\textsf{Email:  xunjianli@ucla.edu}

\vkl
{\small \it $^{\rm b}$Department of Computational Medicine,
           University of California, Los Angeles, CA 90095, USA}
\textsf{Email:  huazhou@ucla.edu}

\vkl
{\small \it $^{\rm c}$Department of Human Genetics, Statistics,
           University of California, Los Angeles,} \\
{\small \it CA 90095, USA}  \ \
\textsf{Email:  klange@ucla.edu}

\vkL
{\small \it $^{\rm d}$Department of Statistics and Data Science,
            Southern University of Science and Technology,} \\
{\small \it Shenzhen 518055, Guangdong Province, P. R. China} \ \
\textsf{Email:  tiangl@sustech.edu.cn}

\vkL
{\small  $^*$Corresponding author's emails: \textsf{xunjianli@ucla.edu, \ \ tiangl@sustech.edu.cn}} \\
\end{flushleft}

\noi {\bf Abstract}. In this paper, we propose a new and broadly applicable root--finding method, called as the \textit{upper--crossing/solution} (US) algorithm, which belongs to the category of non-bracketing (or open domain) methods. The US algorithm is a general principle for iteratively seeking the unique root $\ths$ of a non-linear equation $g(\tha)=0$ and its each iteration consists of two steps: An upper--crossing step (\textsf{U-step}) and a solution step (\textsf{S-step}), where the \textsf{U-step} finds an upper--crossing function or a $U$-function $U(\tha|\tht)$ [whose form depends on $\tht$ being the $t$-th iteration of $\ths$] based on a new notion of so-called changing direction inequality, and the \textsf{S-step} solves the simple $U$-equation $U(\tha|\tht) =0$ to obtain its explicit solution $\thtI$. The US algorithm holds two major advantages: (i) It strongly stably converges to the root $\ths$; and (ii) it does not depend on any initial values, in contrast to Newton's method. The key step for applying the US algorithm is to construct one simple $U$-function $U(\tha|\tht)$ such that an explicit solution to the $U$-equation $U(\tha|\tht) =0$ is available. Based on the first--, second--, third-- and block--derivative of $g(\tha)$, four methods are given for establishing such $U$-functions. Various applications in statistics of the proposed US algorithm and the analysis of the convergence rate of the US algorithm are provided. Furthermore, we develop an acceleration technique for the US algorithm, resulting in a weakly stable convergence. Some numerical experiments and comparisons are also presented. Especially, because of the property of strongly stable convergence, the US algorithm could be one of the powerful tools for solving an equation with multiple roots.

\vkL \noi {\bf Keywords}: Acceleration technique; Changing direction inequality; Non-bracketing methods; Strongly/weakly stable convergence; US algorithm.

\baselineskip 0.30in
\setcounter{equation}{0}
\section{$\!\!\!\!\!\!\!$. Introduction}    

Solving roots of an equation or zeros of a function is the most fundamental problem in computational mathematics. Many problems in science and engineering can be formulated as: Given a continuous function $g(x)$ of a single variable defined in the real line $\bbR \teq (-\infty, \infty)$, find the value $x^* \in \bbR$ such that $g(x^*)=0$. These problems are usually called root--finding problems. Calculating the \textit{maximum likelihood estimator} (MLE) of a parameter in classical statistics and computing the posterior mode in Bayesian statistics can be viewed as finding a solution to a non-linear equation. Let $g(\tha)$ be a non-linear function of a single variable and there exist a unique root $\ths$ to the following equation:
\begin{eqnarray} \label{uclaone1.1}
    g(\tha) = 0, \quad \tha\in \Tha \subseteq \bbR,
\end{eqnarray}
where $\Tha$ is the parameter space or an interval in the real line $\bbR$.

There exist many methods for seeking the root of  (\ref{uclaone1.1}) in numerical analysis and they are roughly categorized into two classes:
\begin{namelist}{0123}
\item[\hspace*{-0.03cm} (i)] Closed domain (bracketing) methods, for example,  the bisection method and the false position method (Wood 1992, Chabert 1999 p.83--112).

\item[\hspace*{-0.03cm} (ii)] Open domain (non-bracketing) methods such as  Newton's method, the secant method, and Muller's method (Wolfe 1959, Boyd \& Vandenberghe 2004 p.484--496, Costabile \Et 2006).
\end{namelist}
The bisection method (Wood 1992) can guarantee to obtain a solution since the root is trapped in the closed interval, while it converges slowly since the information about actual functions behavior is not utilized. By using information about the function $g(\tha)$, the false position method (Chabert 1999) generally converges more rapidly than the bisection method, however it does not give a bound on the error of the root. Non-bracketing methods are not as robust as bracketing methods and may actually diverge. However, they use information about the function itself to refine the approximations of the root. Thus, they are considerably more efficient than bracketing ones.

It is well known that Newton's method or the \textit{Newton--Raphson} (NR) algorithm is a golden standard for seeking the root of (\ref{uclaone1.1}) with the following NR iteration:
\begin{eqnarray} \label{uclaone1.2}
   \thtI = \tht - \frac{ g(\tht) }{ g'(\tht) }, \quad t=0, 1, \dots, \infty,
\end{eqnarray}
where $\tht$ denotes the $t$-th iteration/approximation of the root $\ths$ and $g'(\tha)$ is the first--order derivative of $g(\tha)$. Although the NR algorithm converges quadratically, it is sensitive to initial values (Lindstrom \& Bates 1988). This property hinders its applications in repeated use of the NR algorithm such as in bootstrapping calculation of the standard deviation of an MLE. One advantage of Newton's method is that it is an optimization technique as well as a root--finding technique. In computational statistics, we know that the \textit{expectation--maximization} (EM) algorithm (Dempster \Et 1977) and \textit{minorization--maximization} (MM) algorithm (Lange \Et 2000, Hunter \& Lange 2004) are two popular optimization tools for calculating MLEs of parameters with monotonic convergence, rather than two root--finding tools.

To overcome the drawback of sensitiveness to initial values associated with the NR algorithm and to retain a similar property to the monotonic convergence associated with the EM/MM algorithms, we in this paper propose a new root--finding method, called as the \textit{upper--crossing/solution} (US) algorithm, which belongs to the category of open domain methods. The US algorithm is a general principle for iteratively seeking the root $\ths$ of equation (\ref{uclaone1.1}) with strongly stable convergence (its definition will be given in the beginning of Section 2), and its each iteration consists of an upper--crossing step (\textsf{U-step}) and a solution step (\textsf{S-step}), where
\begin{namelist}{0123456}
\item[\hspace*{-0.03cm} \textsf{U-step}:] Find a $U$-function $U(\tha|\tht)$ satisfying (\ref{uclaone2.5});
\item[\hspace*{-0.03cm} \textsf{S-step}:] Solve the $U$-equation to obtain its root $\thtI$ as showed by (\ref{uclaone2.6}).
\end{namelist}
The two-step process is repeated until convergence occurs.

The rest of the paper is organized as follows. In Section 2, first we define two new relation symbols on the changing direction inequality, and the concept of weakly/strongly stable convergence. Next we introduce new notions such as  upper--crossing function or $U$-function and $U$-equation. Finally each US iteration with two steps (i.e., \textsf{U-step} and \textsf{S-step}) is formulated and the strongly stable convergence of the US algorithm is proved. Four methods for constructing $U$-functions are given in Section 3. Section 4 presents  three applications, including quantile calculation in continuous distributions and calculating MLEs of parameters in two discrete  distributions. In addition, other more applications are provided in the supplementary material. Section 5 provides the analysis of the convergence rate of the US algorithm. Section 6 develops an acceleration technique for the US algorithm, resulting in a weakly stable convergence. Numerical experiments, comparisons and some discussions are given in Sections 7--8. More technical details are put in five Appendices. Here we first present one simple example to give the flavour of the US algorithm.

\vkL \noi \textbf{Example 1} \textsf{(Root of an equation involving a  trigonometric function)}. Suppose that we want to find the root of the equation
$g(x) \teq \cos (0.5 \pi x) - x = 0$ for all $x \in \bbR$. From (\ref{uclaone3.2}), note that
$$
   g'(x) = -0.5 \pi\sin (\pi x/2)-1 \ge - 0.5 \pi -1 \teq b_1,
$$
then the \textsf{U-step} of the US iteration is to find an  upper--crossing function or a $U$-function $U(x|\xt)$, in this example, it is given by (\ref{uclaone3.5}); i.e.,
$$
   U(x|\xt)= g(\xt) + b_1 (x - \xt) = g(\xt) - (0.5 \pi +1)x + (0.5 \pi +1)\xt.
$$
And the \textsf{S-step} is to solve the $U$-equation: $U(x|\xt) =0$, and to obtain its root $x^{(t+1)}$ as
\begin{eqnarray} \label{uclaone1.3}
   x^{(t+1)} = \xt - \frac{g(\xt)}{b_1} =   \frac{(0.5\pi + 1) x^{(t)} + g(x^{(t)})} {0.5\pi + 1}.
\end{eqnarray}
Figure 1 plots both $g(x)$ and $U(x|x^{(t)})$ with the root $x^*=0.594612$. Using the initial value $x^{(0)} = -1$, Figure 1(a) shows that the sequence $\{\xt\}_{t=0}^\infty$ generated by (\ref{uclaone1.3}) strongly stably converges to $x^*$ satisfying $x^{(0)} < x^{(1)} < \cdots < \xt < \cdots \le x^*$; while Figure 1(b) shows that with $x^{(0)}= 2$, the sequence $\{\xt\}_{t=0}^\infty$ generated by (\ref{uclaone1.3}) strongly stably converges to $x^*$ satisfying $x^* \le \cdots < \xt < \cdots < x^{(1)} < x^{(0)}$.

Table 1 displayed these numerical results of the US iterations for two different initial values $x^{(0)}=-1,2$. With the rate of convergence being 0.1069, the US algorithm moved for 10 and 8 steps, respectively. \hfill $\|$

\begin{figure}
  \begin{center}
     $\,$\vskip -3.0cm     \includegraphics[scale=0.6]{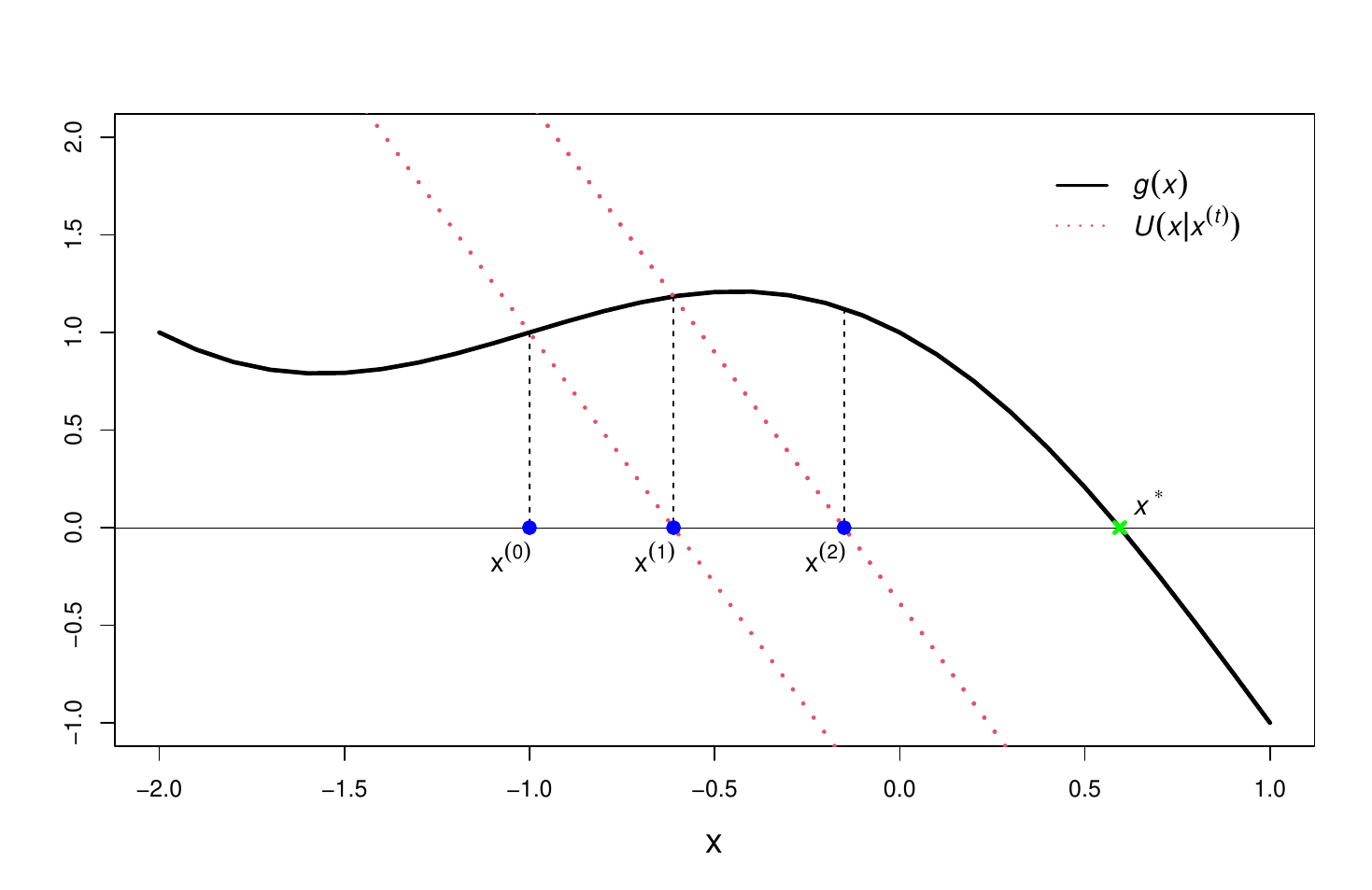}
     \vskip -0.0cm \baselineskip 0.10in {\small  (a)}
     \vskip 0.1cm          \includegraphics[scale=0.6]{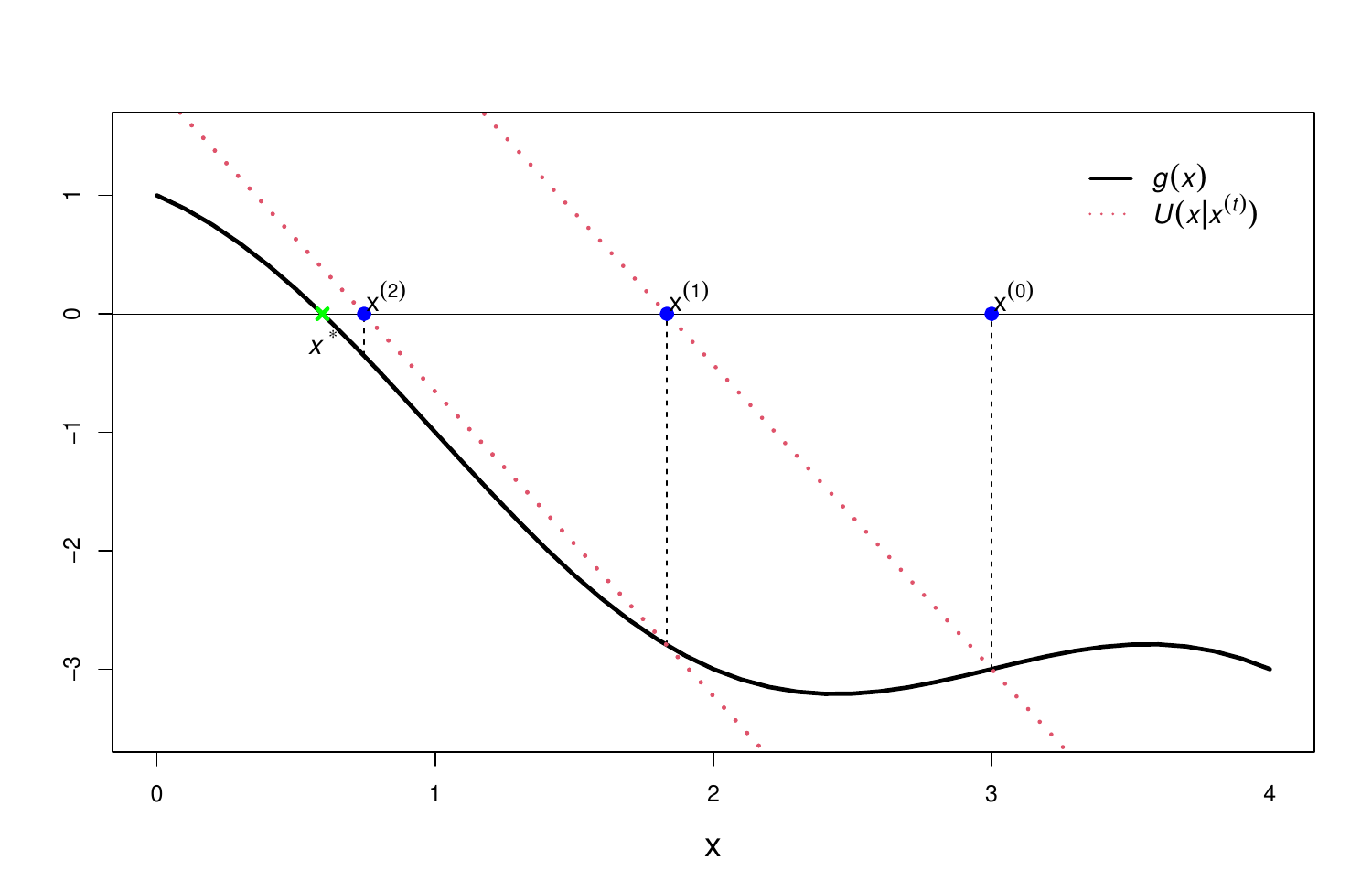}
     \vskip -0.0cm \baselineskip 0.10in {\small  (b)}
  \end{center}
  \vskip -0.0cm \baselineskip 0.10in {\small {\bf Figure 1.} $\;$ Plots of $g(x) = \cos (0.5 \pi x) - x$ and $U(x|x^{(t)}) = - (0.5 \pi +1)x + (0.5 \pi +1)\xt + g(\xt)$. The true root of the equation $g(x)=0$ is $x^*=0.594612$. (a) With the initial value $x^{(0)}=-1$, the sequence $\{\xt\}_{t=0}^\infty$ generated by (\ref{uclaone1.3}) strongly stably converges to $x^*$ satisfying $x^{(0)} < x^{(1)} < \cdots < \xt < \cdots \le x^*$; (b) With the initial value $x^{(0)}= 2$, the sequence $\{\xt\}_{t=0}^\infty$ generated by (\ref{uclaone1.3}) strongly stably converges to $x^*$ satisfying $x^* \le \cdots < \xt < \cdots < x^{(1)} < x^{(0)}$. }
\end{figure}

\vkL
\begin{table}[h]
  \baselineskip 0.10in
  {\bf Table 1.} $\;$ {\rm The US algorithm for solving the equation: $\cos(0.5 \pi x) - x = 0$ with two different initial values $x^{(0)}=-1,2$ }
  \vspace*{-0.2cm}
  \begin{center}
  \renewcommand{\arraystretch}{1.2} \tabcolsep 0.109in \doublerulesep 0.2pt
  \begin{threeparttable}
  \begin{tabular}{l|cccc|cccc} \hline\hline
   $t$   & $x^{(t)}$      &$g(x^{(t)})$  & $\ve^{(t)}$ & \textsf{Rate} & $x^{(t)}$ &$g(x^{(t)})$& $\ve^{(t)}$ & \textsf{Rate}  \\  \hline
    0     & $-$1    & 1     & $-$1.59461 & .7561 & 2     & $-$3    & 1.40539  & .1697 \\
    1     & $-$.611015 & 1.18471 & $-$1.20563 & .6178 & .833046 & $-$.573792 & .238435 & .0639 \\
    2     & $-$.150180 & 1.12248 & $-$.744791 & .4138 & .609850 & $-$.034652 & .015239 & .1154 \\
    3     & \;\;\;.286449 & .614018 & $-$.308163 & .2249 & .596371 & $-$.003983 & .001759 & .1192 \\
    4     & \;\;\;.525293 & .153170 & $-$.069319 & .1405 & .594821 & $-$.000475 & .000210 & .1197 \\
    5     & \;\;\;.584874 & .021967 & $-$.009738 & .1225 & .594637 & $-$.000057 & .000025 & .1197 \\
    6     & \;\;\;.593418 & .002700 & $-$.001193 & .1201 & .594615 & $-$.000007 & .000003 & .1196 \\
    7     & \;\;\;.594468 & .000324 & $-$.000143 & .1198 & .594612 & $-$.000001 & .000000 & .1182 \\
    8     & \;\;\;.594594 & .000039 & $-$.000017 & .1196 & .594612 & \;\;\;.000000 & .000000 & .1069 \\
    9     & \;\;\;.594610 & .000005 & $-$.000002 & .1182 & --      & --            & --            & --   \\
    10    & \;\;\;.594611 & .000001 & \;\;\;.000000 & .1069 & --      & --            & --            & --   \\ \hline \hline
  \end{tabular}
  {\small Note: $\ve^{(t)} \teq x^{(t)}-x^*$ is the difference between the current iteration and the root;  \textsf{Rate} $\teq |\ve^{(t+1)}/\ve^{(t)}|$ is the  rate of convergence. }
  \end{threeparttable}
  \end{center}
\end{table}

\section{$\!\!\!\!\!\!\!$. The US algorithm with strongly stable convergence}    

First, we define two new relation symbols on the \textit{changing direction} (CD) inequalities,  ``$\sge{\sgn(a)}$" and  ``$\sle{\sgn(a)}$", as follows: For two functions $h_1(x)$ and $h_2(x)$ with the same domain $\bbX$,
$$
   h_1(x) \sge{\sgn(a)} h_2(x)  \mbox{ means} \Casesthree{h_1(x) \ge h_2(x),}{\mif a>0,} {h_1(x) = h_2(x),}{\mif a=0,} {h_1(x) \le h_2(x),}{\mif a < 0,}
$$
and
$$
   h_1(x) \sle{\sgn(a)} h_2(x)  \mbox{ means} \Casesthree{h_1(x) \le h_2(x),}{\mif a>0,} {h_1(x) =h_2(x),}{\mif a=0,} {h_1(x) \ge h_2(x),}{\mif a < 0.}
$$
For example, $x^{2n-1} \sge{\sgn(x)} 0$ for any positive integer $n$; $-\sin(x) \sge{\sgn(x - \pi)} 0$ for any $x \in (0, 2\pi)$; and $\cos(x) \sle{\sgn(x - \pi/2)} 0$ for any $x \in (0, \pi)$ are three illustrative examples.

Next, a sequence $\{\tht\}_{t=0}^\infty$ is called to weakly stably converge to a fixed point $\ths$, if
$$
   |\thtI - \ths| < |\tht - \ths|, \qquad \forall \; t=0, 1, \ldots, \infty.
$$
A sequence $\{\tht\}_{t=0}^\infty$ is called to strongly stably converge to a fixed point $\ths$, if
$$
   \tha^{(0)} < \tha^{(1)} < \cdots < \tht < \cdots \le \ths \quad \qor \quad
   \ths \le \cdots < \tht < \cdots < \tha^{(1)} < \tha^{(0)}.
$$
An algorithm is said to possess \textit{weakly (or strongly) stable convergence} if the sequence $\{\tht\}_{t=0}^\infty$ generated by this algorithm weakly (or strongly) stably converges to a fixed point $\ths$. It is clear that a strongly stable convergence implies a weakly stable convergence.

\subsection{Definition of a $U$-function} 

Suppose that directly solving the unique root $\ths$ of the non-linear equation (\ref{uclaone1.1}) is very difficult. Without loss of generality, we assume that
\begin{eqnarray}  \label{uclaone2.1}
    g(\tha) > 0, \quad \forall \; \tha < \ths   \qand
    g(\tha) < 0, \quad \forall \; \tha > \ths.
\end{eqnarray}
If $g(\tha) < 0$ when $\tha<\ths$, then we can multiplies $-1$ on both sides of (\ref{uclaone1.1}) and obtain a new equation $-g(\tha)=0$, which satisfies the assumption (\ref{uclaone2.1}) and has the same root $\ths$.

The idea of the US algorithm is as follows. Let $\tht$ be the $t$-th iteration of $\ths$ and $U(\tha|\tht)$ be a real-valued function of $\tha$ whose form depends on $\tht$. The function $U(\tha|\tht)$ is called a $U$-function for $g(\tha)$ at $\tha = \tht$ if the following three conditions are satisfied:
\begin{eqnarray}
   g(\tha) &\le& U(\tha|\tht), \qquad    \mif \; \tha<\tht, \label{uclaone2.2} \\ [2mm]
   g(\tht) &=& U(\tht|\tht) \hspace*{0.65cm} \mif \; \tha=\tht,  \label{uclaone2.3} 
\end{eqnarray}
\begin{eqnarray}
   g(\tha) &\ge& U(\tha|\tht), \qquad    \mif \; \tha>\tht.   \label{uclaone2.4}
\end{eqnarray}
The three conditions can be equivalently rewritten as
\begin{eqnarray} \label{uclaone2.5}
   g(\tha) \sge{\sgn(\tha - \tht)} U(\tha|\tht), \quad \forall \; \tha, \tht \in \Tha.
\end{eqnarray}
The conditions (\ref{uclaone2.2})--(\ref{uclaone2.3}) indicate that the $U(\tha|\tht)$ function majorizes the objective function $g(\tha)$ at $\tha=\tht$ for all $\tha\le \tht$, and the conditions
(\ref{uclaone2.4})--(\ref{uclaone2.3}) show that the $U(\tha|\tht)$ function minorizes $g(\tha)$ at $\tha=\tht$ for all $\tha \ge \tht$. The geometric interpretations of the objective function $g(\tha)$ (solid curve) and the $U$-function $U(\tha|\tht)$ (dotted curve) are showed in Figure 1.

The assumption (\ref{uclaone2.1}) implies $0 \sge{\sgn(\tha - \ths)} g(\tha)$ for all $\tha, \ths \in \Tha$. By comparing this CD inequality with (\ref{uclaone2.5}), we conclude that $g(\tha)$ is a $U$-function for 0 at $\tha = \ths$. Other properties on $U$-functions are listed in Property 1 below.

\vkL \noi \textbf{Property 1} \textsf{(Basic properties on $U$-functions)}.
\begin{namelist}{0123}
\item[\hspace*{-0.03cm} (a)] Let $g(\tha)$ be defined in $\Tha$, then $g(\tha)$ is a $U$-function for itself at any inner point in $\Tha$;

\item[\hspace*{-0.03cm} (b)] If $U_1(\tha|\tht)$ and $U_2(\tha|\tht)$ are two $U$-functions for $g(\tha)$ at $\tha=\tht$, then $[U_1(\tha|\tht)+U_2(\tha|\tht)]/2$ is a $U$-function for $g(\tha)$ at $\tha=\tht$;

\item[\hspace*{-0.03cm} (c)] For two functions $g_1(\tha)$ and $g_2(\tha)$ defined in $\Tha$, if $U_1(\tha|\tht)$ is a $U$-function for $g_1(\tha)$ at $\tha=\tht$, then $U_1(\tha|\tht) + g_2(\tha)$ is a $U$-function for $g_1(\tha)+g_2(\tha)$ at $\tha = \tht$.   \hfill $\|$
\end{namelist}

\vkL
\subsection{The $U$-equation and the US algorithm} 

If a $U(\tha|\tht)$ function satisfying (\ref{uclaone2.5}) could be found, then one need only solve the simple \textit{U-equation}: $U(\tha |\tht) = 0$  to obtain its solution $\thtI$, denoted by
\begin{eqnarray} \label{uclaone2.6}
   \thtI = \sol\,\{U(\tha |\tht) = 0, \; \forall \; \tha, \tht \in \Tha\},
\end{eqnarray}
which defines a \textit{US algorithm}. Therefore, the US algorithm is a general principle for iteratively seeking the root $\ths$ with strongly stable convergence, and its each iteration consists of an upper--crossing step (\textsf{U-step}) and a solution step (\textsf{S-step}), where
\begin{namelist}{0123456}
\item[\hspace*{-0.03cm} \textsf{U-step}:] Find a $U$-function $U(\tha|\tht)$ satisfying (\ref{uclaone2.5});
\item[\hspace*{-0.03cm} \textsf{S-step}:] Solve the $U$-equation to obtain its root $\thtI$ as showed by (\ref{uclaone2.6}).
\end{namelist}
The two-step process is repeated until convergence occurs.

In general, the solution $\thtI$ in (\ref{uclaone2.6}) has an explicit expression. Usually, the $U$-equation $U(\tha |\tht) = 0$ could be any equation with analytic solution, e.g., a quadratic equation, even a linear equation. Therefore, the US principle is iteratively solving a sequence of tractable surrogate equations $U(\tha |\tht) = 0$ instead of solving the intractable original non-linear equation $g(\tha)=0$.

Theorem 1 below states that the US algorithm holds two characteristics: (i) It strongly stably converges to the root $\ths$; (ii) It does not depend on any initial values in $\Tha$, in contrast to Newton's method. Its proof is put in Appendix A.

\vkL \noi \textbf{Theorem 1} \textsf{(Strongly stable convergence)}. Let $\ths$ be the unique root of the equation (\ref{uclaone1.1}), and $U(\tha|\tht)$ be one $U$-function for $g(\tha)$ at $\tha=\tht$ in $\Tha$. For any initial value $\tha^{(0)}$ in $\Tha$, if $|U(\tha_1|\tht)-U(\tha_2|\tht)|< L |\tha_1-\tha_2|$ for any $\tha_1$ and $\tha_2$ located between $\tha^{(0)}$ and $\ths$, then the sequence $\{\thtI\}_{t=0}^\infty$ specified by (\ref{uclaone2.6}) strongly stably converges to $\ths$, where $L$ is a lipschitz constant relying on $\tha^{(0)}$, $\ths$ and $U(\tha|\tht)$.  \hfill $\|$

\section{$\!\!\!\!\!\!\!$. Four methods for constructing $U$-functions}

For the same objective function $g(\tha)$, it is possible that there exist several $U$-functions. In other words, for the equation (\ref{uclaone1.1}), it may exist more than one US algorithm. We know that the key for a US algorithm is how to find a $U$-function. In this section, based on the objective function $g(\tha)$, we propose four general methods to construct $U$-functions, which can be achieved by utilizing the Taylor expansion of $g(\tha)$ around $\tha=\tht$ with integral remainder (Lange 2016 p.38):
\begin{eqnarray}\label{uclaone3.1}
   g(\tha) = \sum_{k = 0}^n \frac{(\tha- \tht)^k} {k!}g^{(k)}(\tht)
   + \frac{(\tha- \tht)^{n + 1}}{n!}\int_0^1 g^{(n+1)}(\tht + y (\tha- \tht))(1 - y)^n \rd y,
\end{eqnarray}
where $g^{(0)}(\cdot) \teq g(\cdot)$ and $g^{(k)}(\cdot)$ denotes $k$-th derivative of $g(\cdot)$.

\subsection{The FLB function method} 

Suppose the first--derivative $g'(\cdot)$ exists and it is bounded by some function $b(\cdot)$; i.e.,
\begin{eqnarray} \label{uclaone3.2}
   g'(\tha) \ge b(\tha), \quad \forall \; \tha \in \Tha,
\end{eqnarray}
where $b(\tha)$ is called the \textit{first--derivative lower bound} (FLB) function. In (\ref{uclaone3.1}), let $n=0$, we obtain
\begin{eqnarray}
   g(\tha) &=& g(\tht) + (\tha- \tht) \int_0^1 g'(\tht + y (\tha- \tht)) \rd y \non \\ [2mm]
   &=& g(\tht) + \int_{\tht}^{\tha} g'(z) \rd z \non \\ [2mm]
   &\sge{\sgn(\tha-\tht)}& g(\tht) + \int_{\tht}^{\tha} b(z) \rd z \;\teq\; U(\tha|\tht), \quad \forall \; \tha, \tht \in \Tha,  \label{uclaone3.3}
\end{eqnarray}
where the result $g'(z)-b(z)\sge{(\ref{uclaone3.2})} 0$ is used in the above CD inequality.

\subsubsection{Two special cases of the FLB function $b(\tha)$} 

Let $b(\tha)=b_1 + b_2\tha$ be a linear function. Then the $U$-function defined by (\ref{uclaone3.3}) becomes a quadratic function
\begin{eqnarray} \label{uclaone3.4}
   U(\tha|\tht) = \frac{1}{2} b_2\tha^2 + b_1 \tha + \ct{1} \qwith \ct{1} \teq g(\tht) - \frac{1}{2} b_2\tha^{(t)2} - b_1\tht,
\end{eqnarray}
if $b_2\ne 0$. From (\ref{uclaone2.6}), the $U$-equation $U(\tha|\tht)=0$ has two solutions
\begin{eqnarray*}
   \thtI = \frac{-b_1\pm \sqrt{b_1^2 - 2b_2 \ct{1}}}{b_2},
\end{eqnarray*}
but only one is the desired US iteration.

In particular, if $b_2=0$ (i.e., $g'(\tha)\ge b_1$ for all $\tha \in \Tha$, where $b_1$ is called the FLB constant), then (\ref{uclaone3.4}) becomes to a linear function
\begin{eqnarray} \label{uclaone3.5}
   U(\tha|\tht) = g(\tht) + b_1 (\tha - \tht)
\end{eqnarray}
and the corresponding US iteration is given by
\begin{eqnarray*}
   \thtI = \tht - \frac{g(\tht)}{b_1}, \quad b_1\ne 0,
\end{eqnarray*}
which can be viewed as the NR iteration (\ref{uclaone1.2}) by replacing $g'(\tht)$ with $b_1$.

\subsection{The SLUB constants method} 

Suppose that $g''(\cdot)$ exists and there are two constants $\{b_{21}, b_{22}\}$ (which are called the \textit{second--derivative lower--upper bound} (SLUB) constants) such that
\begin{eqnarray} \label{uclaone3.6}
   b_{21}\le  g''(\tha) \le b_{22}, \quad \forall \; \tha \in \Tha.
\end{eqnarray}
Let $I(\cdot)$ denote the indicator function and define a step constant
\begin{eqnarray} \label{uclaone3.7}
   b_2(\tha |\tht) \teq  b_{22} \cdot I(\tha \le \tht) + b_{21} \cdot I(\tha > \tht).
\end{eqnarray}
In (\ref{uclaone3.1}), let $n=1$, we have
\begin{eqnarray}
   g(\tha &=& g(\tht) +(\tha-\tht)g'(\tht) \non \\[2mm]
   & & +\; (\tha- \tht)^2\int_0^1 g''(\tht + y (\tha- \tht))(1 - y) \rd y \label{uclaone3.8} \\ [2mm]
   &\sge{\sgn(\tha-\tht)}& g(\tht) +(\tha-\tht) g'(\tht) + (\tha- \tht)^2\int_{0}^{1} b_2(\tha|\tht)(1 - y) \rd y \non \\[2mm]
   &=& g(\tht) +(\tha-\tht)g'(\tht) + \frac{b_2(\tha| \tht)}{2}  (\tha-\tht)^2 \;\teq\; U(\tha|\tht), \quad \label{uclaone3.9}
\end{eqnarray}
in which the CD inequality holds because of (\ref{uclaone3.6}) and (\ref{uclaone3.7}).
As shown by Figure 2(a), $U(\tha|\tht)$ in (\ref{uclaone3.9}) is a monotonic decreasing function rather than being a real quadratic function, so we could call it a quasi-quadratic function.

\vkL \noi \textbf{Remark 1} \textsf{(A comment on the SLUB constants method)}. The SLUB constants method depends on the assumption specified by (\ref{uclaone3.6}), which is too strong in practice and may be weakened as $b_{21} \le g''(\tha)$ or $g''(\tha) \le b_{22}$. Two special second--derivative bound constant methods are provided in \S 3.2.1 and \S 3.2.2, respectively. \hfill $\|$

\subsubsection{The SLB constant method} 

Suppose that $g''(\cdot)$ exists and there is a \textit{second--derivative lower bound} (SLB) constant $b_{21}$ such that $b_{21}\le  g''(\tha)$. From (\ref{uclaone3.8}), we have
\begin{eqnarray*}
   g(\tha) \ge  g(\tht) +(\tha-\tht)g'(\tht) + \frac{b_{21}}{2} (\tha-\tht)^2
\end{eqnarray*}
with equality iff $\tha=\tht$. It is easy to verify that
\begin{eqnarray}  \label{uclaone3.10}
   U(\tha|\tht) \teq \left[g(\tht) +(\tha-\tht)g'(\tht) + \frac{b_{21}}{2} (\tha-\tht)^2\right] I(\tha \ge \tht) + g(\tha) I(\tha < \tht), \quad
\end{eqnarray}
is a $U$-function for $g(\tha)$ at $\tha = \tht$. The $U$-function specified by (\ref{uclaone3.10}) is a mergence of a half-quadratic function when $\tha \ge \tht$ with the original $g(\tha)$ when $\tha < \tht$. Because the half-quadratic function has a closed--form zero value for $\tht \le \ths$, we define the SLB constant method by updating
\begin{eqnarray} \label{uclaone3.11}
   \thtI &=& \sol\,\left\{g(\tht) +(\tha-\tht) g'(\tht) + \frac{b_{21}}{2} (\tha-\tht)^2 = 0, \; \forall \; \tha \ge \tht, \ths \ge \tht \right\},  \qquad \non \\ [2mm]
   & = & \left\{   \begin{array}{ll}
   \tht - \dis \frac{g'(\tht) +\sqrt{g^{'2}(\tht) - 2 b_{21} g(\tht)}}{ b_{21} }, & \mif  b_{21} < 0, \\ [2mm]
   \tht - \dis \frac{g(\tht)}{ g'(\tht) }, & \mif  b_{21} = 0.
   \end{array} \right.
\end{eqnarray}
According to the strongly stable convergence property of the US algorithm, we have $\tht \le \ths$ as long as the initial value $\tha^{(0)}$ satisfying $g(\tha^{(0)}) \ge 0$ because $g(\tha) \ge 0$ iff $\tha \le \ths$ for all $\tha \in \Tha$.

\subsubsection{The SUB constant method} 

Suppose that $g''(\cdot)$ exists and there is a \textit{second--derivative upper bound} (SUB) constant $b_{22}$ satisfying inequality $g''(\tha) \le b_{22}$. From (\ref{uclaone3.8}), we have
\begin{eqnarray*}
   g(\tha) \le  g(\tht) +(\tha-\tht)g'(\tht) + \frac{b_{22}}{2} (\tha-\tht)^2,
\end{eqnarray*}
with equality iff $\tha=\tht$. It is easy to verify that
\begin{eqnarray} \label{uclaone3.12}
   U(\tha|\tht) \teq \left[g(\tht) +(\tha-\tht)g'(\tht) + \frac{b_{22}}{2} (\tha-\tht)^2\right] I(\tha \le \tht) + g(\tha) I(\tha > \tht), \quad
\end{eqnarray}
is a $U$-function for $g(\tha)$ at $\tha = \tht$. The $U$-function specified by (\ref{uclaone3.12}) is a mergence of a half-quadratic function when $\tha \le \tht$ with the original $g(\tha)$ when $\tha > \tht$. Because the half-quadratic function has a closed--form zero value for $\tht \ge \ths$, we define the SUB constant method by updating
\begin{eqnarray}
   \thtI & = & \sol\,\left\{g(\tht) +(\tha-\tht) g'(\tht) + \frac{b_{22}}{2} (\tha-\tht)^2 = 0, \; \forall \; \tha \le \tht, \ths \le \tht \right\} \qquad \non 
\end{eqnarray}
\begin{eqnarray} \label{uclaone3.13}
   &=& \left\{   \begin{array}{ll}
   \tht - \dis \frac{g'(\tht) +\sqrt{g^{'2}(\tht) - 2 b_{22} g(\tht)}}{ b_{22} }, & \mif  b_{22} > 0, \\ [2mm]
   \tht - \dis \frac{g(\tht)}{ g'(\tht) }, & \mif  b_{22} = 0.
   \end{array} \right.
\end{eqnarray}
According to the strongly stable convergence property of the US algorithm, we have $\tht \ge \ths$ as long as the initial value $\tha^{(0)}$ satisfying $g(\tha^{(0)}) \le 0$ because $g(\tha) \le 0$ iff $\tha \ge \ths$ for all $\tha \in \Tha$.

\vkL \noi \textbf{Remark 2} \textsf{(Connection between SLB \& SUB methods and Newton's method)}. When $b_{21}=0$, the second iteration in (\ref{uclaone3.11}) is Newton's iteration. Hence, we conclude that if $g''(\tha) \ge 0$, Newton's iteration is strongly stable convergence to the root of the equation $g(\tha) = 0$ provided that the initial value $\tha^{(0)}$ satisfies $g(\tha^{(0)}) \ge 0$.

When $b_{22}=0$, the second iteration in (\ref{uclaone3.13}) is Newton's iteration. Hence, we conclude that if $g''(\tha) \le 0$, Newton's iteration is strongly stable convergence to the root of the equation $g(\tha) = 0$ provided that the initial value $\tha^{(0)}$ satisfies $g(\tha^{(0)}) \le 0$.  \hfill $\|$

\subsection{The TLB constant method}  

Suppose that $g'''(\tha)$ exists and is bounded by some constant $b_3$; i.e.,
\begin{eqnarray} \label{uclaone3.14}
   g'''(\tha) \ge b_3, \quad \forall \; \tha \in \Tha,
\end{eqnarray}
where $b_3$ is called the \textit{third--derivative lower bound} (TLB) constant. In (\ref{uclaone3.1}), let $n=2$, we have
\begin{eqnarray}
   g(\tha) &=& \sum_{k = 0}^2 \frac{(\tha- \tht)^k }{k!} g^{(k)}(\tht)+ \frac{(\tha- \tht)^3}{2}\int_0^1 g'''(\tht + y (\tha- \tht)) (1 - y)^2 \rd y \non \\ [2mm]
   &=& \sum_{k = 0}^2 \frac{(\tha- \tht)^k} {k!}g^{(k)}(\tht) + \frac{(\tha- \tht)^3}{2} \int_0^1 b_3 \times (1 - y)^2 \rd y \non \\ [2mm]
   & &+\; \frac{(\tha- \tht)^3}{2}\int_0^1 \left[ g'''(\tht + y (\tha- \tht)) - b_3 \right] (1 - y)^2 \rd y \non 
\end{eqnarray}
\begin{eqnarray} \label{uclaone3.15}
   &\sge{\sgn(\tha-\tht)}& \sum_{k = 0}^2 \frac{(\tha- \tht)^k} {k!}g^{(k)}(\tht) + \frac{b_3}{6}(\tha-\tht)^3
   \teq U(\tha|\tht), \quad \forall \; \tha, \tht \in \Tha. \qquad
\end{eqnarray}

\subsection{The fixed--block method}  

Assume that there are $m$ univariate functions $\{u_j(\tha)\}_{j=1}^m$ such that $g(\tha)$ can be written as $G(u_1(\tha), \ldots, u_m(\tha))$, where $u_j(\tha)$ is also called the $j$-th \textit{block  function} of $g(\tha)$. If there exist $1 \le j_1 <\cdots< j_r \le m$ with $1\le r < m$, such that
\begin{eqnarray}  \label{uclaone3.16}
   \frac{\pa G(u_1(\tha), \ldots, u_m(\tha))}{\pa u_j(\tha)}  u_j'(\tha) \ge 0, \quad \forall\, j = j_1, \ldots, j_r, \quad \forall\, \tha \in \Tha.
\end{eqnarray}
Without loss of generality, let $j_1=1, \ldots, j_r=r$, then
\begin{eqnarray}  \label{uclaone3.17}
   U(\tha|\tht) \teq U_{1\cdots r}(\tha|\tht) = G(u_1(\tht), \ldots, u_r(\tht), u_{r+1}(\tha), \ldots, u_m(\tha))
\end{eqnarray}
for all $\tha, \tht \in \Tha$, is a $U$-function for $g(\tha)$ at $\tha = \tht$.

\vkL \noi\textbf{Proof of (\ref{uclaone3.17})}. We only prove the conclusion for the case of $r=1$. For the general case, we could prove the conclusion similarly. When $r=1$, we only need to prove that
\begin{eqnarray} \label{uclaone3.18}
   g(\tha) \sge{\sgn(\tha - \tht)} U_1(\tha|\tht), \quad \forall\, \tha, \tht \in \Tha.
\end{eqnarray}
From (\ref{uclaone3.16}), we first consider Case I:
\begin{eqnarray}  \label{uclaone3.19}
   \frac{\pa G(u_1(\tha),\ldots, u_m(\tha))}{\pa u_1(\tha)} \ge 0 \qand u_1'(\tha) \ge 0,
\end{eqnarray}
i.e., $u_1(\cdot)$ is increasing on $\tha$. (i) When $\tha < \tht$, we obtain $u_1(\tha) \le u_1(\tht)$. And then we have
$$
   g(\tha) = G(u_1(\tha), u_2(\tha), \ldots, u_m(\tha)) \sle{(\ref{uclaone3.19})} G(u_1(\tht), u_2(\tha), \ldots, u_m(\tha)) = U_1(\tha|\tht).
$$
(ii) When $\tha > \tht$, we obtain $u_1(\tha) \ge u_1(\tht)$. And then we have
$$
   g(\tha) = G(u_1(\tha), u_2(\tha), \ldots, u_m(\tha)) \sge{(\ref{uclaone3.19})} G(u_1(\tht), u_2(\tha), \ldots, u_m(\tha)) = U_1(\tha|\tht).
$$
(iii) When $\tha = \tht$ , we have $g(\tht) = G(u_1(\tht),\ldots, u_m(\tht)) = U_1(\tht|\tht)$. Thus, for Case I, the CD inequality (\ref{uclaone3.18}) is true. Similarly, the CD inequality (\ref{uclaone3.18}) holds for Case \II: $\pa G(u_1(\tha), \ldots, u_m(\tha))/\pa u_1(\tha) \le 0$ and $u_1'(\tha) \le 0$.   \hfill $\Box$

\vkL For illustration, we first present two simple examples showing how to construct the $U$-function in the US algorithm based on the fixed--block method.

\vkL \noi \textbf{Example 2} \textsf{(Roots of two cubic equations)}. Suppose that we want to find the root of the cubic equation $g(\tha) \teq \tha^3 - 2\tha + 1 = 0$ for all $\tha \ge 0$. First, it is easy to obtain $g(0)=1>0$, indicating that  (\ref{uclaone2.1}) is satisfied. We can write $g(\tha)= u_1(\tha) + u_2(\tha)$ with $u_1(\tha)=\tha^3$ and $u_2(\tha) = - 2\tha + 1$. Since $[\pa g(\tha)/\pa u_1(\tha)] u_1'(\tha) = 1 \times 3\tha^2 \ge 0$, from (\ref{uclaone3.17}), we have $U(\tha|\tht) =  u_1(\tht) + u_2(\tha) = \tha^{(t)3} - 2\tha + 1$. Let $U(\tha|\tht)=0$, we obtain the US iteration: $\thtI= (\tha^{(t)3} + 1)/2$.

Next, assume that we want to find the root of the cubic equation $g(\tha) \teq -\tha^3 + 2\tha + 2 = 0$ for all $\tha \ge 0$. First, it is easy to obtain $g(0)=2 >0$, indicating that  (\ref{uclaone2.1}) is satisfied. We can write $g(\tha)= u_1(\tha) + u_2(\tha)$ with $u_1(\tha)=-\tha^3$ and $u_2(\tha) = 2\tha + 2$. Since $[\pa g(\tha)/\pa u_2(\tha)] u_2'(\tha) = 1 \times 2 > 0$, from (\ref{uclaone3.17}), we have $U(\tha|\tht) =  u_1(\tha) + u_2(\tht) = -\tha^3 + 2\tht + 2$. Let $U(\tha|\tht)=0$, we obtain the US iteration: $\thtI= (2\tht + 2)^{1/3}$.  \hfill $\|$

\vkL \noi \textbf{Remark 3} \textsf{(A comment on the fixed--block method)}. Unlike the FLB function, SLUB constants and TLB constant methods that rely on the first-- to third--derivative of $g(\tha)$, the fixed--block method only requires that some block functions of the $g(\tha)$ are increasing or decreasing with respect to $\tha$ such that (\ref{uclaone3.16}) is satisfied, which largely weakens the assumptions on $g(\tha)$. By fixing these block functions $u_1(\tha), \ldots, u_r(\tha)$ at $\tha=\tht$, we obtain the $U$-function $U_{1\cdots r}(\tha|\tht)$ from $g(\tha) = G(u_1(\tha), \ldots, u_m(\tha))$.                      \hfill $\|$

\vkL
Furthermore, in Section 4.3, we will employ the fixed--block method to calculate the MLEs of parameters in the generalized Poisson distribution. Other cases (e.g., the unit power--logarithmic distribution and the Geeta distribution) are provided in the supplementary material.

\section{$\!\!\!\!\!\!\!$. Applications}  

In this section, we give three models to illustrate the application of the US algorithm. Additional examples demonstrating its applications are provided in the supplementary material due to space limitation.

\subsection{Calculation of quantile in continuous distributions} 

Let $f(x|\bth)$ and $F(x|\bth)$ respectively denote the pdf and cdf of a population distribution. Given $\bth$ and a real number $p\in (0,1)$, the $p$-th quantile $\xi_p$ of the distribution can be calculated as
$$
   \xi_p = \sol\, \Big\{ g(x|\bth) \teq p - F(x|\bth) = 0, \; \forall\; x \in \bbX \Big\},
$$
where we need to solve an integral equation. In this subsection, we will demonstrate that the US algorithm can be employed to iteratively find the $p$-th quantile $\xi_p$.

Assume that $g'(x|\bth) \teq \rd g(x|\bth)/ \rd x $ is bounded by some function $b(x|\bth)$; i.e.,
\begin{eqnarray} \label{uclaone4.1}
   g'(x|\bth) = - f(x|\bth) \ge b(x|\bth), \quad \forall\; x \in \bbX.
\end{eqnarray}
Then, from (\ref{uclaone2.6}), the US iteration is to solve
\begin{eqnarray} \label{uclaone4.2}
   x^{(t+1)} = \sol \left\{U(x| \xt, \bth) \sr{(\ref{uclaone3.3})} g(\xt |\bth) + \int_{\xt}^x b(z|\bth) \rd z =0, \; \forall \; x, \xt \in \bbX \right\}.
\end{eqnarray}
In this subsection, we will consider two cases: (i) The mode $x_{\rm mod}$ of the density $f(x|\bth)$ exists; (ii) The mode $x_{\rm mod}$ of $f(x|\bth)$ does not exist.

\subsubsection{The mode of the density exists} 

Let $x_{\rm mod}$ denote the mode of the density $f(x|\bth)$, then we have $- f(x|\bth) \ge - f(x_{\rm mod}|\bth)$ for all $x\in \bbX$. From
(\ref{uclaone4.2}), the US iteration is to solve
\begin{eqnarray} \label{uclaone4.3}
   x^{(t+1)} &=& \sol \Big\{ g(\xt | \bth) - f(x_{\rm mod}|\bth) (x - \xt) =0, \; \forall \; x, \xt \in \bbX \Big\}  \non \\ [2mm]
   &=& \xt - \frac{F(\xt |\bth) - p}{f(x_{\rm mod}|\bth)}.
\end{eqnarray}

\vkL \noi \textbf{Example 3} \textsf{(Quantile of the normal distribution)}. To find the $p$-th quantile of the normal distribution with mean $\mu$ and variance $\si^2$, we denote its cdf and pdf by $\Phi(x|\mu, \si^2)$ and $\phi(x|\mu, \si^2) = [1/(\sqrt{2\pi} \si)]\exp[-(x-\mu)^2/(2\si^2)]$, respectively. Note that the mode of $\phi(x|\mu, \si^2)$ is given by $x_{\rm mod}=\mu$, then (\ref{uclaone4.3}) becomes
\begin{eqnarray} \label{uclaone4.4}
   x^{(t+1)} = \xt - \sqrt{2\pi} \si \Big[ \Phi(\xt|\mu, \si^2) - p \Big].
\end{eqnarray}

\vkL \noi \textbf{Example 4} \textsf{(Quantile of the skew normal distribution)}. To find the $p$-th quantile of the skew normal distribution (Azzalini 1985) with pdf
\begin{eqnarray} \label{uclaone4.5}
   f(x|\mu, \si^2, \al) &=& \frac{2}{\si\sqrt{2\pi}} \exp{\left[-\frac{(x-\mu)^2}{2\si^2}\right]} \Phi\left( \al \frac{x-\mu}{\si} \right),
\end{eqnarray}
where $x, \mu \in \bbR$, $\si^2>0$, $\al \ge 0$ and $\Phi(\cdot)$ is the cdf of $N(0, 1)$, in Appendix B, we proposed an MM algorithm to iteratively calculate the mode $x_{\rm mod}$ of the skew normal density as shown in (\ref{uclaoneB.1}). From (\ref{uclaone4.3}), the US iteration for calculating the $p$-th quantile $\xi_p$ is
$$
   x^{(t+1)} = \xt - \frac{F(\xt |\mu, \si^2, \al ) - p}{f(x_{\rm mod}|\mu, \si^2, \al)},
$$
where $F(x |\mu, \si^2, \al)$ is the cdf of the skew normal distribution.

\subsubsection{The mode of the density does not exist} 

Let $F(x |\al, \be)$ denote the cdf of the beta distribution $\mBeta(\al, \be)$ and its pdf is
\begin{eqnarray*}
   f(x|\al, \be) &=& \frac{1}{B(\al, \be)} x^{\al-1}(1-x)^{\be-1}, \quad x\in (0, 1),\; \al>0, \; \be >0.
\end{eqnarray*}
Given $\al\in (0,1), \be\in (0,1)$ and $p\in (0,1)$, suppose that we want to calculate the $p$-th quantile of the distribution as $\xi_p = \sol\, \big\{ g(x|\al, \be) \teq p - F(x|\al, \be) = 0, \; \forall\; x \in (0, 1) \big\}$. We know that if at least one of $\al$ and $\be$ is less than 1, the mode of the beta density $f(x|\al, \be)$ does not exist. Thus, the US iteration specified by (\ref{uclaone4.3}) can not be used for such cases.

We employ (\ref{uclaone4.1})--(\ref{uclaone4.2}) to calculate $\xi_p$. First, it is easy to verify the following inequality:
\begin{eqnarray}  \label{uclaone4.6}
   \frac{1}{x(1-x)} \le \frac{1}{2x^2} + \frac{1}{2(1-x)^2}, \quad \forall\; x\in (0,1),
\end{eqnarray}
where the equality holds iff $x=1/2$. Next,
\begin{eqnarray*}
   g'(x|\al, \be) &=& - f(x| \al, \be) = - f(x|\al+1, \be+1) \frac{B(\al+1, \be+1)}{B(\al, \be)} \cdot \frac{1}{x(1-x)} \\ [2mm]
   &\sge{(\ref{uclaone4.6})} & - f(x| \al+1, \be+1) \frac{B(\al+1, \be+1)}{B(\al, \be)} \cdot \left[ \frac{1}{2x^2}+\frac{1}{2(1-x)^2} \right]\\ [2mm]
   &\ge& - f(x_{\rm mod, 1} | \al+1, \be+1) \frac{B(\al+1, \be+1)}{2 B(\al, \be)} \left[ \frac{1}{x^2} + \frac{1}{(1-x)^2} \right] \\ [2mm]
   &\sr{(\ref{uclaone4.7})} & -\frac{a}{x^2} - \frac{a}{(1-x)^2} \;\teq\; b(x |\al, \be), \quad \forall \; x \in (0, 1),
\end{eqnarray*}
where $x_{\rm mod, 1} \teq \al/(\al+\be)$ denotes the mode of the beta density $f(x|\al+1, \be+1)$ and
\begin{eqnarray}  \label{uclaone4.7}
   a \teq f(x_{\rm mod, 1}| \al+1, \be+1)\frac{B(\al+1, \be+1)} {2B(\al, \be)}.
\end{eqnarray}
Finally, from (\ref{uclaone4.2}), the US iteration for calculating the $p$-th quantile $\xi_p$ is
\begin{eqnarray*}
   x^{(t+1)} &=& \sol \left\{g(\xt |\al, \be) + \int_{\xt}^x b(z|\al, \be) \rd z =0, \; \forall \; x, \xt \in (0,1) \right\} 
\end{eqnarray*}
\begin{eqnarray*}
   &=& \sol \left\{g(\xt |\al, \be) + \frac{a}{x} - \frac{a}{1-x} - \frac{a}{\xt} + \frac{a}{1-\xt}  =0, \; \forall \; x, \xt \in (0,1) \right\} \\ [2mm]
   &=& \sol \left\{a_1^{(t)}x^2 + a_2^{(t)} x + a =0, \; \forall \; x, \xt \in (0,1) \right\},
\end{eqnarray*}
where $a_1^{(t)} \teq F(\xt |\al, \be) -p + a/\xt - a/(1-\xt)$ and $a_2^{(t)} \teq -a_1^{(t)} -2 a$.

\subsection{MLE of $\tha$ in Yule--Simon distribution} 

The Yule--Simon distribution (Yule 1925, Garcia 2011) can be employed to model species among genera, words frequencies in texts, numbers of papers published by researchers, cities by population, income by size (Leisen \Et 2017). Let $\{X_i\}_{i=1}^n \iid \mbox{YS}(\tha)$ with pmf $\tha B(x, \tha +1)$, $x=1, 2, \ldots,\infty$, where $\tha\, (>0)$ is the shape parameter. The log-likelihood function is $\ell(\tha) = n\log (\tha) + n \log \Ga(\tha +1) - \sum_{i=1}^n \log \Ga(x_i+\tha+1) + \mbox{constant}$. Thus, the MLE of $\tha$ is the root of the following equation:
\begin{eqnarray*}
   0 &=& \ell'(\tha) = \frac{n}{\tha}  + \sum_{i=1}^n [\psi(\tha+1)-  \psi(x_i+\tha+1)] \\ [2mm]
   &\sr{(\ref{uclaone4.10})}& \frac{n}{\tha}  + \sum_{i=1}^n  \left(  - \sum_{m=0}^{\infty} \frac{1}{m+\tha+1} +  \sum_{m=0}^{\infty} \frac{1}{m +x_i +\tha +1} \right)   \qquad [\mbox{Let } m+x_i=k]
   \\ [2mm]
   &=& \frac{n}{\tha}  + \sum_{i=1}^n  \left(  - \sum_{m=0}^{\infty} \frac{1}{m+\tha+1} +  \sum_{k=x_i}^{\infty} \frac{1}{k +\tha +1} \right) \non \\ [2mm]
   &=& \frac{n}{\tha}  + \sum_{i=1}^n  \left( -\sum_{m=0}^{x_i-1} \frac{1}{m +\tha +1} \right) \;\teq\;  g(\tha).
\end{eqnarray*}
Since
\begin{eqnarray*}
   g'(\tha) &=& - \frac{n}{\tha^2}  + \sum_{i=1}^n \left[ \sum_{m=0}^{x_i-1} \frac{1}{(m+\tha+1)^2}  \right] \\ [2mm]
   &=& - \frac{n}{\tha^2}  + \sum_{i=1}^n \left[ \frac{1}{(\tha+1)^2} + I(x_i\ge 2) \sum_{m=1}^{x_i-1} \frac{1}{(m+\tha+1)^2}\right]
   \ge - \frac{n}{\tha^2}  + \frac{ n}{(\tha+1)^2} \teq  b(\tha),
\end{eqnarray*}
the US iteration for calculating the MLE $\hth$ is
\begin{eqnarray}
   \thtI &\sr{(\ref{uclaone3.3})} & \sol \left\{ g(\tht) + \int_{\tht}^{\tha} b(z) \rd z = 0, \; \forall \; \tha, \tht >0 \right\}  \non \\ [2mm]
   &=& \sol \left\{ g(\tht) + \frac{n}{\tha}  - \frac{ n}{\tha+1} - \frac{n}{\tht} + \frac{ n}{\tht+1} = 0, \; \forall \; \tha, \tht >0  \right\} \non \\ [2mm]
   &=& \sol \left\{ a_9^{(t)} \tha^2 + a_9^{(t)}\tha +n =0, \; \forall \; \tha, \tht >0  \right\}, \label{uclaone4.15}
\end{eqnarray}
where $a_9^{(t)} \teq g(\tht) - n/\tht + n/(\tht+1)$.

\subsection{MLE of $\tha$ in generalized Poisson distribution}

Let $\{X_i\}_{i=1}^n\iid \mbox{GP}(\la, \tha)$ with pmf $\la(\la + \tha x)^{x-1} \e^{-\la -\tha x}/x!$ for $x=0, 1, \ldots,\infty$, where $\la\, (>0)$ and $\max (-1,-\la / r)<\tha \le 1$ are two parameters and $r\ (\ge 4)$ is the largest positive integer for which $\la+\tha r>0$ when $\tha<0$ (Consul \& Jain 1973; Consul \& Famoye 1992). The log-likelihood is $\ell(\la,\tha) = \sum_{i=1}^n  \left[ \log \la + (x_i-1) \log(\la+\tha x_i) -\la-\tha x_i \right] + \mbox{constant}$. There is an MM algorithm for calculating the MLE of $\tha$  for the case of $\tha \in [0, 1]$, while, to the best of our knowledge, there is no efficient algorithm to calculate $\hth$ for the case of $\tha<0$. In this subsection, we only consider the case of $\tha<0$. Let
\begin{eqnarray}
   0 &=& \frac{\pa \ell(\la,\tha)}{\pa \la} = \frac{n}{\la} + \sum_{i=1}^n \frac{x_i-1}{\la+\tha x_i} -n  \qand \label{uclaone4.16} \\ [2mm]
   0 &=& \frac{\pa \ell(\la,\tha)}{\pa \tha} =  - n\Bx + \sum_{i=1}^n \frac{(x_i-1)x_i} {\la+\tha x_i} \;\teq\; g(\tha|\la). \label{uclaone4.17}
\end{eqnarray}
On the one hand, by multiplying $\la$ and $\tha$ to the two partial differential equations (\ref{uclaone4.16}) and (\ref{uclaone4.17}), respectively, and then adding them together, we obtain $\la = (1-\tha)\Bx$. That is, given $\tha$, the MLE of $\la$ has an explicit expression.

On the other hand, given $\la$, we define
$$
   u_{1i}(\tha) = \frac{(x_i-1)x_i} {\la+\tha x_{(n)}} \qand
   u_{2i}(\tha) = \frac{\la+\tha x_{(n)}}{\la+\tha x_i}, \quad i=1, \ldots, n,
$$
where $x_{(n)} = \max(x_1, \ldots, x_n)$. Note that $u_{2i}'(\tha)=\la(x_{(n)} - x_i)/(\la+\tha x_i)^2 \ge 0$, then $u_{2i}(\tha)$ is increasing with respect to $\tha$. By applying the fixed--block method, given $\la$, we can construct a $U$-function for $g(\tha|\la)$ in (\ref{uclaone4.17}) as
\begin{eqnarray*}
   g(\tha|\la) =  - n\Bx + \sum_{i=1}^n u_{1i}(\tha)   u_{2i}(\tha)
   \sge{\sgn(\tha - \tht)}  - n\Bx + \sum_{i=1}^n u_{1i}(\tha)  u_{2i}(\tht) \;\teq\;  U(\tha|\tht, \la).
\end{eqnarray*}
Thus, the US iteration for calculating the MLE $\hth$ is
$$
   \thtI \sr{(\ref{uclaone2.6})} \sol \Big\{ U(\tha|\tht, \la) = 0, \; \forall \, \tha, \tht <0 \Big\}  = \frac{\la+\tht x_{(n)}} {n\Bx x_{(n)}} \sum_{i=1}^n  \frac{(x_i-1)x_i}{\la+\tht x_i} -\frac{\la}{x_{(n)}}.
$$

\section{$\!\!\!\!\!\!\!$. Analysis of the convergence rate} 

Let $\ths$ be the unique root of the equation $g(\tha)=0$ for all $\tha \in \Tha$. To obtain the rate of convergence for the US iteration $\thtI \teq h(\tht)$, we consider the first--order Taylor expansion of $h(\tht)$ around $\ths$: $h(\tht) = h(\ths) +(\tht-\ths) h'(\ths) + 0.5(\tht-\ths)^2 h''(\tth)$,
where $\tth$ is a point between $\tht$ and $\ths$. The rate of convergence for the US algorithm is defined by
\begin{eqnarray*}
   \lim_{t\to \infty} \frac{|\thtI - \ths|}{|\tht - \ths|} &=& \lim_{t\to \infty} \frac{|h(\tht) - h(\ths)|} {|\tht - \ths|} \\ [2mm]
   &=&  \lim_{t\to \infty} |h'(\ths) + 0.5(\tht-\ths) h''(\tilde{\tha})| = |h'(\ths)|.
\end{eqnarray*}
If $|h'(\ths)| \in (0, 1)$ or $|h'(\ths)| = 0$, then the US algorithm is said to have a linear rate of convergence or a super-linear rate of convergence.

In this section, we investigate the convergence rate of the US algorithm for the four methods of constructing $U$-functions presented in Section 3. Theorem 2 below states that the US algorithm based on the FLB function method holds a linear convergence rate, just like the EM and MM algorithms with monotonic convergence and linear convergence rate in maximization of an objective function. Theorems 3--4 below show that the US algorithm holds a quadratic and even a cubic convergence rates, corresponding to the SLUB constants method and the TLB constant method. Theorem 5 indicates that the US algorithm based on the fixed--block method possesses a linear convergence rate.

\vkL \noi \textbf{Theorem 2} \textsf{(Linear convergence rate for the FLB function method)}. Let the FLB function $b(\tha)$ be defined by (\ref{uclaone3.2}), and the $U$-function $U(\tha|\tht)$ be given by (\ref{uclaone3.3}). Then the US iteration
\begin{eqnarray} \label{uclaone5.1}
   \thtI = \sol\, \left\{U(\tha |\tht) = g(\tht) + \int_{\tht}^{\tha} b(z) \rd z = 0, \; \forall \; \tha, \tht \in \Tha \right\} \;\teq\; h(\tht)
\end{eqnarray}
holds a linear convergence rate given by
\begin{eqnarray} \label{uclaone5.2}
   \lim_{t\to \infty} \frac{|\thtI - \ths|}{|\tht - \ths|} = |h'(\ths)| = \left| 1- \frac{g'(\ths)}{b(\ths)} \right| = 1- \frac{g'(\ths)}{b(\ths)} \in [0, 1).
\end{eqnarray}
\hfill $\|$

\vkl \noi \textbf{Theorem 3} \textsf{(Quadratic convergence rate for the SLUB constants method)}. Let two SLUB constants $\{b_{21}, b_{22}\}$ be defined by (\ref{uclaone3.6}), and the $U$-function $U(\tha|\tht)$ be given by (\ref{uclaone3.9}). Then the US iteration
\begin{eqnarray} \label{uclaone5.3}
   \thtI = \sol \, \{U(\tha|\tht) =0, \quad \forall \; \tha, \tht \in \Tha \}  \teq  h(\tht)
\end{eqnarray}
has a quadratic convergence rate given by
\begin{eqnarray} \label{uclaone5.4}
   \lim_{t\to \infty} \frac{|\thtI - \ths|}{|\tht - \ths|^2} = \frac{|h''(\ths)|}{2} = \left|\frac{b_{22} - g''(\tha^*)}{2 g'(\tha^*)} \right| \in (0,\infty).
\end{eqnarray}
\hfill $\|$

\vkl \noi \textbf{Theorem 4} \textsf{(Cubic convergence rate for the TLB constant method)}. Let the TLB constant $b_3$ be defined by (\ref{uclaone3.14}), and the $U$-function $U(\tha|\tht)$ be given by (\ref{uclaone3.15}). Then the US iteration
\begin{eqnarray} \label{uclaone5.5}
   \thtI = \sol \, \{U(\tha|\tht) =0, \quad \forall \; \tha, \tht \in \Tha \} \teq  h(\tht)
\end{eqnarray}
has a cubic convergence rate given by
\begin{eqnarray} \label{uclaone5.6}
   \lim_{t\to \infty} \frac{|\thtI - \ths|}{|\tht - \ths|^3} = \frac{|h'''(\ths)|}{6} = \bigg|\frac{b_3 - g'''(\ths)}{6 g'(\ths)}\bigg| \in (0, \infty).
\end{eqnarray}
\hfill $\|$

\vkl \noi \textbf{Theorem 5} \textsf{(Linear convergence rate for the fixed--block method)}. Let the conditions in (\ref{uclaone3.16}) be satisfied, and the $U$-function $U_{1\cdots r}(\tha|\tht)$ given by (\ref{uclaone3.17}) satisfy the lipschitz condition $|U_{1\cdots r}(\tha_1|\tht)-U_{1\cdots r}(\tha_2|\tht)|\le L |\tha_1 - \tha_2|$ for any $\tha_1$ and $\tha_2$ located between $\tha^{(0)}$ and $\ths$.
Then the US iteration
\begin{eqnarray} \label{uclaone5.7}
   \thtI = \sol\, \{U_{1\cdots r}(\tha|\tht) = 0, \quad \forall \; \tha, \tht \in \Tha \} \teq h(\tht)
\end{eqnarray}
holds a linear convergence rate given by
\begin{eqnarray*}
   \lim_{t\to \infty} \frac{|\thtI - \ths|}{|\tht - \ths|} = |h'(\ths)| = \left| 1- \frac{g'(\ths)}{U_{1\cdots r}'(\ths|\ths)} \right| =  1- \frac{g'(\ths)}{U_{1\cdots r}'(\ths|\ths)} \in [0, 1),
\end{eqnarray*}
where
\begin{eqnarray*}
   U_{1\cdots r}'(\ths|\ths) & \teq & \frac{\rd U_{1\cdots r}(\tha|\ths)}{\rd \tha}\bigg|_{\tha=\ths} \\ [2mm]
   &=& \sum_{j=r+1}^m  \frac{\pa G \big( u_1(\ths), \ldots, u_r(\ths), u_{r+1}(\tha), \ldots, u_m(\tha) \big)} {\pa u_j(\tha)} u_j'(\tha)\bigg|_{\tha=\ths}
\end{eqnarray*}
is given by (\ref{uclaoneA.16}).  \hfill $\|$

\section{$\!\!\!\!\!\!\!$. Fast US algorithms with weakly stable convergence} 

Theorems 2 and 5 in Section 5 stated that US algorithms based on both the FLB function method and the fixed--block method hold a linear convergence rate only. To accelerate the convergence speed of such US algorithms, in this section we propose an acceleration technique with weakly stable convergence for US algorithms based on the FLB function method and the fixed--block method. For convenience, we define the following notations:
$$
   U'(\tht|\tht) \teq \frac{\rd U(\tha|\tht)}{\rd \tha}\bigg|_{\tha=\tht} \qand g'(\tht) \teq \frac{\rd g(\tha)}{\rd \tha}\bigg|_{\tha=\tht}.
$$

\vkL \noi \textbf{Theorem 6} \textsf{(A fast US algorithm with weakly stable convergence)}. Let $\ths$ be the unique root of the equation (\ref{uclaone1.1}), and $U(\tha|\tht)$ be specified by (\ref{uclaone3.3}) or (\ref{uclaone3.17}). A fast US iteration is defined as
\begin{eqnarray} \label{uclaone6.1}
   \thtI = \tht  + s(\tht) (\tth^{(t+1)} - \tht) = [1-s(\tht)] \tht + s(\tht)\tth^{(t+1)},
\end{eqnarray}
where $s(\tht)$ is said to be  the step length, $\tth^{(t+1)} - \tht$ is called a direction,
\begin{eqnarray}
   s(\tht) & \teq & \min\left\{ \frac{U'(\tht|\tht)}{g'(\tht)} I(g'(\tht)<0) + I(g'(\tht)\ge 0),\;  2\right\} \qand  \label{uclaone6.2} \\ [2mm]
   \tth^{(t+1)} &\teq& \sol \left\{U(\tha|\tht) =0, \; \forall\, \tha,\tht \in \Theta  \right\}.  \label{uclaone6.3}
\end{eqnarray}
Thus,
\begin{namelist}{01234}
\item[\hspace*{-0.03cm} (i)] For any initial value $\tha^{(0)}$ in $\Tha$, the sequence $\{\thtI\}_{t=0}^\infty$ specified by (\ref{uclaone6.1}) weakly converges to $\ths$. In other words, the fast US algorithm does not depend on any initial values in $\Tha$.

\item[\hspace*{-0.03cm} (ii)] If $b(\ths)/g'(\ths) > 2$, then the convergence speed ($2g'(\ths)/b(\ths)$) of the fast US algorithm is twice of the convergence speed ($g'(\ths)/b(\ths)$) of the US algorithm based on the FLB function method. If $1 \le b(\ths)/g'(\ths) \le 2$, then the fast US algorithm has a super--linear convergence rate comparing with a linear convergence rate for the US algorithm based on the FLB function method.

\item[\hspace*{-0.03cm} (iii)] If $U'(\ths |\ths)/g'(\ths) > 2$, then the convergence speed ($2g'(\ths)/U'(\ths |\ths)$) of the fast US algorithm is twice of the convergence speed ($g'(\ths)/U'(\ths |\ths)$) of the US algorithm based on the fixed--block method. If $1 \le U'(\ths |\ths)/g'(\ths) \le 2$, then the fast US algorithm has a super--linear convergence rate comparing with a linear convergence rate for the US algorithm based on the fixed--block method. \hfill $\|$
\end{namelist}

\vkL \noi \textbf{Remark 4} \textsf{(A comment on Theorem 6)}. Obviously, for the SLUB constants and TLB constant methods, both the $U$-functions in (\ref{uclaone3.9}) and (\ref{uclaone3.15}) satisfy $U'(\tht|\tht)/g'(\tht)=1$, so that $\thtI = \tth^{(t+1)}$. Therefore, our acceleration technique is useless for the US algorithm based on both the SLUB constants and TLB constant methods. \hfill $\|$

\section{$\!\!\!\!\!\!\!$. Numerical experiments} 

To compare the proposed US algorithm with Newton's method and bisection method by using numerical experiments, we consider two root-finding problems: Roots of a high-order polynomial equation and the $p$-th quantile of a normal distribution. In Section 7.3, we illustrate that the US algorithm can be applied to solve an equation with multiple roots. In Section 7.4, we investigate the computational efficiency of the proposed fast US algorithm. For further experimentation, our R code can be freely downloaded from the website https: \url{https://github.com/Xunjian-Li/US-algorithm.git}

\subsection{Root of a high-order polynomial equation} 

Assume that there exists a unique real root for the following polynomial equation in $(0, a)$:
\begin{eqnarray}  \label{uclaone7.1}
   0 = g(\tha) \teq a_3 \tha^m + a_2 \tha^2+ a_1 \tha+a_0 \teq g_1(\tha) + a_2 \tha^2+ a_1 \tha+a_0, \quad \tha\in (0, a),
\end{eqnarray}
where $m\; (\ge 3)$ is a real number (unnecessary positive integer) and $a_0>0$. Based on
\begin{eqnarray*}
   g_1''(\tha)  &=& a_3 m (m-1)\tha^{m-2} \; \in\;  (a_3 m (m-1)a^{m-2},0), \quad \mif a_3 < 0,     \\ [2mm]
   g_1'''(\tha) &=& a_3 m (m-1)(m-2) \tha^{m-3} \; >\;  0, \quad\quad\quad\quad\quad\quad \mif a_3 > 0,
\end{eqnarray*}
and Property 1(c), we can construct the following two quadratic $U$-functions:
\begin{eqnarray}
   U_2(\tha|\tht) &\sr{(\ref{uclaone3.9})} & a_3 [\tht]^m + a_3 m [\tht]^{m-1} (\tha-\tht) + 0.5 b_2(\tha|\tht) (\tha-\tht)^2 \non \\ [2mm]
   & & + \; a_2 \tha^2+ a_1 \tha+a_0 \non \\ [2mm]
   &=& \Big[ 0.5 b_2(\tha|\tht) + a_2 \Big] \tha^2 + \Big\{ a_3 m [\tht]^{m-1} - b_2(\tha|\tht) \tht + a_1 \Big\} \tha \non \\ [2mm]
   & & +\; a_3(1-m) [\tht]^m + 0.5 b_2(\tha|\tht) [\tht]^2 + a_0 \label{uclaone7.2}  \qand \\ [2mm]
   U_3(\tha|\tht) &\sr{(\ref{uclaone3.15})} & a_3 [\tht]^m + a_3 m [\tht]^{m-1} (\tha-\tht) + 0.5 a_3 m(m-1) [\tht]^{m-2}(\tha-\tht)^2   \non \\ [2mm]
   & & +\; a_2 \tha^2 + a_1 \tha+a_0 \non \\  [2mm]
   &=& \Big \{ 0.5 a_3 m(m-1) [\tht]^{m-2} + a_2 \Big\} \tha^2 + \Big\{ a_3m(2- m) [\tht]^{m-1} + a_1 \Big\} \tha \non \\ [2mm]
   & &+\; a_3(m-1) (0.5m -1) [\tht]^m + a_0, \label{uclaone7.3}
\end{eqnarray}
where $b_2(\tha| \tht) = 0 \cdot I(\tha \le \tht) + a_3 m (m-1)a^{m-2} \cdot I(\tha > \tht)$.

\vkL
\begin{table}[h]
  \baselineskip 0.10in
  {\bf Table 2.} $\;$ {\rm Comparisons of three algorithms in terms of the percentage of converged algorithms, the average number of iterations, and the total computation time in second}
  \vspace*{-0.4cm}
  \begin{center}
  \renewcommand{\arraystretch}{1.2} \tabcolsep 0.12in \doublerulesep 0.2pt
  \begin{threeparttable}
  \begin{tabular}{l|ccc|ccc} \hline\hline
    Coefficients & \multicolumn{3}{c|}{$\{a_0, a_1, a_2, a_3, m\} = \{1, -1, 1, -1, 3\}$}    & \multicolumn{3}{c} {$\{a_0, a_1, a_2, a_3, m\} = \{1, -1, -3, 1, 3\}$}  \\ \hline
    Algorithms   &$\mbox{US}_2$ & Newton     &Bisection  &$\mbox{US}_3$ &Newton      &Bisection    \\  \hline
    \textsf{Percentage}   & 100.0\%     & 100.0\%   & 100.0\%   & 100.0\%      & 100.0\%    & 100.0\%     \\
    \textsf{N}            & 7.0000      & 6.3672    & 41.000    & 7.0000       & 6.3630     & 40.000      \\
    \textsf{Time(s)}      & 1.6931      & 1.5313    & 10.255    & 1.6931      & 1.6206     & 10.184      \\  \hline \hline
  \end{tabular}
  {\small Note: $\mbox{US}_2$ denotes the US algorithm based on $U_2(\tha|\tht)$ in (\ref{uclaone7.2}); $\mbox{US}_3$ denotes the US algorithm based on $U_3(\tha|\tht)$ in (\ref{uclaone7.3});
  \textsf{Percentage} is the percentage of converged algorithms among 100,000 repetitions for different initial values randomly chosen from $U(0, 2)$; \textsf{N} is the average number of iterations for converged algorithms among 100,000 repetitions; \textsf{Time(s)} is the total computation time in second of calculating $\ths$ for 100,000 repetitions. }
  \end{threeparttable}
  \end{center}
\end{table}

For the purpose of numerical experiments, we set $\{a_0, a_1, a_2, a_3, m\} = \{1, -1, 1, -1, 3\}$ in (\ref{uclaone7.2}), $\{a_0, a_1, a_2, a_3, m\} = \{1, -1, -3, 1, 3\}$ in (\ref{uclaone7.3}), and $a=2$. We know that the unique root of (\ref{uclaone7.1}) in the interval (0, 2) is $\ths=0.4608111$ and $\ths=1$ for the two cases. We use three algorithms (i.e., the US algorithm, Newton's method and bisection method) to seek the root of (\ref{uclaone7.1}) by running each algorithm for 100,000 repetitions with initial values randomly chosen from the uniform distribution on the interval $(0, 2)$. The percentage of converged algorithms among 100,000 repetitions (denoted by \textsf{Percentage}), the average number of iterations for converged algorithms among 100,000 repetitions (denoted by \textsf{N}), and the total computation time in second of calculating $\ths$ for 100,000 repetitions [denoted by \textsf{Time(s)}] are recorded and reported in Table 2.

From Table 2, we can see that all three algorithms have 100\% convergence performance according to \textsf{Percentage}. The US algorithm and Newton's methods have similar performance in terms of \textsf{N} and \textsf{Time(s)}.  The bisection method performs worst among the three for the two cases.

\subsection{Calculation of the $p$-th quantile of normal distribution} 

In Example 3 of subsection 4.1.1, we have applied the US algorithm to find the $p$-th quantile of the normal distribution $N(\mu, \si^2)$; i.e., to iteratively calculate the solution to the integral equation $0 = g(x) \teq p - \Phi(x|\mu, \si^2)$ with US iteration given by (\ref{uclaone4.4}). Theorem 2 shows that the US iteration (\ref{uclaone4.4}) only has a linear rate of convergence. To speed the convergence, we first derive the second--derivative to the fourth--derivative of $g(x)$ as follows:
\begin{eqnarray*}
   g''(x)     &=& \frac{x-\mu}{\si^2} \phi(x|\mu, \si^2), \\ [2mm]
   g'''(x)    &=& \frac{1}{\si^2}\left[1 - \frac{(x-\mu)^2}{\si^2} \right] \phi(x|\mu, \si^2) \qand  \\ [2mm]
   g^{(4)}(x) &=& -\frac{x-\mu}{\si^4}\left[3 - \frac{(x-\mu)^2}{\si^2} \right]\phi(x|\mu, \si^2).
\end{eqnarray*}
Note that there are three roots $\{\mu, \mu \pm \sqrt{3} \si\}$ for the equation $0 = g^{(4)}(x)$, then $g'''(x)$ achieves its minimum $- 2 (\sqrt{2\pi} \si^3 \e^{3/2} )^{-1}$ when $x = \mu \pm \sqrt{3} \sigma$. Similarly, since $\mu\pm \sigma$ are two roots of the equation $0 = g'''(x)$, we know that $g''(x)$ achieves its minimum $-(\sqrt{2\pi} \si^2 \e^{1/2})^{-1}$ and maximum $(\sqrt{2\pi} \si^2 \e^{1/2})^{-1}$ when $x = \mu- \sigma$ and $\mu+ \sigma$, respectively. Thus, we obtain
\begin{eqnarray*}
   g''(x) & \in & \left[-(\sqrt{2\pi} \si^2 \e^{1/2} )^{-1},\; (\sqrt{2\pi} \si^2 \e^{1/2})^{-1}\right] \qand     \\ [2mm]
   g'''(x) & \ge & - 2 (\sqrt{2\pi} \si^3 \e^{3/2} )^{-1}, \quad \forall \; x\in \bbR,
\end{eqnarray*}
so that we can construct the following two $U$-functions:
\begin{eqnarray}
   U_2(x|\xt) & \sr{(\ref{uclaone3.9})} & g(\xt) + g'(\xt)(x- \xt) + 0.5 b_2(x|\xt) (x-\xt)^2   \qand   \label{eqnfive7.4} 
\end{eqnarray}
\begin{eqnarray}
   U_3(x|\xt) &\sr{(\ref{uclaone3.15})} & g(\xt) + g'(\xt)(x-\xt) + 0.5 g''(\xt) (x-\xt)^2  \non \\  [2mm]
   & & -\; (3\sqrt{2\pi} \si^3 \e^{3/2} )^{-1}(x-\xt)^3, \label{uclaone7.5}
\end{eqnarray}
where $b_2(x|\xt) = (\sqrt{2\pi} \si^2 \e^{1/2} )^{-1} \times [I(x \le \xt) -  I(x > \xt)]$.

Given $p=0.05$, $\mu=1$ and $\sigma=1$, Figure 2(a) respectively plots $g(x) = p-\Phi(x|\mu, \sigma^2)$ with black solid curve and $U_2(x|\xt)$ in (\ref{eqnfive7.4}) with the red dashed curve, which decreases monotonically rather than being a quadratic curve. Given $\xt = -1.5$, by taking $b_2(x|\xt) = -(\sqrt{2\pi} \si^2 \e^{1/2} )^{-1}$, we have $g(\xt)>0$ and $x^{(t+1)}$ would be the unique root of the $U$-equation $U_2(x|\xt )=0$, resulting in $\xt < x^{(t+1)} \le x^*$. On the other hand, if given another $\xt$ satisfying $g(\xt )<0$, by taking $b_2(x|\xt) = (\sqrt{2\pi} \si^2 \e^{1/2} )^{-1}$, then $x^{(t+1)}$ would be the unique root of the $U$-equation $U_2(x|\xt )=0$, resulting in $x^* \le x^{(t+1)} < \xt$.

\begin{figure}
  \begin{center}
     $\,$\vskip -3.0cm     \includegraphics[scale=0.6]{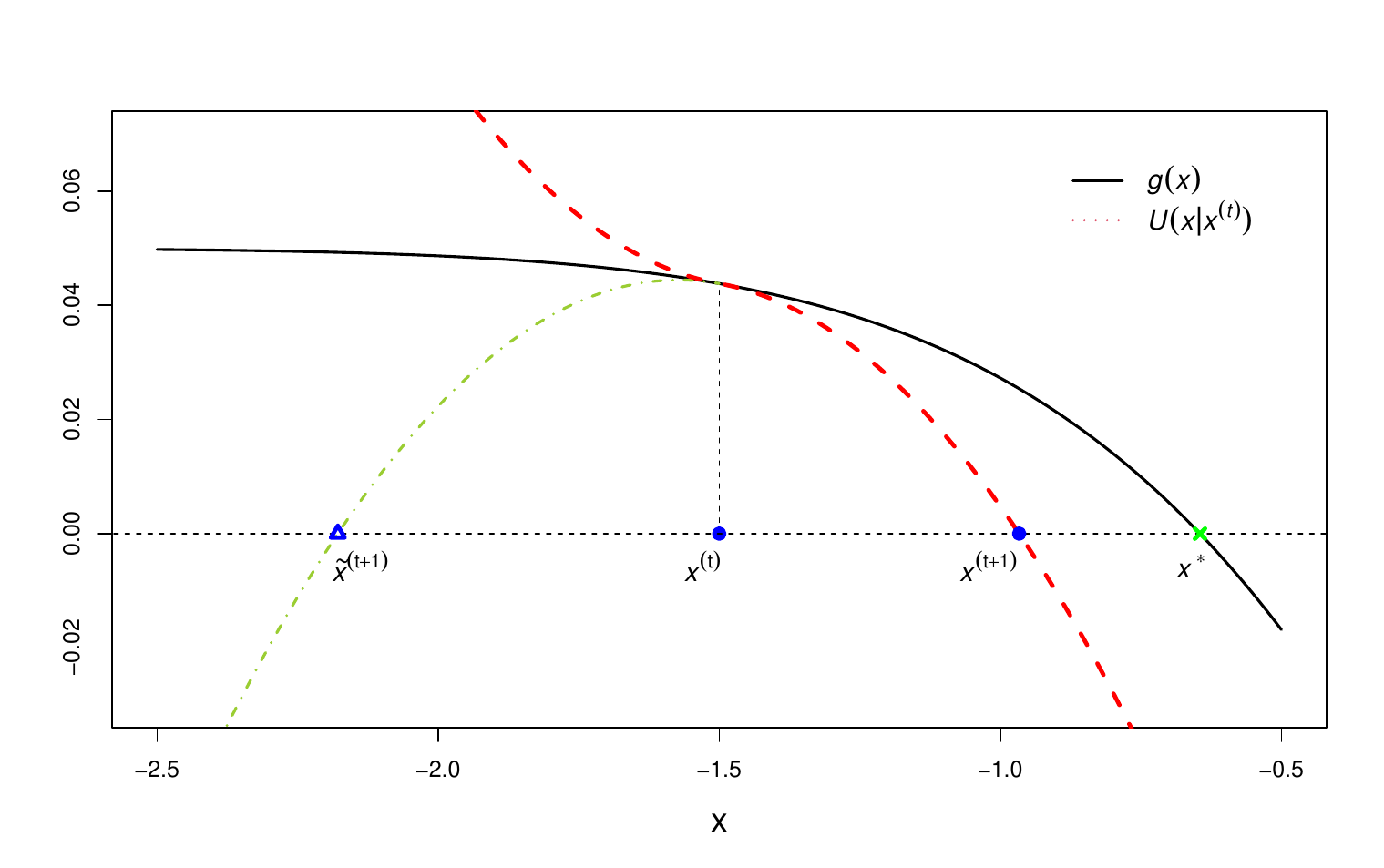}
     \vskip -0.0cm \baselineskip 0.10in {\small  (a)}
     \vskip 0.1cm          \includegraphics[scale=0.6]{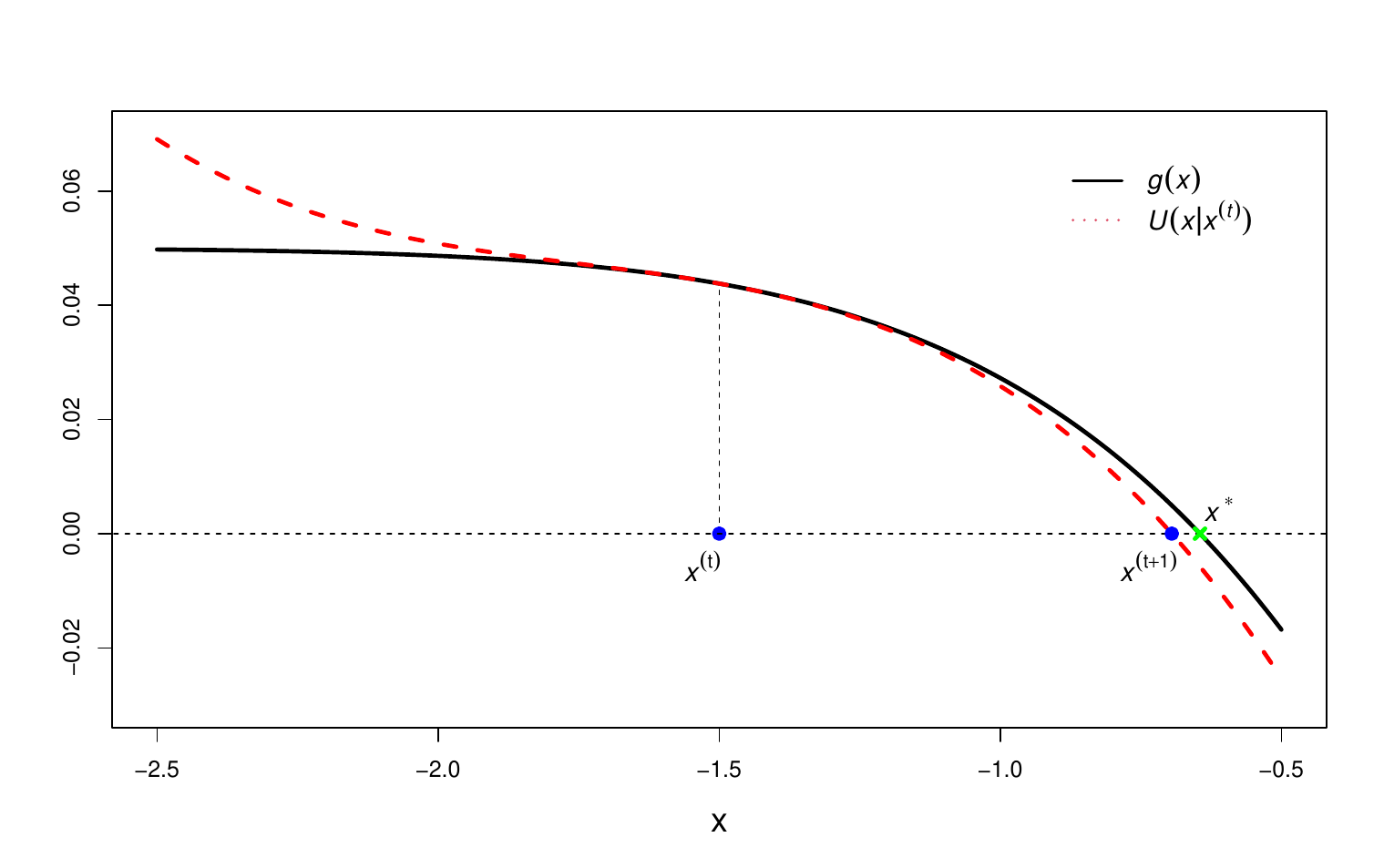}
     \vskip -0.0cm \baselineskip 0.10in {\small  (b)}
  \end{center}
  \vskip -0.0cm \baselineskip 0.10in {\small {\bf Figure 2.} $\;$ Plots of $g(x) = 0.05 - \Phi(x|\mu, \sigma^2)$ with $\mu=1$ and $\sigma=1$. The $0.05$-th quantile of the normal distribution is $x^*=-0.6448536$ and the $t$-th iteration $\xt =-1.5$. (a) From $g(\xt )>0$, the $(t+1)$-th iteration $x^{(t+1)}$ would be the unique root of the $U$-equation $0 = U_2(x|\xt )$ in (\ref{eqnfive7.4}) by taking $b_2(x|\xt) = -(\sqrt{2\pi} \si^2 \e^{1/2} )^{-1}$, resulting in $\xt < x^{(t+1)} \le x^*$. (b) From $g(\xt )>0$, the $x^{(t+1)}$ would be the unique root of the $U$-equation $0 = U_3(x|\xt )$ in (\ref{uclaone7.5}), resulting in $\xt < x^{(t+1)} \le x^*$. }
\end{figure}

\vkL
\begin{table}[h]
  \baselineskip 0.10in
  {\bf Table 3.} $\;$ {\rm Comparisons of four algorithms in terms of the percentage of converged algorithms, the average number of iterations, and the total computation time in second }
  \vspace*{-0.4cm}
  \begin{center}
  \renewcommand{\arraystretch}{1.2} \tabcolsep 0.07in \doublerulesep 0.2pt
  \begin{threeparttable}
  \begin{tabular}{l|cccc|cccc} \hline\hline
    Parameters   &\multicolumn{4}{c|}{$\{p, \mu\} = \{0.01, -2\}$}
                 & \multicolumn{4}{c}{$\{p, \mu\} = \{0.01, 2\}$}  \\ \hline
    Algorithms   &$\mbox{US}_2$&$\mbox{US}_3$&Newton      &Bisection    & $\mbox{US}_2$         & $\mbox{US}_3$          &Newton      &Bisection   \\  \hline
    \textsf{Percentage}   & 100.0\%   & 100.0\%  & {\bf 39.83 \%} & 100.0\%   & 100.0\%      & 100.0\%  & {\bf 55.13\%}  & 100.0\%   \\
    \textsf{N}            & 10.542    & 5.0683   & 6.6596        & 33.000    & 10.289       & 4.7208   & 6.3622         & 36.000  \\
    \textsf{Time(s)}      & 6.8351    & 6.6289   & 1.9440        & 9.9709    & 6.7457       & 6.0023   & 2.1058         & 10.714   \\  \hline\hline
    Parameters   &\multicolumn{4}{c|}{$\{p, \mu\} = \{0.9, -2\}$}  & \multicolumn{4}{c}{$\{p, \mu\} = \{0.9, 2\}$}  \\ \hline
    Algorithms   &$\mbox{US}_2$&$\mbox{US}_3$&Newton      &Bisection    & $\mbox{US}_2$          & $\mbox{US}_3$          &Newton      &Bisection   \\  \hline
    \textsf{Percentage}   & 100.0\%   & 100.0\%  & {\bf 42.76\%}   & 100.0\%   & 100.0\%   & 100.0\%  & {\bf 38.51\%}  & 100.0\%   \\
    \textsf{N}            & 5.9950    & 3.8007   & 5.0339    & 39.000    & 6.7522   & 4.2315    & 4.9217   & 40.000 \\
    \textsf{Time(s)}      & 4.1958    & 5.0533   & 1.7797    & 11.934    & 4.6499   & 5.5095    & 1.7026   & 11.849  \\  \hline \hline
  \end{tabular}
  {\small Note: $\mbox{US}_2$ denotes the US algorithm based on $U_2(\tha|\tht)$ in (\ref{eqnfive7.4}); $\mbox{US}_3$ denotes the US algorithm based on $U_3(\tha|\tht)$ in (\ref{uclaone7.5});
  \textsf{Percentage} is the percentage of converged algorithms among 100,000 repetitions for different initial values randomly chosen from $U(-4,4)$; \textsf{N} is the average number of iterations for converged algorithms among 100,000 repetitions; \textsf{Time(s)} is the total computation time in second of calculating $x^*$ for 100,000 repetitions. }
  \end{threeparttable}
  \end{center}
\end{table}

In numerical experiments, we set $p=0.01, \, 0.9$, $\mu=-2, \, 2$ and $\si=1$. For each combination of the four cases of $\{p, \mu \}$, we use four algorithms (i.e., two US algorithms based on (\ref{eqnfive7.4}) and (\ref{uclaone7.5}), Newton's method and bisection method) to seek the root of (\ref{uclaone7.1}) by running each algorithm for 100,000 repetitions with initial values randomly chosen from $U(-4, 4)$. The percentage of converged algorithms among 100,000 repetitions (denoted by \textsf{Percentage}), the average number of iterations for converged algorithms among 100,000 repetitions (denoted by \textsf{N}), and the total computation time in second of calculating $x^*$ for 100,000 repetitions [denoted by \textsf{Time(s)}] are recorded and displayed in Table 3.

In Table 3, we can see that Newton's method has the worst convergence performance according to \textsf{Percentage}. In terms of \textsf{N}, the $\mbox{US}_3$ performs the best. In terms of \textsf{Time(s)}, although Newton's method only takes about 2 seconds for 100,000 repetitions for Case 1, we can see that there are only about 40\% valid computations, hence we prefer $\mbox{US}_2$ and $\mbox{US}_3$.

\subsection{Solving an equation with multiple roots} 

Thanks for the strongly stable convergence, the US algorithm could be one of the powerful tools for solving an equation with multiple roots. For the purpose of illustration, we consider the following toy equation:
\begin{eqnarray} \label{uclaone7.6}
   0 = g(x) \teq -0.5x-2\sin(x) +1, \quad x \in \bbR.
\end{eqnarray}
Figure 3(a) shows the equation (\ref{uclaone7.6}) with three roots, denoted by $x_1^* < x_2^* < x_3^*$, respectively. First, based on the first--order derivative $g'(x) = -0.5 - 2\cos(x) \ge -2.5$, we can construct the first $U$-function for $g(x)$ at $x=\xt$ and obtain the corresponding US iterations as follows:
\begin{eqnarray} \label{uclaone7.7}
   U_{1,1}(x|\xt) & \sr{(\ref{uclaone3.5})} & g(\xt) -2.5 (x- \xt) \qand \non \\ [2mm]
   x^{(t+1)} &=& \xt + \frac{g(\xt)}{2.5}  = \xt + \frac{-0.5\xt -2 \sin(\xt) +1 } {2.5}.
\end{eqnarray}
Given an initial value $x^{(0)}\, (<x_1^*)$, say $x^{(0)}=0$, the US iterations $\{x^{(t+1)}\}_{t=0}^{\infty}$ defined by (\ref{uclaone7.7}) converged to the first root $x_1^* = 0.4090497$ as shown in Figure 3(b1).

\begin{figure}
  \begin{center}
     $\,$\vskip -3.0cm     \includegraphics[scale=0.6]{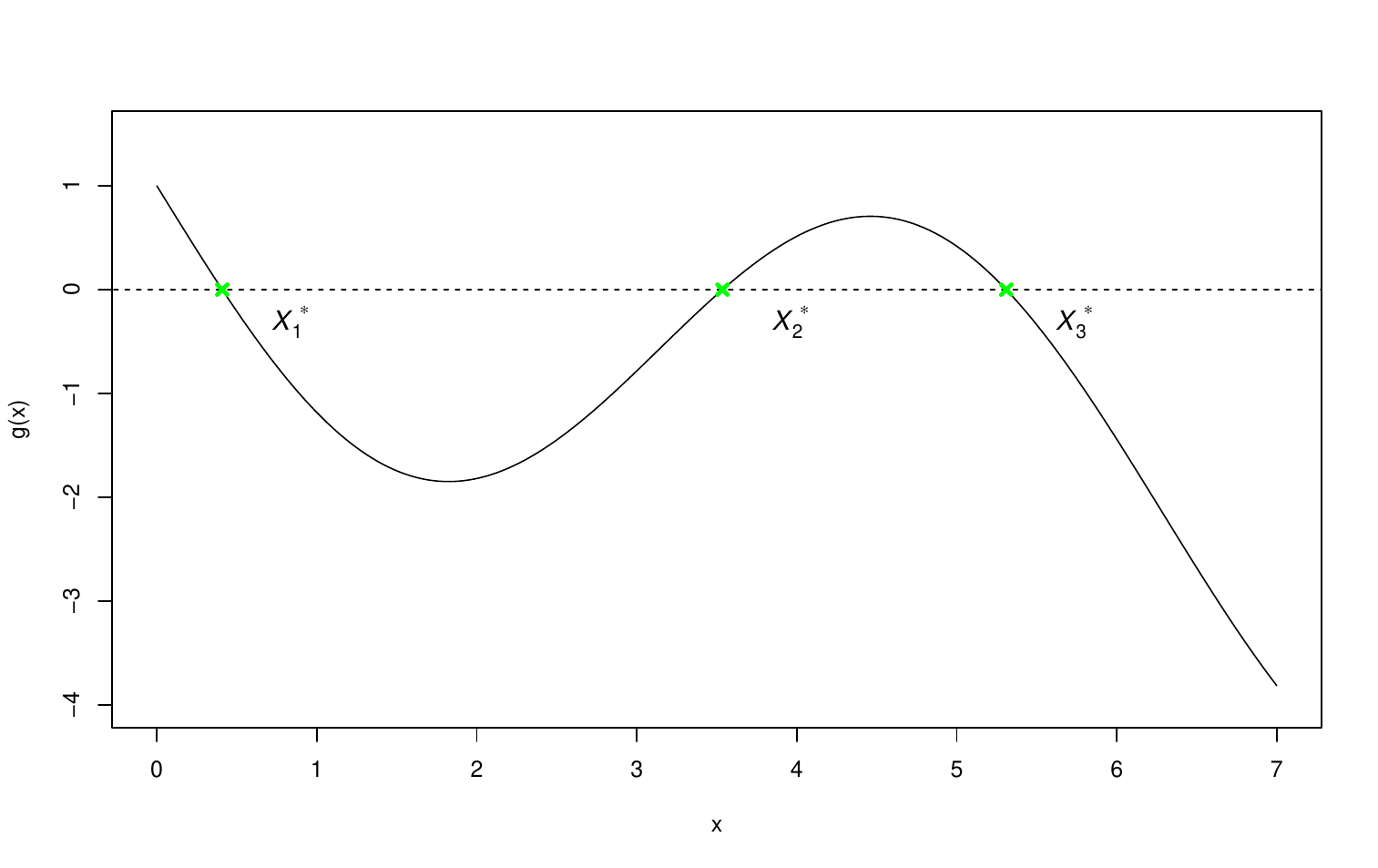}
     \vskip -0.0cm \baselineskip 0.10in {\small  (a)}
     \vskip 0.6cm          \includegraphics[scale=0.6]{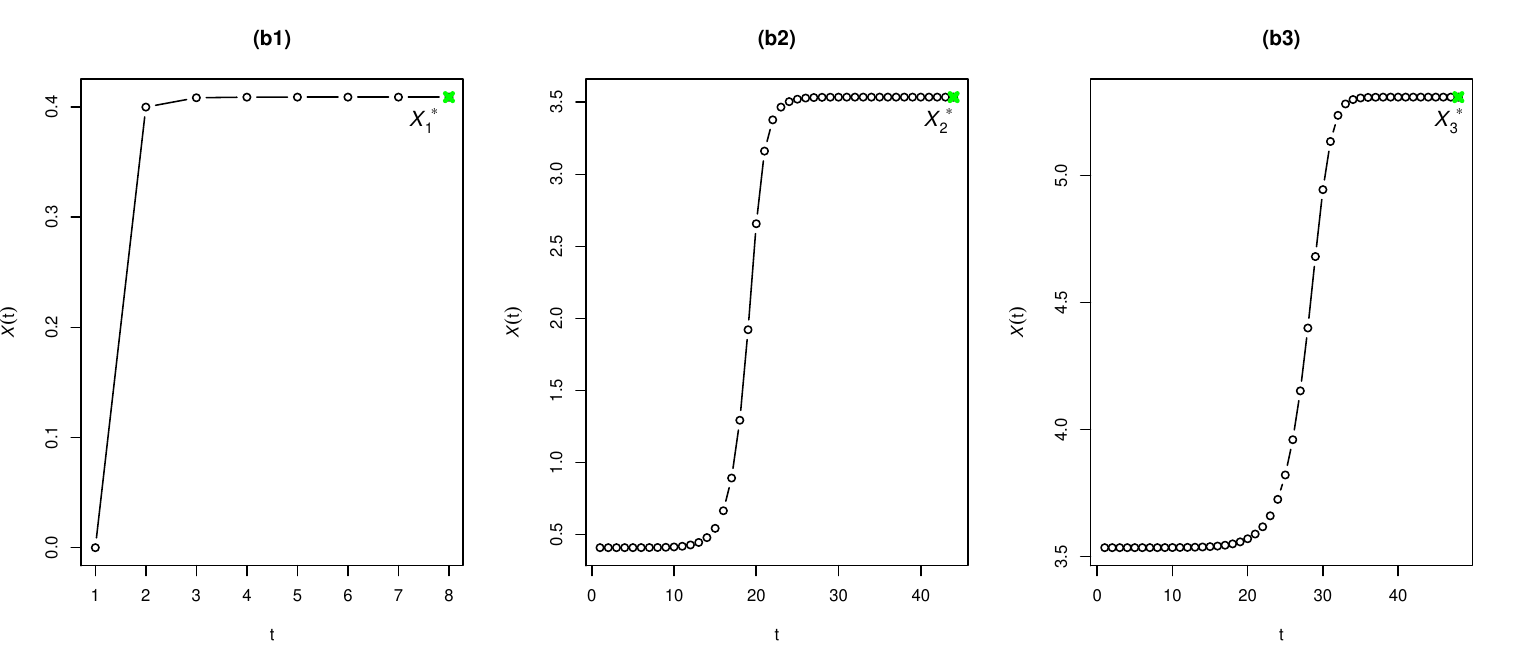}
     \vskip -0.0cm \baselineskip 0.10in {\small  (b)}
  \end{center}
  \vskip 0.5cm \baselineskip 0.10in {\small {\bf Figure 3.} $\;$ (a) Plots of $0 = g(x) = -0.5x-2\sin(x) +1$ with three roots $x_1^*=0.4090497$, $x_2^*=3.535612$ and $x_3^*=5.308993$. (b) With the initial values $x^{(0)}=0$, $x^{(0)}=x_1^*+10^{-6}$ and $x^{(0)}=x_2^*+10^{-6}$, the sequence $\{x^{(t+1)}\}_{t=0}^\infty$ generated respectively by (\ref{uclaone7.7}), (\ref{uclaone7.8}), (\ref{uclaone7.7}) strongly stably converges to $\{x_k^*\}_{k=1}^3$, satisfying $x^{(0)} < x^{(1)} < \cdots < \xt < \cdots \le x_k^*$ for $k=1,2,3$.}
\end{figure}

Second, once we obtained $x_1^*$, we define a new function
$$
   g_2(x) \teq -g(x) = 0.5x + 2\sin(x) - 1, \quad x \in \bbR.
$$
Based on the first--order derivative $g'_2(x) = 0.5 + 2\cos(x) \ge -1.5$, we can construct the second $U$-function for $g_2(x)$ at $x=\xt$ and obtain the corresponding US iterations as follows:
\begin{eqnarray} \label{uclaone7.8}
   U_{1,2}(x|\xt) & \sr{(\ref{uclaone3.5})} & g_2(\xt) -1.5 (x- \xt) \qand \non \\ [2mm]
   x^{(t+1)} &=& \xt + \frac{g_2(\xt)}{1.5}  = \xt + \frac{0.5\xt +2 \sin(\xt) +1 } {2.5}.
\end{eqnarray}
Set $x^{(0)} \teq x_1^* + \ve$ as the initial value such that $g_2(x^{(0)}) >0$, where $\ve$ is a small positive real number, say $10^{-6}$. Then the US iterations $\{x^{(t+1)}\}_{t=0}^{\infty}$ defined by (\ref{uclaone7.8}) converged to the second root $x_2^* = 3.535612$ as shown in Figure 3(b2).

Third, once we obtained $x_2^*$, we set $x^{(0)} \teq x_2^* + \ve$ as the initial value such that $g(x^{(0)}) >0$, then the US iterations $\{x^{(t+1)}\}_{t=0}^{\infty}$ defined by (\ref{uclaone7.7}) converged to the second root $x_3^* =  5.308993$ as shown in Figure 3(b3).

\subsection{Computational efficiency of the proposed fast US algorithm} 

To investigate the computational efficiency of the proposed fast US algorithm defined by (\ref{uclaone6.1}), we consider again the calculation of the MLE of $\tha$ in the Yule--Simon distribution (\S 4.2), which is equivalent to seeking the root of the following equation:
\begin{eqnarray} \label{uclaone7.9}
   0 = u_1(\tha) + u_2(\tha), \qwh u_1(\tha) \teq \sum_{i=1}^n  \left( -\sum_{m=0}^{x_i-1} \frac{1}{m +\tha +1} \right), \; u_2(\tha) \teq \frac{n}{\tha}.
\end{eqnarray}
Garcia (2011) provided a fixed--point algorithm to solve the root of (\ref{uclaone7.9}) by replacing $\tha$ in $u_1(\tha)$ with $\tht$, and obtained the following iteration:
\begin{eqnarray} \label{uclaone7.10}
   \thtI  = \frac{n}{-u_1(\tht)}.
\end{eqnarray}
Note that
$$
   \frac{\pa [u_1(\tha) + u_2(\tha)]}{\pa u_1(\tha)} u_1'(\tha) = 1 \times \sum_{i=1}^n \sum_{m=0}^{x_i-1} \frac{1}{(m +\tha +1)^2} > 0
$$
satisfying the condition in (\ref{uclaone3.16}), then the fixed--point algorithm (\ref{uclaone7.10}) is a special US algorithm based on the fixed--block method with $U$-function $U(\tha|\tht) = u_1(\tht) +n/\tha$. Since every US algorithm has a strongly stable convergence, of course, we can also construct a fast fixed--point algorithm defined by (\ref{uclaone6.1}) to accelerate its convergence speed.

\begin{table}[!h]
  \baselineskip 0.10in
  {\bf Table 4.} $\;$ {\rm Comparisons of five algorithms in terms of the percentage of converged algorithms, the average number of iterations, and the total computation time in second }
  \vspace*{-0.2cm}
  \begin{center}
  \renewcommand{\arraystretch}{1.2} \tabcolsep 0.17in \doublerulesep 0.2pt
  \begin{threeparttable}
  \begin{tabular}{l|ccccc} \hline\hline
    Parameters   &\multicolumn{5}{c}{$\tha = 0.5$}  \\ \hline
    Algorithms   &$\mbox{US}$&$\mbox{F--US}$&Newton&$\mbox{fixed--point}$        & $\mbox{F--fixed--point}$   \\  \hline
    \textsf{Percentage}   &100.00\% & 100.00\% & {\bf 71.71\%} & 100.00\% & 100.00\%                     \\
    \textsf{N}            & 8.570 & 5.229 & 10.205 & 12.276 & 6.252                     \\
    \textsf{Time(s)}      & 48.099 & 53.154 & {\bf 3676.3} & 68.552 & 64.276              \\ \hline\hline
    Parameters   &\multicolumn{5}{c}{$\tha = 1$}  \\ \hline
    Algorithms   &$\mbox{US}$&$\mbox{F--US}$&Newton&$\mbox{fixed--point}$        & $\mbox{F--fixed--point}$   \\  \hline
    \textsf{Percentage}   & 100.00\% & 100.00\% & {\bf 70.04\%} & 100.00\% & 100.00\%                     \\
    \textsf{N}            & 10.913 & 5.689 & 6.905 & 19.191 & 6.443                  \\
    \textsf{Time(s)}      & 39.716 & 35.527 & {\bf 1552.2} & 69.279 & 40.340        \\  \hline\hline
    Parameters   &\multicolumn{5}{c}{$\tha = 5$}  \\ \hline
    Algorithms   &$\mbox{US}$&$\mbox{F--US}$&Newton&$\mbox{fixed--point}$        & $\mbox{F--fixed--point}$   \\  \hline
    \textsf{Percentage}   & 100.00\%     & 100.00\%   & 100.00 \%      & 100.00\%       & 100.0\%                       \\
    \textsf{N}            & 20.170 & 5.465 & 6.771 & {\bf 83.400} & 39.005                  \\
    \textsf{Time(s)}      & 60.379 & 27.462 & 33.085 & {\bf 245.40} & {\bf 189.63}      \\  \hline \hline
    Parameters   &\multicolumn{5}{c}{$\tha = 10$}  \\ \hline
    Algorithms   &$\mbox{US}$&$\mbox{F--US}$&Newton&$\mbox{fixed--point}$        & $\mbox{F--fixed--point}$   \\  \hline
    \textsf{Percentage}   & 100.00\%     & 100.00\%   & 100.00 \%      & 100.00\%       & 100.00\%                       \\
    \textsf{N}            & 25.400 & 6.663 & 8.952 & {\bf 179.87} & {\bf 89.48}                \\
    \textsf{Time(s)}      & 73.708 & 32.203 & 42.005 & {\bf 515.50} & {\bf 418.23}              \\ \hline
    \hline
  \end{tabular}
  {\small Note: US denotes the US iteration in (\ref{uclaone4.15}); F-US denotes the fast US algorithm based on (\ref{uclaone4.15}) \& (\ref{uclaone6.1}); \textsf{Percentage} is the percentage of converged algorithms among 10,000 repetitions for different initial values randomly chosen from $U(1,5)$; \textsf{N} is the average number of iterations for converged algorithms among 10,000 repetitions; \textsf{Time(s)} is the total computation time in second of calculating $\hth$ for 10,000 repetitions. }
  \end{threeparttable}
  \end{center}
\end{table}

To compare the computational efficiency among the US algorithm in (\ref{uclaone4.15}), the fast US algorithm based on (\ref{uclaone4.15}) \& (\ref{uclaone6.1}), Newton's method, fixed--point algorithm (\ref{uclaone7.10}) and the fast fixed--point algorithm (\ref{uclaone7.10}) \& (\ref{uclaone6.1}), we choose the true values of $\tha$ to be 0.5, 1, 5, 10 and sample size $n=$ 400. We apply the five algorithms for each scenario to seek the root of (\ref{uclaone7.9}) by conducting each algorithm for 10,000 repetitions with initial values randomly chosen from the uniform distribution $U(1, 5)$. The percentage of converged algorithms among 10,000 repetitions (denoted by \textsf{Percentage}), the average number of iterations for converged algorithms among 10,000 repetitions (denoted by \textsf{N}), and the total computation time in second of calculating $\hth$ for 10,000 repetitions [denoted by \textsf{Time(s)}] are recorded and displayed in Table 4.

In Table 4, we can see that Newton's method has the worst convergence performance according to \textsf{Percentage}, \textsf{N} and \textsf{Time(s)} for $\tha=0.5$ and 1, and the fixed--point algorithm has the worst convergence performance according to \textsf{N} and \textsf{Time(s)} for $\tha=5$ and 10. The fast fixed--point algorithm converges faster than the fixed--point algorithm with a less average number of iterations for each algorithm and total computation time. Specifically, for small $\tha$, the fast fixed--point algorithm enjoys similarly convergence performance to the fast US algorithm, which possesses the best convergence performance among the five algorithms in all cases. For larger $\tha$, the fast fixed--point algorithm reduces the average number of iterations for each fixed--point algorithm by at least half. Except for the fast US algorithm, the US algorithm performs better than the rest algorithms for $\tha \ge 5$. So, it is safe to conclude that the fast US algorithm is of good computation efficiency with the least average number of iterations and total computation time in seconds.

\section{$\!\!\!\!\!\!\!$. Discussions} 

To solve the root of a non-linear equation $g(\tha)=0$, in this paper, we have established a general framework of a new root--finding method, called as the US algorithm. Its each iteration consists of a \textsf{U-step} and an \textsf{S-step}, where the \textsf{U-step} is to construct one $U$-function $U(\tha |\tht)$ based on a new notion of changing direction inequality for the original function $g(\tha)$, while the \textsf{S-step} solves the $U$-equation $U(\tha | \tht)=0$ with respect to the left argument to obtain the next iterate $\thtI$. Unlike the NR algorithm, which is sensitive to initial values and may diverge if a poor initial value is chosen, from Theorem 1, we know that the US algorithm has the property of strongly stable convergence; that is, it does not depend on any initial values in $\Tha$ and strongly stably converges to the root $\ths$. Especially, because of the property of strongly stable convergence, the US algorithm could be one of the powerful tools for solving an equation with multiple roots as shown in Section 6.3.

The critical requirement for applying the US algorithm is to construct one $U$-function $U(\tha|\tht)$ for $g(\tha)$ at $\tha=\tht$ (i.e., to satisfy the CD inequality $g(\tha) - U(\tha|\tht) \sge{\sgn(\tha-\tht)} 0$) such that an explicit solution to the surrogate equation $U(\tha|\tht)=0$ is available. By using Taylor expansion, we have developed four methods for constructing $U$-functions based on the FLB function, SLUB constants, TLB constant and fixed--block. The resulting US algorithms have linear/quadratic/cubic/linear convergence rates as shown in Theorems 2--5, respectively.

Although they enjoy strongly stable convergence, the US algorithms based on both the FLB function method and the fixed--block method hold a linear convergence rate only. In Theorem 6, we propose a simple acceleration technique as shown in (\ref{uclaone6.1}) for the two US algorithms, resulting in a weakly stable convergence. In addition, Theorem 6 showed that the convergence speed of the proposed fast US algorithm is twice of the convergence speed of the US algorithms based on the FLB function method and fixed--block method, even has a super-linear convergence rate.

In the future studies, we will develop other methods for constructing $U$-functions. We would also like to accelerate the US
algorithm with slow convergence rate if the bound functions are not tight enough, particularly for the FLB function-based US algorithms. Other accelerating techniques for the proposed US algorithms may significantly spark the method's potential.

\section*{Acknowledgements}

Guo-Liang TIAN's research was partially supported by National Natural Science Foundation of China (No. 12171225).

\section*{References}
\begin{description} \itemsep=-\parsep \itemindent=-1.3cm
  \item[ ] Arbous AG and Kerrich JE (1951). Accident statistics and the concept of accident--proneness. \textit{Biometrics} \textbf{7}(4): 340--432.

  \item[ ] Azzalini A (1985). A class of distributions which includes the normal ones. \textit{Scandinavian Journal of Statistics} \textbf{12}(2): 171--178.

  \item[ ] Bates GE and Neyman J (1952). Contributions to the theory of accident proneness. 1. an optimistic model of the correlation between light and severe accidents.  \textit{University of California Publications in Statistics} \textbf{1}(9): 215--254.

  \item[ ] Boyd S and Vandenberghe L (2004). \textit{Convex Optimization}. Cambridge: Cambridge University Press.

  \item[ ] Chabert JL (1999). \textit{A History of Algorithms: From the Pebble to the Microchip}, Chapter 3: Methods of false position.  Berlin: Springer.

  \item[ ] Chesneau C (2021). A note on an extreme left skewed unit distribution: Theory, modelling and data fitting. \textit{Open Statistics} {\bf 2}(1): 1--23.

  \item[ ] Consul PC (1990). Geeta distribution and its properties. \textit{Communications in Statistics --- Theory and Methods} {\bf 19}(8): 3051--3068.

  \item[ ] Consul PC and Famoye F (1992). Generalized Poisson regression model. \textit{Communications in Statistics---Theory and Methods} {\bf 21}(1): 89--109.

  \item[ ] Consul PC and Jain GC (1973). A generalization of the Poisson distribution. {\it Technometrics} {\bf 15}(4): 791--799.

  \item[ ] Costabile F, Gualtieri MI and Luceri R (2006). A modification of Muller's method. \textit{Calcolo}, {\bf 43}(1): 39--50.

  \item[ ] Dempster AP, Laird NM and Rubin DB (1977). Maximum likelihood from incomplete data via the EM algorithm (with discussions). \textit{Journal of the Royal Statistical Society: Series B (Methodological)} {\bf 39}: 1--38.

  \item[ ] Garcia JMG (2011). A fixed--point algorithm to estimate the Yule--Simon distribution parameter. \textit{Applied Mathematics and Computation} \textbf{217}(21): 8560--8566.

  \item[ ] Hallinan JAJ (1993). A review of the Weibull distribution. \textit{Journal of Quality Technology} \textbf{25}(2): 85--93.

  \item[ ] Hunter DR and Lange K (2004). A tutorial on MM algorithms. \textit{The American Statistician} \textbf{58}(1): 30--37.

  \item[ ] Lange K (2016). {\it MM Optimization Algorithms.} Philadelphia: SIAM.

  \item[ ] Lange K, Hunter DR and Yang I (2000). Optimization transfer using surrogate objective functions (with discussions). \textit{Journal of Computational and Graphical Statistics} {\bf 9}, 1--20.

  \item[ ] Leisen F, Rossini L and Villa C (2017). A note on the posterior inference for the Yule--Simon distribution. \textit{Journal of Statistical Computation and Simulation} \textbf{87}(6), 1179--1188.

  \item[ ] Nair U, Sankaran PG and Balakrishnan N (2018). \textit{Reliability Modelling and Analysis in Discrete Time} Chapter 3: Discrete Lifetime Models. London: Academic Press.

  \item[ ] Rinne H (2008). \textit{The Weibull Distribution: A Handbook}. Florida: CRC Press.

  \item[ ] Wolfe P (1959). The secant method for simultaneous nonlinear equations. \textit{Communications of the ACM} \textbf{2}(12): 12--13.

  \item[ ] Wood GR (1992). The bisection method in higher dimensions. \textit{Mathematical Programming} \textbf{55}: 319--337.

  \item[ ] Yule GU (1925). A mathematical theory of evolution, based on the conclusions of Dr. JC Willis, FR S. \textit{Philosophical Transactions of the Royal Society of London. Series B, Containing Papers of A Biological Character} \textbf{213}(402--410): 21--87.
\end{description}

\newpage
\baselineskip 0.30in
\renewcommand{\theequation}{\arabic{section}.\arabic{equation}}
\setcounter{equation}{0}
\renewcommand{\theequation}{A.\arabic{equation}}
\renewcommand{\thesubsection}{A.\arabic{subsection}}
\section*{Appendix A: \ Proofs of Theorems 1--6}

\noi\textbf{Proof of Theorem 1}. (i) We first prove that the US algorithm strongly stably converges to the unique root $\ths$. Let $\tht$ be the $t$-th iteration of $\ths$, then we only need to prove that
\begin{eqnarray}
   \tht < \thtI \le \ths, & \mif \tht < \ths, \label{uclaoneA.1} \\ [2mm]
   \ths \le \thtI < \tht, & \mif \tht > \ths. \label{uclaoneA.2}
\end{eqnarray}
For the case of $\tht < \ths$, note that $U(\tha|\tht)$ is a $U$-function for $g(\tha)$ at $\tha=\tht$, then we have
\begin{eqnarray}
   U(\tha|\tht) & \sge{(\ref{uclaone2.2})\, \& \, (\ref{uclaone2.3})} & g(\tha) \stackrel{(\ref{uclaone2.1})}{>} 0, \quad \mif \tha \le \tht (<\ths), \label{uclaoneA.3} \\ [2mm]
   U(\tha|\tht) & \sle{(\ref{uclaone2.4})} & g(\tha), \hspace*{1.45cm} \mif  \tha > \tht. \label{uclaoneA.4}
\end{eqnarray}
Let $\thtI$ be the $(t+1)$-th US iteration; i.e., $U(\thtI|\tht)=0$. Then we can show
\begin{eqnarray} \label{uclaoneA.5}
   \thtI > \tht.
\end{eqnarray}
In fact, if (\ref{uclaoneA.5}) is not true (i.e., $\thtI \le \tht$ is true), by replacing $\tha$ in (\ref{uclaoneA.3}) with $\thtI$, we have $0 = U(\thtI|\tht) \ge g(\thtI) > 0$, which is contradictious. So (\ref{uclaoneA.5}) is true.

Next, we can show that
\begin{eqnarray} \label{uclaoneA.6}
   \thtI \le \ths.
\end{eqnarray}
In fact, if (\ref{uclaoneA.6}) is not true (i.e., $\thtI > \ths$ is true), by replacing $\tha$ in (\ref{uclaoneA.4}) with $\thtI$, we have $g(\thtI) \ge U(\thtI|\tht)= 0$, which
is contradictory to the fact that $g(\tha) < 0$ for any $\tha>\ths$ as shown by (\ref{uclaone2.1}). So (\ref{uclaoneA.6}) is true. By combining (\ref{uclaoneA.5}) with (\ref{uclaoneA.6}), we obtain (\ref{uclaoneA.1}).

Finally, we prove that
$$
\lim_{t\rightarrow \infty} \tht = \ths.
$$
If $\lim_{t\rightarrow \infty} \tht < \ths$, there must exist a positive constant $c$ such that $g(\tha)\ge c >0$ for any $\tha^{(0)} \le \tha \le \lim_{t\rightarrow \infty} \tht$, then, by using the lipschitz condition, we have
\begin{eqnarray*}
   g(\tht) & \sle{(\ref{uclaone2.2})} & U(\tht|\tht) \; \sr{(\ref{uclaone2.6})} \;  U(\tht|\tht) - U(\thtI|\tht) \\
   &\le & |U(\thtI|\tht) - U(\tht|\tht)| \; \sle{ \footnotesize \mbox{Lipschitz}}\; L (\thtI - \tht),
\end{eqnarray*}
and then
$$
   \tht - \tha^{(t-1)} \;\ge\; \frac{g(\tha^{(t-1)})}{L} \;\ge\; \frac{c}{L} , \quad \forall \; t = 1,\ldots,\infty.
$$
Since
$$
   \tht - \tha^{(0)} \;=\; \sum_{i=1}^{t} (\tha^{(i)} - \tha^{(i-1)}) \;\ge\; \sum_{i=1}^{t} \frac{c}{L} \;=\; \frac{c}{L} t,
$$
we have
$$
   \tht \;\ge\; \tha^{(0)} + \frac{c}{L} t \;\rightarrow \; \infty \mbox{ when } t \;\rightarrow\; \infty,
$$
which is contradictory to the fact that $\tht \le \ths < \infty$. For the case of $\tht > \ths$, similarly, we can prove (\ref{uclaoneA.2}) and $\lim_{t\rightarrow \infty} \tht = \ths$.

(ii) From (\ref{uclaoneA.1}) and (\ref{uclaoneA.2}), we can see that the US algorithm does not depend on any initial values in $\Tha$. \hfill $\Box$

\vkU\noi\textbf{Proof of Theorem 2}. From (\ref{uclaone3.3}), we have
$$
   |U'(\tha|\tht)| \;=\; |b(\tha)| \; \le \; \max_{\tha\in \Tha_1} |b(\tha)| \;\teq\; L_1,
$$
where $\Tha_1$ is the interval between $\tha^{(0)}$ and $\ths$. Thus, the $U$-function $U(\tha|\tht)$ in (\ref{uclaone3.3}) satisfies the lipschitz condition $|U(\tha_1|\tht) - U(\tha_2|\tht)| \le L_1 |\tha_1 - \tha_2|$ for any $\tha_1,\tha_2 \in \Tha_1$. From (\ref{uclaone5.1}), we have
$$
   0 = U(h(\tht)|\tht)  = g(\tht) + \int_{\tht}^{h(\tht)} b(z) \rd z.
$$
By taking derivative with respect to $\tht$ on the both sides of the above equation, we obtain $0 = g'(\tht) + b(h(\tht)) h'(\tht) - b(\tht)$,
so that
$$
   h'(\tht) =  \frac{b(\tht) - g'(\tht)}{b(h(\tht))}.
$$
Note that the rate of convergence for the US algorithm is
\begin{eqnarray*}
   |h'(\ths)| &=& \lim_{t\to \infty} |h'(\tht)| = \lim_{t\to \infty} \left| \frac{b(\tht) - g'(\tht)}{b(h(\tht))} \right| \\ [2mm]
   &=& \left| \frac{b(\ths) - g'(\ths)}{b(h(\ths))} \right| = \left| \frac{b(\ths) - g'(\ths)}{b(\ths)} \right|  = \left| 1- \frac{g'(\ths)}{b(\ths)} \right|,
\end{eqnarray*}
we only need to show that
\begin{eqnarray} \label{uclaoneA.7}
   \left| 1- \frac{g'(\ths)} {b(\ths)} \right| = 1- \frac{g'(\ths)} {b(\ths)} \in [0, 1).
\end{eqnarray}
If the condition, $b(\ths) \le g'(\ths) <0$, is true, then it is easy to check the correctness of (\ref{uclaoneA.7}). From (\ref{uclaone3.2}), we of course have $b(\ths) \le g'(\ths)$.

Now, we only need to prove $g'(\ths) < 0$. For any $\ve >0$, from the assumption (\ref{uclaone2.1}), we have $g(\ths+\ve)<0$ and $g(\ths-\ve)>0$ so that we obtain
\begin{eqnarray}
   g_+'(\ths) &\teq & \lim_{\ve \to 0^+}\frac{g(\ths+\ve)-g(\ths)}{\ve} = \lim_{\ve \to 0^+} \frac{g(\ths+\ve)}{\ve} < 0 \qand \label{uclaoneA.8} \\ [2mm]
   g_-'(\ths) &\teq & \lim_{\ve \to 0^+}\frac{g(\ths-\ve)-g(\ths)}{-\ve} =
   \lim_{\ve \to 0^+}\frac{g(\ths-\ve)}{-\ve} < 0. \label{uclaoneA.9}
\end{eqnarray}
Since the first--order derivative $g'(\cdot)$ exists, we have $g'(\ths) = g_+'(\ths) = g_-'(\ths) < 0$.   \hfill $\Box$

\vkU\noi\textbf{Proof of Theorem 3}.  From (\ref{uclaone3.9}), we have
$$
   |U'(\tha|\tht)| \;=\; |g'(\tht) + b_2(\tha| \tht) (\tha-\tht)| \; \le \; \max_{\tha\in \Tha_2} |g'(\tht) + b_2(\tha| \tht) (\tha-\tht)| \;\teq\; L_2,
$$
where $\Tha_2$ is the interval between $\tha^{(0)}$ and $\ths$. Thus, the $U$-function $U(\tha|\tht)$ in (\ref{uclaone3.9}) satisfies the lipschitz condition $|U(\tha_1|\tht) - U(\tha_2|\tht)| \le L_2 |\tha_1 - \tha_2|$ for any $\tha_1,\tha_2 \in \Tha_2$.
From (\ref{uclaone5.3}), we have
\begin{eqnarray} \label{uclaoneA.10}
   g(\tht) + g'(\tht) [h(\tht)-\tht] + \frac{b_2(h(\tht), \tht)} {2} [h(\tht)-\tht]^2 =0.
\end{eqnarray}
By taking the first-derivative about $\tht$ on the both sides of (\ref{uclaoneA.10}), we obtain
$$
    g'(\tht) h'(\tht) + \Big\{ g''(\tht) + b_2(h(\tht), \tht) [h'(\tht)-1] \Big\}[h(\tht)-\tht] = 0.
$$
By letting $t\to  \infty$, we have $g'(\ths) h'(\ths) = 0$ because
\begin{eqnarray} \label{uclaoneA.11}
   \lim_{t\to \infty} [h(\tht)-\tht] = h(\ths) - \ths = 0.
\end{eqnarray}
Usually, we know that $g'(\ths)\neq 0$ so that
\begin{eqnarray} \label{uclaoneA.12}
   h'(\ths)=0.
\end{eqnarray}
By taking the second-derivative with respect to $\tht$ on the both sides of (\ref{uclaoneA.10}), we have
\begin{eqnarray*}
   0 &=&g''(\tht) h'(\tht) + g'(\tht) h''(\tht) + \Big\{g'''(\tht) + b_2(h(\tht), \tht) h''(\tht) \Big\}[h(\tht)-\tht] \\ [2mm]
   & & +\; \Big\{g''(\tht)+ b_2(h(\tht), \tht) [h'(\tht)-1] \Big\} [h'(\tht)-1].
\end{eqnarray*}
By letting $t\to \infty$ and using (\ref{uclaoneA.11}) \& (\ref{uclaoneA.12}), we obtain
$$
   |h''(\ths)| = \lim_{t\to \infty} |h''(\tht)| =   \left|\frac{b_2(\ths, \ths)-g''(\tha^*)}{g'(\tha^*)} \right| \; \sr{(\ref{uclaone3.7})} \; \left|\frac{b_{22} - g''(\tha^*)}{g'(\tha^*)} \right|,
$$
indicating that the US iteration (\ref{uclaone5.3}) has a quadratic convergence rate given by (\ref{uclaone5.4}). \hfill $\Box$

\vkU\noi\textbf{Proof of Theorem 4}.  From (\ref{uclaone3.15}), we have
\begin{eqnarray*}
   |U'(\tha|\tht)| &=& |g'(\tht) + g''(\tht) (\tha-\tht) + \frac{1}{3} b_3 (\tha-\tht)^2| \\
   &\le& \max_{\tha\in \Tha_3} |g'(\tht) + g''(\tht) (\tha-\tht) + \frac{1}{3} b_3 (\tha-\tht)^2| \;\teq\; L_3,
\end{eqnarray*}
where $\Tha_3$ is the interval between $\tha^{(0)}$ and $\ths$. Thus, the $U$-function $U(\tha|\tht)$ in (\ref{uclaone3.15}) satisfies the lipschitz condition $|U(\tha_1|\tht) - U(\tha_2|\tht)| \le L_3 |\tha_1 - \tha_2|$ for any $\tha_1,\tha_2 \in \Tha_3$.
From (\ref{uclaone5.5}), we have
\begin{eqnarray}  \label{uclaoneA.13}
   0 = g(\tht) + g'(\tht)[h(\tht)-\tht] + \frac{g''(\tht)}{2}  [h(\tht)-\tht]^2  + \frac{b_3}{6}  [h(\tht)-\tht]^3. \quad
\end{eqnarray}
By taking the first-derivative about $\tht$ on the both sides of (\ref{uclaoneA.13}), we obtain
\begin{eqnarray*}
    0 &=& g'(\tht) h'(\tht) + g''(\tht)[h(\tht)-\tht]+ \frac{g'''(\tht)} {2} [h(\tht)-\tht]^2  \\ [2mm]
    & & + \; g''(\tht) [h(\tht)-\tht][h'(\tht)-1]  + \frac{b_3}{2} [h(\tht)-\tht]^2 [h'(\tht)-1].
\end{eqnarray*}
By letting $t\to \infty$, we have $h'(\ths) = 0$ when $g'(\ths) \ne 0$, because
\begin{eqnarray}  \label{uclaoneA.14}
   \lim_{t\to \infty} [h(\tht)-\tht] = h(\ths) - \ths = 0.
\end{eqnarray}
By taking the second-derivative about $\tht$ on the both sides of (\ref{uclaoneA.13}), we obtain
\begin{eqnarray*}
    0 &=& g''(\tht) h'(\tht) + g'(\tht) h''(\tht) + g'''(\tht)[h(\tht)-\tht]+ g''(\tht)[h'(\tht)-1] \\ [2mm]
    & & + \; \frac{g^{(4)}(\tht)}{2} [h(\tht)-\tht]^2 + g'''(\tht) [h(\tht)-\tht][h'(\tht)-1] \\ [2mm]
    & & + \; \left\{g'''(\tht) + \frac{b_3}{2} [h'(\tht)-1] \right\} \cdot [h(\tht)-\tht][h'(\tht)-1] \\ [2mm]
    & & + \; \left\{g''(\tht) + \frac{b_3}{2} [h(\tht)-\tht] \right\} \cdot \Big\{[h'(\tht)-1]^2+[h(\tht)-\tht]h''(\tht) \Big\}.
\end{eqnarray*}
By letting $t\to \infty$, we have $h''(\ths) = 0$, because $h'(\ths) = 0$ and (\ref{uclaoneA.14}). By taking the third-derivative about $\tht$ on the both sides of equation (\ref{uclaoneA.13}), we obtain
\begin{eqnarray*}
    0 &=& g'''(\tht) h'(\tht) + g''(\tht) h''(\tht) + g''(\tht) h''(\tht)+ g'(\tht) h'''(\tht) \\ [2mm]
    & & + \; g^{(4)}(\tht)[h(\tht)-\tht]+ g'''(\tht)[h'(\tht)-1] + g'''(\tht)[h'(\tht)-1] \\ [2mm]
    & &+\;  g''(\tht)h''(\tht) + \frac{1}{2} g^{(5)}(\tht) [h(\tht)-\tht]^2 +  g^{(4)}(\tht) [h(\tht)-\tht][h'(\tht)-1] \\ [2mm]
    & &+\; \left\{2g^{(4)}(\tht) + \frac{b_3}{2} h''(\tht)\right\}\cdot [h(\tht)-\tht][h'(\tht)-1] \\ [2mm]
    & &+ \;\left\{2g'''(\tht) + \frac{b_3}{2}[h'(\tht)-1]\right\}\cdot \left\{[h'(\tht)-1]^2+[h(\tht)-\tht]h''(\tht)\right\} \\ [2mm]
    & &+ \;\left\{g'''(\tht)   + \frac{b_3}{2} [h'(\tht)-1]\right\}\cdot \left\{[h'(\tht)-1]^2 + [h(\tht)-\tht]h''(\tht)\right\} \\  [2mm]
    & &+ \;\left\{g''(\tht)   + \frac{b_3}{2} [h(\tht)-\tht]\right\}\cdot \left\{3[h'(\tht)-1]h''(\tht) +[h(\tht)-\tht]h'''(\tht)\right\}.
\end{eqnarray*}
By letting $t\to \infty$ and using $h'(\ths) = 0$, (\ref{uclaoneA.14}) and $h''(\ths) = 0$, we obtain
$$
   |h'''(\ths)| =  \bigg|\frac{b_3 - g'''(\ths)}{g'(\ths)}\bigg|,
$$
indicating that the US iteration (\ref{uclaone5.5}) has a cubic convergence rate given by (\ref{uclaone5.6}). \hfill $\Box$

\vkU\noi\textbf{Proof of Theorem 5}. For convenience, we define
\begin{eqnarray} \label{uclaoneA.15}
  G \Big(u_1(\tha), \ldots, u_r(\tha), u_{r+1} \big( h(\tha) \big), \ldots, u_m \big(h(\tha) \big) \Big) \teq  G[\tha, h(\tha)].
\end{eqnarray}
From (\ref{uclaone5.7}), we know that $\thtI = h(\tht)$ is the root of $U_{1\cdots r}(\thtI |\tht)$ $=0$; i.e., we have
\begin{eqnarray*}
   0 = U_{1\cdots r}(h(\tht)|\tht) & \sr{(\ref{uclaone3.17})} &
   G \Big(u_1(\tht), \ldots, u_r(\tht), u_{r+1}\big(h(\tht) \big), \ldots, u_m \big(h(\tht) \big) \Big) \\ [2mm]
   & \sr{(\rm \ref{uclaoneA.15})} & G[\tht, h(\tht)].
\end{eqnarray*}
By taking derivative with respect to $\tht$ on the both sides of the above equation, we obtain
\begin{eqnarray*}
   0 &=& \sum_{j=1}^r \frac{\pa G[\tht, h(\tht)]}{\pa u_j(\tht)} u_j'(\tht) + \sum_{j=r+1}^m \frac{\pa G[\tht, h(\tht)]}{\pa u_j \big(h(\tht) \big)} u_j' \big( h(\tht) \big) h'(\tht),
\end{eqnarray*}
where $u_j'(\tha) \teq \rd u_j(\tha) /\rd \tha$. Thus, we have
\begin{eqnarray*}
   h'(\tht) &=& -\left[\sum_{j=r+1}^m \frac{\pa G[\tht, h(\tht)] }{\pa u_j \big( h(\tht) \big)} u_j' \big( h(\tht) \big)\right]^{-1} \sum_{j=1}^r \frac{\pa G[\tht, h(\tht)] }{\pa u_j(\tht)} u_j'(\tht).
\end{eqnarray*}
Note that the rate of convergence for the US algorithm is
\begin{eqnarray*}
   |h'(\ths)| &=& \lim_{t\to \infty} |h'(\tht)|
   = \left| \left[\sum_{j=r+1}^m  \frac{\pa G[\ths, h(\ths)] }{\pa u_j \big(h(\ths) \big)} u_j' \big(h(\ths) \big)\right]^{-1} \sum_{j=1}^r  \frac{\pa G[\ths, h(\ths)] }{\pa u_j(\ths)} u_j'(\ths) \right| \\ [2mm]
   &=& \left| \left[\sum_{j=r+1}^m  \frac{\pa G[\ths, \ths] }{\pa u_j(\ths)} u_j'(\ths)\right]^{-1} \sum_{j=1}^r  \frac{\pa G[\ths, \ths] }{\pa u_j(\ths)} u_j'(\ths) \right| \\ [2mm]
   &=& \left| \frac{g'(\ths)- U_{1\cdots r}'(\ths|\ths)}{U_{1\cdots r}'(\ths|\ths)} \right| = \left|\frac{g'(\ths)}{U_{1\cdots r}'(\ths|\ths)} - 1 \right|,
\end{eqnarray*}
where $G[\ths, \ths] = G\big(u_1(\ths), \ldots, u_m(\ths) \big)$,
\begin{eqnarray}
   U_{1\cdots r}'(\ths|\ths) &\teq & U_{1\cdots r}'(\tha|\ths)\big|_{\tha=\ths} \non \\ [2mm]
   &=& \sum_{j=r+1}^m  \frac{\pa G \big( u_1(\ths), \ldots, u_r(\ths), u_{r+1}(\tha), \ldots, u_m(\tha) \big)} {\pa u_j(\tha)} u_j'(\tha)\bigg|_{\tha=\ths} \non \\ [2mm]
   &=&  \sum_{j=r+1}^m  \frac{\pa G (u_1(\ths), \ldots, u_m(\ths))}{\pa u_j(\ths)} u_j'(\ths) = \sum_{j=r+1}^m  \frac{\pa G [\ths, \ths]} {\pa u_j(\ths)} u_j'(\ths), \qquad  \label{uclaoneA.16}  \\ [2mm]
   g'(\ths) &\teq& g'(\tha)\big|_{\tha=\ths} = \sum_{j=1}^m  \frac{\pa G \big(u_1(\tha), \ldots, u_m(\tha) \big)}{\pa u_j(\tha)} u_j'(\tha) \bigg|_{\tha=\ths} \non \\ [2mm]
   &=& \sum_{j=1}^m  \frac{\pa G (u_1(\ths), \ldots, u_m(\ths))}{\pa u_j(\ths)} u_j'(\ths) = \sum_{j=1}^m  \frac{\pa G [\ths, \ths]} {\pa u_j(\ths)} u_j'(\ths).  \label{uclaoneA.17}
\end{eqnarray}
We only need to prove
\begin{eqnarray} \label{uclaoneA.18}
   \left|\frac{g'(\ths)} {U_{1\cdots r}'(\ths|\ths)} -1 \right| = 1- \frac{g'(\ths)} {U_{1\cdots r}'(\ths|\ths)} \in [0, 1).
\end{eqnarray}
Note that
\begin{eqnarray*}
   0 \srg {(\rm \ref{uclaoneA.8})\, \& \, (\rm \ref{uclaoneA.9})} g'(\ths) & \sr{(\rm \ref{uclaoneA.17})} & \sum_{j=1}^r  \frac{\pa G [\ths, \ths]} {\pa u_j(\ths)} u_j'(\ths) + \sum_{j=r+1}^m  \frac{\pa G [\ths, \ths]} {\pa u_j(\ths)} u_j'(\ths)   \\ [2mm]
   & \sge{(\ref{uclaone3.16})} & \sum_{j=r+1}^m  \frac{\pa G [\ths, \ths]} {\pa u_j(\ths)} u_j'(\ths) \sr{(\rm \ref{uclaoneA.16})} U_{1\cdots r}'(\ths|\ths),
\end{eqnarray*}
then it is easy to check the correctness of (\ref{uclaoneA.18}). \hfill $\Box$

\vkU\noi\textbf{Proof of Theorem 6}. (i) From  (\ref{uclaone2.2})--(\ref{uclaone2.4}), for any $\ve >0$,  we have
\begin{eqnarray*}
   g(\tht-\ve) &\le & U(\tht-\ve|\tht), \quad g(\tht)= U(\tht|\tht)
   \qand \\ [2mm]
   g(\tht+\ve) &\ge & U(\tht+\ve|\tht),
\end{eqnarray*}
so that we obtain
\begin{eqnarray*}
   g_-'(\tht) &\teq & \lim_{\ve \to 0^+}\frac{g(\tht-\ve)-g(\tht)}{-\ve} \\ [2mm]
   &\ge & \lim_{\ve \to 0^+}\frac{U(\tht-\ve|\tht)-U(\tht|\tht)}{-\ve} \teq U_-'(\tht|\tht) \qand \\ [2mm]
   g_+'(\tht) &\teq & \lim_{\ve \to 0^+} \frac{g(\tht +\ve)-g(\tht)}{\ve} \\ [2mm]
   &\ge & \lim_{\ve \to 0^+}\frac{U(\tht+\ve|\tht)-U(\tht|\tht)}{\ve} \teq U_+'(\tht|\tht),
\end{eqnarray*}
implying
\begin{eqnarray} \label{uclaoneA.19}
   U'(\tht|\tht) \le g'(\tht).
\end{eqnarray}
Hence
$$
   s(\tht) \sr{(\ref{uclaone6.2})} \Cases{\min(1, 2)=1} {\mif g'(\tht)\ge 0} {\min \left\{\dis\frac{U'(\tht|\tht)}{g'(\tht)} \, \left(\sge{\rm{(\ref{uclaoneA.19})}} 1 \right), \; 2\right\} } {\mif g'(\tht)< 0,}
$$
so that
\begin{eqnarray} \label{uclaoneA.20}
   1 \le s(\tht) \le 2.
\end{eqnarray}

We consider two cases: $\tht<\ths$ and $\tht >\ths$. For Case I: $\tht<\ths$, from (\ref{uclaone6.3}), we have
\begin{eqnarray} \label{uclaoneA.21}
   \tht<\tth^{(t+1)} < \ths.
\end{eqnarray}
Thus,
\begin{eqnarray*}
   \thtI-\ths & \sr{(\ref{uclaone6.1})} & \tht + s(\tht) (\tth^{(t+1)}- \tht) -\ths  \sge{\rm{(\ref{uclaoneA.20})}} \tht + (\tth^{(t+1)}- \tht) -\ths \\ [2mm]
   &=& \tth^{(t+1)} -\ths \srg{\rm{(\ref{uclaoneA.21})}} \tht -\ths \qand 
\end{eqnarray*}
\begin{eqnarray*}
   \thtI-\ths & \sr{(\ref{uclaone6.1})} & \tht + s(\tht) (\tth^{(t+1)}- \tht) -\ths \sle{\rm{(\ref{uclaoneA.20})}}  \tht + 2(\tth^{(t+1)}- \tht)-\ths \\ [2mm]
   &=& (\tth^{(t+1)} -\tht) + (\tth^{(t+1)} -\ths)  \srl{\rm{(\ref{uclaoneA.21})}}  \tth^{(t+1)} -\tht \srl{\rm{(\ref{uclaoneA.21})}} \ths -\tht,
\end{eqnarray*}
so that we have
\begin{eqnarray} \label{uclaoneA.22}
   |\thtI-\ths| < |\tht-\ths|, \quad \forall \, t=0, 1, \ldots, \infty,
\end{eqnarray}
implying that the sequence $\{\thtI\}_{t=0}^\infty$ specified by (\ref{uclaone6.1}) weakly converges to $\ths$.

For Case \II: $\tht>\ths$, from (\ref{uclaone6.3}), we have
\begin{eqnarray} \label{uclaoneA.23}
   \ths < \tth^{(t+1)} < \tht.
\end{eqnarray}
Thus,
\begin{eqnarray*}
   \thtI-\ths & \sr{\rm{(\ref{uclaone6.1})}} & \tht + s(\tht) (\tth^{(t+1)}- \tht) -\ths \sle{\rm{(\ref{uclaoneA.20})}}  \tht + (\tth^{(t+1)}- \tht)-\ths \\ [2mm]
   &=&  \tth^{(t+1)} -\ths \srl{\rm{(\ref{uclaoneA.23})}} \tht -\ths  \qand \\ [2mm]
   \thtI-\ths & \sr{\rm{(\ref{uclaone6.1})}} & \tht + s(\tht) (\tth^{(t+1)}- \tht) -\ths  \sge{\rm{(\ref{uclaoneA.20})}} \tht + 2 (\tth^{(t+1)}- \tht) -\ths \\ [2mm]
   &=& (\tth^{(t+1)} -\tht) + (\tth^{(t+1)} -\ths) \srg{\rm{(\ref{uclaoneA.23})}} \tth^{(t+1)} -\tht \srg{\rm{(\ref{uclaoneA.23})}} -(\tht -\ths),
\end{eqnarray*}
so that we have (\ref{uclaoneA.22}) again.

(ii) To derive the convergence rate for the fast US algorithm, we define $\thtI \teq h(\tht)$ and $\tth^{(t+1)} \teq \tilde{h}(\tht)$. By taking derivative about $\tht$ on the both sides of (\ref{uclaone6.1}), we obtain
\begin{eqnarray} \label{uclaoneA.24}
   h'(\tht) = 1 + s'(\tht) (\tth^{(t+1)}-\tht) +s(\tht) \big[\tilde{h}'(\tht) - 1 \big].
\end{eqnarray}
Note that the convergence rate for the fast US algorithm is
\begin{eqnarray} \label{uclaoneA.25}
   |h'(\ths)| = \lim_{t\to \infty} |h'(\tht)|
   \sr{\rm{(\ref{uclaoneA.24})}} \big| 1 + s(\ths)[\tilde{h}'(\ths) - 1] \big|.
\end{eqnarray}
From (\ref{uclaoneA.19}), we have $U'(\ths|\ths) \le g'(\ths) <0$, so $s(\ths) \sge{(\ref{uclaoneA.20})} 1$. For the US algorithm based on the FLB function method, we have
\begin{eqnarray}
   \tilde{h}'(\ths) & \sr{\rm{(\ref{uclaone5.2})}} & 1- \frac{g'(\ths)}{b(\ths)} \in [0, 1), \label{uclaoneA.26} \\ [2mm]
   U'(\ths|\ths) & \sr{\rm{(\ref{uclaone5.1})}} & b(\ths) \qand \label{uclaoneA.27} \\ [2mm]
   s(\ths) & \sr{\rm{(\ref{uclaone6.2})\, \& \, (\ref{uclaoneA.27})}} & \min\left\{ \frac{b(\ths)}{g'(\ths)} ,\; 2\right\}. \label{uclaoneA.28}
\end{eqnarray}
If $b(\ths)/g'(\ths) > 2$; i.e., $b(\ths) < 2 g'(\ths) < 0$, then $s(\ths) \sr{(\ref{uclaoneA.28})} 2$, and we have
\begin{eqnarray*}
   h'(\ths) & \sr{\rm{(\ref{uclaoneA.25}) \, \& \, (\ref{uclaoneA.26})}} &  1 + 2\left[1- \frac{g'(\ths)}{b(\ths)} - 1\right] = 1 -  \frac{2g'(\ths)}{b(\ths)} \in (0, 1) \\ [2mm]
   & < &   1 - \frac{g'(\ths)}{b(\ths)}   \; \sr{\rm{(\ref{uclaoneA.26})}} \; \tilde{h}'(\ths),
\end{eqnarray*}
indicating that the convergence speed of the fast US algorithm, defined as $1-h'(\ths) = 2g'(\ths)/b(\ths)$,  is twice of the convergence speed of the US algorithm based on the FLB function method, defined as $1- \tilde{h}'(\ths) = g'(\ths)/b(\ths)$.

If $1 \le b(\ths)/g'(\ths) \le 2$, then $s(\ths) \sr{\rm{(\ref{uclaoneA.28})}} b(\ths)/g'(\ths)$, and we have
$$
   h'(\ths) \sr{\rm{(\ref{uclaoneA.25}}) \, \& \, (\ref{uclaoneA.26})}   1 + \frac{b(\ths)}{g'(\ths)} \left[1- \frac{g'(\ths)}{b(\ths)} - 1\right]  = 0
   \le   1 - \frac{g'(\ths)}{b(\ths)}  \; \sr{\rm{(\ref{uclaoneA.26})}} \;  \tilde{h}'(\ths),
$$
indicating that the fast US algorithm has a super--linear convergence rate comparing with the US algorithm based on the FLB function method.

(iii) By replacing $b(\ths)$ with $U'(\ths |\ths)$ in the proof process of (ii), we can obtain a similar conclusion. \hfill $\Box$

\baselineskip 0.30in
\renewcommand{\theequation}{\arabic{section}.\arabic{equation}}
\setcounter{equation}{0}
\renewcommand{\theequation}{B.\arabic{equation}}
\renewcommand{\thesubsection}{B.\arabic{subsection}}
\section*{Appendix B: \  Mode of the skew normal distribution}

Let $x_{\rm mod}$ denote the mode of the skew normal distribution with pdf given by (\ref{uclaone4.5}), in this appendix, we will propose an MM algorithm to iteratively calculate
$$
   x_{\rm mod} =  \arg\; \max_{x\in \bbR}  \log[f(x |\mu, \si^2, \al)] = \arg\; \max_{x\in\bbR} \left\{-\frac{(x-\mu)^2}{2\si^2} + \log\left[\Phi \left(\frac{x-\mu}{\si_*}\right) \right]\right\},
$$
where $\si_* \teq \si/\al$. Let $\bth =(\mu, \si_*^2)^{\T}$ and define the left-truncated normal density as
$$
   h(x|x^{(s)}, \bth) \teq  \phi\left(x+x^{(s)}-\mu |\bth\right) \cdot \Phi^{-1} \left( \frac{x^{(s)}-\mu}{\si_*}\right) \cdot I(x<\mu) ,
$$
where $\phi(\cdot |\bth)$ denotes the pdf of $N(\mu, \si_*^2)$. By applying the integral version of Jensen's inequality, we have
\begin{eqnarray*}
   \log\left[\Phi \left(\frac{x-\mu}{\si_*} \right) \right]
   &=& \log \left[\int_{-\infty}^x \phi(y |\bth) \rd y \right]  = \log \left[\int_{-\infty}^{\mu} \phi (z+x-\mu |\bth )\rd z\right] \\ [2mm]
   &=& \log \left[\int_{-\infty}^{\mu} \frac{\phi (z+x-\mu |\bth)} {h(z|x^{(s)}, \bth)} \cdot h(z|x^{(s)}, \bth) \rd z\right] \\ [2mm]
   &\ge & \int_{-\infty}^{\mu} \log \left[ \frac{\phi (z+x-\mu |\bth )} {h(z|x^{(s)}, \bth)} \right] \cdot h(z|x^{(s)},\bth)\rd z,
\end{eqnarray*}
where the equality holds iff $x = x^{(s)}$. Thus, the minorizing function can be constructed as
\begin{eqnarray*}
   Q(x|x^{(s)}) &\teq & - \frac{(x-\mu)^2}{2\si^2} + \int_{-\infty}^{\mu} \log \left[ \frac{\phi (z+x-\mu |\bth )} {h(z|x^{(s)}, \bth)}\right] \cdot h(z|x^{(s)},\bth)\rd z \\ [2mm]
   &=&- \;\frac{(x-\mu)^2}{2\si^2} - \int_{-\infty}^{\mu} \frac{(z+x-2\mu)^2}{2\si_*^2} \cdot h(z|x^{(s)},\bth)\rd z + \mbox{constant}\\[3mm]
   &=&- \;\left(\frac{1}{2\si^2} + \frac{1}{2\si_*^2} \right) x^2 + \left[ \frac{\mu}{\si^2} + \frac{2\mu}{\si_*^2} - \frac{1}{\si_*^2} \int_{-\infty}^{\mu} z \cdot h(z|x^{(s)},\bth)\rd z \right] x + \mbox{constant}.
\end{eqnarray*}
Let $\rd Q(x|x^{(s)})/\rd x =0$, we obtain the $(s+1)$-th MM iteration:
\begin{eqnarray} \label{uclaoneB.1}
   x^{(s+1)} = \frac{1}{1+\al^2}\left[ \mu(1 + 2\al^2 )- \al^2 \int_{-\infty}^{\mu} z \cdot h(z| x^{(s)},\bth)\rd z \right].
\end{eqnarray}

\baselineskip 0.30in
\renewcommand{\theequation}{\arabic{section}.\arabic{equation}}
\setcounter{equation}{0}
\renewcommand{\theequation}{C.\arabic{equation}}
\renewcommand{\thesubsection}{C.\arabic{subsection}}
\section*{Appendix C: \  Inequalities on Riemann's zeta function}

First, for any $\tha>1$, we have
\begin{eqnarray}   \label{uclaoneC.1}
   \left(\tha- \frac{3}{2} \right)^2 + \frac{3}{4} >0 \Leftrightarrow \tha^2 - 3\tha + 3 >0 \Leftrightarrow \frac{1}{\tha -1} < \frac{1}{(\tha -1)^2} +1.
\end{eqnarray}
Next, we prove the following key inequality:
\begin{eqnarray} \label{uclaoneC.2}
   \int_1^{\infty} \frac{\log(z+1)}{z^\tha} \rd z &=& \frac{1}{1-\tha} \int_1^{\infty} \log(z+1) \rd z^{1-\tha} \non \\ [2mm]
   &=& \frac{\log(2)}{\tha-1} + \frac{1}{\tha-1} \int_1^{\infty} \frac{z^{1-\tha}}{z+1} \rd z  \non \\ [2mm]
   &< & \frac{\log(2)}{\tha-1} + \frac{1}{\tha-1} \int_1^{\infty} \frac{z^{1-\tha}}{z} \rd z  \non\\ [2mm]
   &=&  \frac{\log(2)}{\tha-1} + \frac{1}{(\tha-1)^2} \; \stackrel{\rm (\ref{uclaoneC.1})} {<} \; \frac{\log(2)+1}{(\tha-1)^2} + \log(2).
\end{eqnarray}
Finally, we obtain
\begin{eqnarray}
   Z(\tha) &=& \sum_{x=1}^{\infty} \frac{1}{x^{\tha}} \ge \int_1^{\infty} \frac{1}{x^{\tha}} \rd x = \frac{1}{\tha-1}, \label{uclaoneC.3} \\ [2mm]
   Z'(\tha)&=& -\sum_{x=2}^{\infty} \frac{\log(x)}{x^\tha} \ge  - \int_1^{\infty} \frac{\log(z+1)}{z^\tha} \rd z \; \stackrel{\rm (\ref{uclaoneC.2})} {>} \; - \frac{\log(2)+1}{(\tha-1)^2} - \log(2), \non \\ [2mm]
   Z''(\tha) &=& \sum_{x=2}^{\infty} \frac{\log^2(x)}{x^\tha} \le \int_1^{\infty} \frac{\log^2(z+1)}{z^\tha} \rd z = \int_1^{\infty} \frac{\log^2(z+1)}{1-\tha} \rd z^{1-\tha} \non \\ [2mm]
   &=&  \frac{\log^2(2)}{\tha-1} + \frac{2}{\tha-1} \int_1^{\infty} \frac{z}{z+1} \cdot z^{-\tha} \log(z+1) \rd z  \non 
\end{eqnarray}
\begin{eqnarray}
   &\le&  \frac{\log^2(2)}{\tha-1} + \frac{2}{\tha-1} \int_1^{\infty} \frac{\log(z+1)}{z^\tha} \rd z \non \\ [2mm]
   &\stackrel{\rm (\ref{uclaoneC.2})} {<}&  \frac{\log^2(2)}{\tha-1} + \frac{2}{\tha-1} \left[\frac{\log(2)+1}{(\tha-1)^2} + \log(2)\right]. \label{uclaoneC.4}
\end{eqnarray}

\baselineskip 0.24in
\renewcommand{\theequation}{\arabic{section}.\arabic{equation}}
\setcounter{equation}{0}
\renewcommand{\theequation}{D.\arabic{equation}}
\renewcommand{\thesubsection}{D.\arabic{subsection}}

\section*{Appendix D: \  \textsf{R} codes with US algorithm for solving the \\ \hspace*{4.0cm} roots of a high--order polynomial equation}

\begin{verbatim}
root.HPE = function(a3, a2, a1, a0, a.max = 10, m, x0 = 0)
{ # Function name: root.HPE(a3, a2, a1, a0, a.max = 10, m, x0 = 0)
  # ---------------------------Aim------------------------------------------
  # Solve the equation: 0 = a3*x^m + a2*x^2 + a1*x + a0, 0 < x < a.max
  # ---------------------------Input----------------------------------------
  #    a3: The coefficient of x^m
  #    a2: The coefficient of x^2
  #    a1: The coefficient of x
  #    a0: The constant, that is a positive real number
  # a.max: The maximum of x
  #     m: The highest order of the equation
  #    x0: An initial value and the default value is 0
  # ---------------------------Output---------------------------------------
  # The unique root of the polynomial equation in (0, a.max)
  ##########################################################################
  xt.save = c(); yt.save = c();
  xt = x0; y1t = a3*xt^m + a2*xt^2 + a1*xt + a0
  xt.save[1] = xt; yt.save[1] = y1t
  k = 1; esp = 1
  while(esp>1e-8){
    if(a3<0){                ### When a3<0, we apply the SLUF method
      if(y1t>0){ b2 = a3*m*(m-1)*a.max^(m-2) }else{b2 = 0}
      A2 = b2/2+a2
      A1 = a3*m*xt^(m-1)-b2*xt + a1
      A0 = a3*(1-m)*xt^m+b2/2*xt^2 + a0
    }
    else{                  ### When a3>0, we apply the TLF method
      A2 = a3*m*(m-1)/2*xt^(m-2)+a2
      A1 = a3*m*(2-m)*xt^(m-1)+a1
      A0 = a3*(m-1)*(m/2-1)*xt^m+a0
    }
    Delta.det = A1^2-4*A2*A0
    if(Delta.det >= 0){
      xt.save[k+1] = -(A1+sqrt(Delta.det))/(2*A2)
      yt.save[k+1] = a3*xt.save[k+1]^m + a2*xt.save[k+1]^2 +
        a1*xt.save[k+1]+a0
      esp = abs(yt.save[k+1])
      xt = xt.save[k+1]
      y1t = yt.save[k+1]
      k = k+1
    }
    else{
      stop(paste('There exists no root in (0,',a.max,')'))
    }
  }
  return(xt)
}
# --------------------------------------------------------------------------
# --------Example: x^3-3x^2+x+1 = 0-----------------------------------------
> a3 = 1; a2 = -3; a1 = 1; a0 = 1
> root.HPE(a3, a2, a1, a0, a.max=2, m=3, x0=0)
        [,1]      [,2]       [,3]         [,4] [,5]
xt.save    0 0.7675919 0.99418291 9.999999e-01    1
yt.save    1 0.4522631 0.01163398 1.968315e-07    0
\end{verbatim}

\vkU
\baselineskip 0.24in
\renewcommand{\theequation}{\arabic{section}.\arabic{equation}}
\setcounter{equation}{0}
\renewcommand{\theequation}{E.\arabic{equation}}
\renewcommand{\thesubsection}{E.\arabic{subsection}}

\section*{Appendix E: \  \textsf{R} codes with analytical solutions \\ \hspace*{4.0cm} for solving the roots of a cubic equation}

Let $a_3 \ne 0$. Solving the roots of the following cubic equation
\begin{eqnarray} \label{uclaoneE.1}
   a_3 x^3 + a_2 x^2 + a_1 x + a_0 = 0
\end{eqnarray}
is equivalent to finding the roots of the cubic equation without quadratic term
$$
   y^3 + p y + q = 0,
$$
where $y = x + a_2/(3a_3)$,
\begin{eqnarray}  \label{uclaoneE.2}
   p = \frac{3 a_3 a_ 1 - a_2^2}{3a_3^2} \qand q = \frac{2a_2^3 - 9a_3a_2a_1 + 27 a_3^2 a_0}{27 a_3^3}.
\end{eqnarray}

Define
$$
   \de =  \left(\frac{q}{2}\right)^2 + \left(\frac{p}{3}\right)^3.
$$
The discriminant rule is as follows:
\begin{namelist}{0123}
\item[\hspace*{-0.03cm} (a)] If $\de >0$, then there exists only one real root $y_1$.

\item[\hspace*{-0.03cm} (b)] If $\de =0$, then there exist three multiple roots $y_1=y_2=y_3$, so that when $\de \ge 0$, the real roots of the equation (\ref{uclaoneE.2}) and the equation (\ref{uclaoneE.1}) are
    \begin{eqnarray*}
       y_1 = u_1 + u_2 \qand x_1 = y_1 - \frac{a_2}{3a_3},
    \end{eqnarray*}
    respectively, where
    $$
       u_1 = \left[ -q/2 + \sqrt{\de } \right]^{1/3} \qand u_2 = \left[ -q/2 - \sqrt{\de } \right]^{1/3}.
    $$

\item[\hspace*{-0.03cm} (c)] If $\de <0$, then there exist three different real roots of the equation (\ref{uclaoneE.2}) and the equation (\ref{uclaoneE.1}), which can be expressed by
    \begin{eqnarray*}
       y_1 &=& 2r^{1/3} \cos \tha ,\quad y_2 \;=\; 2r^{1/3} \cos \left[\tha+ 2/(3\pi)\right], \quad y_3 \;=\; 2r^{1/3} \cos \left[\tha+ 4/(3\pi)\right], \\
       x_1 &=& y_1 - a_2/(3a_3) ,\quad x_2 \;=\; y_2 - a_2/(3a_3) \qand x_3 \;=\; y_3 - a_2/(3a_3),
    \end{eqnarray*}
    where
    $$
       r = \sqrt{ -\left(p/3\right)^3} \qand \tha = \frac{\arccos \left(-q/2r\right)} {3}.
    $$
\end{namelist}

\vkL The \textsf{R} codes for solving the cubic equation (\ref{uclaoneE.1}) with explicit solutions is provided as follows:
\begin{verbatim}
root_ploynomial_3 = function(a3, a2, a1, a0)
{ # Function name:(a3, a2, a1, a0)#########################################
  # ---------------------------Aim-----------------------------------------
  # Calculating the equation: a3*x^3+a2*x^2+a1*x+a0, -inf < x < inf
  # ---------------------------Input---------------------------------------
  #    a3: the coefficient of x^3
  #    a2: the coefficient of x^2
  #    a1: the coefficient of x
  #    a0: the constant, that is a positive real number
  # ---------------------------Output---------------------------------------
  # The real roots
  ##########################################################################
  nroot <- function(x,n)
  { # Function for calculating x^(1/n), n = 1,3,5,...
    abs(x)^(1/n)*sign(x)
  }
  p = (3*a3*a1-a2^2)/(3*a3^2)
  q = (2*a2^3-9*a3*a2*a1+27*a3^2*a0)/(27*a3^3)
  disterm = (q/2)^2+(p/3)^3
  if(disterm>=0){    ## When disterm>=0, there exists one real root
    u1 = nroot(-q/2+sqrt(disterm),3)
    u2 = nroot(-q/2-sqrt(disterm),3)
    t1 = u1+u2
    x1 = t1 - a2/(3*a3)
    return(x1)
  }else{            ## When disterm<0, there exist three real roots
    r = sqrt(-(p/3)^3)
    theta = 1/3*acos(-q/(2*r))
    t1 = 2*nroot(r,3)*cos(theta)
    t2 = 2*nroot(r,3)*cos(theta+2/3*pi)
    t3 = 2*nroot(r,3)*cos(theta+4/3*pi)
    x1 = t1 - a2/(3*a3)
    x2 = t2 - a2/(3*a3)
    x3 = t3 - a2/(3*a3)
    return(c(x1,x2,x3))
  }
}
root_ploynomial_3(a3,a2,a1,a0)
# --------------------------------------------------------------------------
# --------Example: x^3-3x^2+x+1 = 0-----------------------------------------
> a3 = 1; a2 = -3; a1 = 1; a0 = 1
> root_ploynomial_3(1, -3, 1, 1)
[1]  2.4142136 -0.4142136  1.0000000
\end{verbatim}

\newpage
\baselineskip 0.30in
\renewcommand{\theequation}{\arabic{section}.\arabic{equation}}
\setcounter{equation}{0}
\renewcommand{\theequation}{1.\arabic{equation}}
\renewcommand{\thesubsection}{1.\arabic{subsection}}

\section*{Supplementary material: Other more applications of the \\ \hspace*{7.4cm} US algorithm in statistics}  

\subsection{Calculation of exact $p$-values for skew null distributions} 

\subsubsection{The chi-squared null distribution}  

Let $\{X_i\}_{i=1}^n \iid N(\mu, \si^2)$ with unknown $\{\mu, \si^2\}$ and we want to test $H_0$: $\si^2 =\si_0^2$ against $H_1$: $\si^2 \ne \si_0^2$. The test statistic is $\chi =\nu S^2/\si_0^2$ with corresponding observed value $\chi_{\rm obs} \teq \nu s^2/\si_0^2$, where $S^2 = (1/\nu) \sum_{i=1}^n (X_i - \BX)^2$ is the sample variance, $s^2$ is its observed value and $\nu \heq n-1$. Under $H_0$, we have $\chi \sim \chi^2(\nu)$, which is called chi-squared null distribution. We denote the pdf of $\chi^2(\nu)$ by
$h_\nu(x) = [2^{\nu/2} \Ga(\nu/2)]^{-1}  x^{\nu/2 - 1} \e^{-x/2}$ for $x>0$, whose mode is $\nu-2 \teq m_0$.

\vkL\noi\textbf{(a) Case I: $\chi_{\rm obs} < m_0$}\vkL\noi
The exact $p$-value can be calculated as $p\mbox{-value} =  \Pr\{\chi^2(\nu) \le \chi_{\rm obs}\} + \Pr\{\chi^2(\nu) \ge x_u\}$, where $x_u \,(> m_0)$ satisfies $h_{\nu}(x_u) = h_{\nu}(\chi_{\rm obs})$; i.e., $m_0 \log (x_u) - x_u = m_0\log (\chi_{\rm obs}) - \chi_{\rm obs} \;\teq\; c_0$. Thus, finding $x_u$ is equivalent to finding the root to the equation $g(x_u) =  m_0 \log (x_u) - x_u - c_0 = 0$. Since $g'(x_u) = m_0/x_u - 1 > -1 \;\teq\; b(x_u)$, the US iteration for calculating $x_u$ is
\begin{eqnarray*}
   x_u^{(t+1)} &\sr{(\ref{uclaone3.3})} & \sol \left\{ g(\xt_u) + \tex\int_{\xt_u}^{x_u} b(z) \rd z = 0, \; \forall \; x_u, \xt_u > m_0 \right\} \\ [2mm]
   &=& \sol \Big\{ g(\xt_u) - x_u + \xt_u =0, \; \forall \; x_u, \xt_u > m_0  \Big\} \\ [2mm]
   &=& g(\xt_u) + \xt_u = m_0 \log \big(\xt_u/\chi_{\rm obs} \big) + \chi_{\rm obs}.
\end{eqnarray*}

\vkl\noi\textbf{(b) Case \II: $\chi_{\rm obs} > m_0$}\vkL\noi
The exact $p$-value can be calculated as $p\mbox{-value} =  \Pr\{\chi^2(\nu) \le x_l\} + \Pr\{\chi^2(\nu) \ge \chi_{\rm obs}\}$, where $x_l \in (0, m_0)$ satisfies $h_{\nu}(x_l) = h_{\nu}(\chi_{\rm obs})$; i.e., $m_0 \log (x_l) - x_l = m_0 \log (\chi_{\rm obs}) - \chi_{\rm obs} \;\teq\; c_0$. Thus, finding $x_l$ is equivalent to finding the root to the equation $g(x_l) = x_l - m_0 \log (x_l) + c_0 = 0$. Since $g'(x_l) = 1 - m_0/x_l > -m_0/x_l \;\teq\; b(x_l)$, the US iteration for calculating $x_l$ is
\begin{eqnarray*}
   x_l^{(t+1)} &\sr{(\ref{uclaone3.3})} & \sol \left\{ g(\xt_l) + \tex\int_{\xt_l}^{x_l} b(z) \rd z = 0, \; \forall \; x_l, \xt_l  \in (0, m_0) \right\} \\ [2mm]
   &=& \sol \Big\{ g(\xt_l) -m_0 \log (x_l) + m_0 \log (\xt_l) =0, \; \forall \; x_l, \xt_l  \in (0, m_0) \Big\} \\ [2mm]
   &=& \xt_l \exp \left[ g(\xt_l)/m_0 \right].
\end{eqnarray*}

\subsubsection{The \textit{F} null distribution}  

Let $\{X_{ij}\}_{j=1}^{n_i} \iid N(\mu_i, \si_i^2)$ for $i=1, 2$, and the two random samples be independent, where $\{\mu_i, \si_i^2\}$ are unknown. We want to test $H_0$: $\si_1^2 =\si_2^2$ against $H_1$: $\si_1^2 \ne \si_2^2$, and the test statistic is $F =S_1^2/S_2^2$ with corresponding observed value $F_{\rm obs} \teq s_1^2/s_2^2$, where $S_i^2$ is the $i$-th sample variance, $s_i^2$ is its observed value and $\nu_i \teq n_i-1$. Under $H_0$, we have $F \sim F(\nu_1, \nu_2)$, which is called $F$ null distribution. We denote the pdf of $F(\nu_1, \nu_2)$ by
$$
   h_{\nu_1, \nu_2}(x) = \frac{\Ga(\nu_{12}/2) \nu^{\nu_1/2} } {\Ga(\nu_1/2) \Ga(\nu_2/2)} x^{\nu_1/2 -1 } (1 + \nu x )^{- \nu_{12}/2}, \quad x>0, \;\; \nu \teq \frac{\nu_1}{\nu_2}, \;\; \nu_{12} \teq \nu_1+\nu_2.
$$
The mode of $h_{\nu_1, \nu_2}(x)$ is $m_1\nu_2 /[\nu_1(\nu_2+2)] \teq m_0$, where $m_1 \teq \nu_1-2$.

\vkL\noi\textbf{(a) Case I: $F_{\rm obs} < m_0$}\vkL\noi
The exact $p$-value can be calculated as $p\mbox{-value} =  \Pr\{F(\nu_1, \nu_2) \le F_{\rm obs}\} + \Pr\{F(\nu_1, \nu_2) \ge x_u\}$, where $x_u \,(> m_0)$ satisfies $h_{\nu_1, \nu_2}(x_u) = h_{\nu_1, \nu_2}(F_{\rm obs})$; i.e.,
$m_1 \log (x_u) - \nu_{12} \log (1+ \nu x_u) = m_1 \log (F_{\rm obs}) - \nu_{12} \log (1+ \nu F_{\rm obs}) \teq c_0$. Thus, finding $x_u$ is equivalent to finding the root to the equation $g(x_u) = m_1 \log (x_u) - \nu_{12} \log (1+ \nu x_u) - c_0 = 0$. Since $g'(x_u) = m_1/x_u - \nu_{12} \nu/(1 + \nu x_u) > m_1/x_u - \nu_{12}/x_u = -(2+\nu_2)/x_u \;\teq\; b(x_u)$, the US iteration for finding $x_u$ is
\begin{eqnarray*}
   x_u^{(t+1)} &\sr{(\ref{uclaone3.3})} & \sol \left\{ g(\xt_u) + \tex\int_{\xt_u}^{x_u} b(z) \rd z = 0, \; \forall \; x_u, \xt_u > m_0 \right\} \\ [2mm]
   &=& \sol \Big\{ g(\xt_u) -(2+\nu_2) \log (x_u / \xt_u) =0, \; \forall \; x_u, \xt_u > m_0 \Big\} 
\end{eqnarray*}
\begin{eqnarray*}
   &=& \xt_u \exp [ g(\xt_u)/(2+\nu_2)].
\end{eqnarray*}

\vkl\noi\textbf{(b) Case \II: $F_{\rm obs} > m_0$}\vkL\noi
The exact $p$-value can be calculated as $p\mbox{-value} =  \Pr\{F(\nu_1, \nu_2) \le x_l\} + \Pr\{F(\nu_1, \nu_2) \ge F_{\rm obs}\}$, where $x_l \in (0, m_0)$ satisfies $h_{\nu_1, \nu_2}(x_l) = h_{\nu_1, \nu_2}(F_{\rm obs})$; i.e.,
$m_1 \log (x_l) - \nu_{12} \log (1+ \nu x_l) = m_1 \log (F_{\rm obs}) - \nu_{12} \log (1+ \nu F_{\rm obs}) \teq c_0$. Thus, finding $x_l$ is equivalent to finding the root to the equation $g(x_l) =  c_0 - m_1 \log (x_l) + \nu_{12} \log (1+ \nu x_l) = 0$. Since $g'(x_l) = -m_1/x_l + \nu_{12} \nu/(1 + \nu x_l) > -m_1/x_l \;\teq\; b(x_l)$, the US iteration for calculating $x_l$ is
\begin{eqnarray*}
   x_l^{(t+1)} &\sr{(\ref{uclaone3.3})} & \sol \left\{ g(\xt_l) + \tex\int_{\xt_l}^{x_l} b(z) \rd z = 0, \; \forall \; x_l, \xt_l \in (0, m_0) \right\} \\ [2mm]
   &=& \sol \Big\{ g(\xt_l) -m_1 \log (x_l) + m_1 \log (\xt_l) =0, \; \forall \; x_l, \xt_l \in (0, m_0) \Big\} \\ [2mm]
   &=& \xt_l \exp [ g(\xt_l)/m_1 ].
\end{eqnarray*}

\subsection{MLEs of parameters in a class of continuous distributions}

\subsubsection{MLE of $\al$ in gamma distribution} 

Let $\{X_i\}_{i=1}^n \iid \mGamma(\al,\be)$ with pdf $\be^{\al} x^{\al-1}\e^{-\be x} /\Ga(\al)$, $x>0$, where $\{\al, \be\}$ are two positive parameters. The log-likelihood function $\ell(\al,\be) = nG(\ibx)(\al-1) - n \Bx \be + n \left[ \al \log \be - \log \Ga(\al) \right]$,
where $G(\ibx) = G(x_1, \ldots, x_n) \teq (1/n)\sum_{i=1}^n \log (x_i)$ and $\Bx = (1/n)\sum_{i=1}^n  x_i$. Let
\begin{eqnarray}
   0 = \frac{\pa \ell(\al,\be)}{\pa \al} &=& nG(\ibx) + n \log (\be) - n \psi(\al), \label{uclaone4.8} \\ [2mm]
   0 = \frac{\pa \ell(\al,\be)}{\pa \be} &=& - n\Bx + \frac{n\al}{\be}. \label{uclaone4.9}
\end{eqnarray}
Given $\al$, from (\ref{uclaone4.9}), we have $\be = \al/\Bx$. Given $\be$,  from (\ref{uclaone4.8}), the MLE of $\al$ is the root of the equation: $g(\al) = c_0 - \psi(\al) = 0$ for $\al>0$, where $c_0 \teq G(\ibx) + \log (\be)$ is a known constant, $\psi(\cdot)$ is the digamma function defined as
\begin{eqnarray} \label{uclaone4.10}
   \psi(\al) \teq \frac{\Ga'(\al)}{\Ga(\al)} = -\ga +\sum_{m=0}^{\infty}\left(\frac{1}{m+1}-\frac{1}{m+\al}\right),
\end{eqnarray}
here $\ga = \lim_{n\to \infty} \left(\sum_{m=1}^n m^{-1} - \log n  \right) \approx 0.5772$ is the Euler--Mascheroni constant. Since
\begin{eqnarray*}
   g'(\al) &=& - \psi'(\al) = -\sum_{m=0}^{\infty} \frac{1}{(m+\al)^2} = -\frac{1}{\al^2} - \sum_{m=1}^{\infty} \frac{1}{(m+\al)^2} \\ [2mm]
   & > & -\frac{1}{\al^2} - \sum_{m=1}^{\infty} \frac{1}{m^2} = -\frac{1}{\al^2} - \frac{\pi^2}{6} \;\teq \; b(\al),
\end{eqnarray*}
the US iteration for calculating the MLE $\hal$ is
\begin{eqnarray*}
   \altI 
   &=& \sol \left\{ g(\alt) +\frac{1}{\al} - \frac{\pi^2}{6}\al - \frac{1}{\alt}+\frac{\pi^2}{6}\alt = 0, \; \forall \; \al, \alt >0  \right\} \\ [2mm]
   &=& \sol \left\{ (\pi^2/6)\al^2 - a_3^{(t)} \al - 1=0, \; \forall \; \al, \alt >0  \right\},
\end{eqnarray*}
where $a_3^{(t)} = c_0 -\psi(\alt) + (\pi^2/6)\alt - 1/\alt$.

By replacing $\be$ with $\al/\Bx$, the profile log-likelihood function of $\al$ is
$\ell(\al) =  n \al \log \al + nG(\ibx)(\al-1) - n \al - n\al \log \Bx - n\log \Ga(\al)$. Let
$$
   0 \;=\; \frac{1}{n} \ell'(\al) \;=\; \log \al + G(\ibx) - \log \Bx - \psi(\al) \;\teq\; g_1(\al) + g_2(\al) \;\teq\; g(\al),
$$
where $g_1(\al) \teq \log \al$ and $g_2(\al) \teq G(\ibx) - \log \Bx - \psi(\al)$.
Since
\begin{eqnarray*}
   g_1'(\al) &=& \al^{-1}, \quad g_1''(\al) \;=\; -\al^{-2}, \quad g_1'''(\al) \;=\; 2\al^{-3} \;>\; 0 \;\teq\; b_3  \\
   g_2'(\al) &=& - \psi'(\al), \quad g_2''(\al) \;=\; - \psi''(\al) \;>\;0 \;\teq\; b_{21},
\end{eqnarray*}
we can construct a $U$-fnction for $g_1(\al)$ by using TLB method and construct a $U$-fnction for $g_2(\al)$ by using SLB method, namely,
\begin{eqnarray*}
   U_1(\al|\alt) &\sr{(\ref{uclaone3.15})}& \log \alt + \frac{1}{\alt}(\al-\alt) - \frac{1}{2\al^{2(t)}}(\al-\alt)^2 \qand \\
   U_2(\al|\alt) &\sr{(\ref{uclaone3.10})}& -\left[\psi(\alt) + \psi'(\alt)(\al-\alt) \right] I(\al \ge \alt) - \psi(\al) I(\al < \alt).
\end{eqnarray*}
Therefore, the $U$-function for $g(\al)$ is
\begin{eqnarray*}
   U(\al|\alt) &\teq& U_1(\al|\alt) + U_2(\al|\alt) \\
   &=& \log \alt + \frac{1}{\alt}(\al-\alt) - \psi(\al) I(\al < \alt) - \frac{1}{2\al^{2(t)}}(\al-\alt)^2 \\[2mm]
   & &-\;\left[\psi(\alt) + \psi'(\alt)(\al-\alt) \right] I(\al \ge \alt).
\end{eqnarray*}
Given initial value $\al^{(0)}$ such that $g(\al^{(0)})\ge 0$, the US iteration for calculating the MLE $\hal$ is
\begin{eqnarray*}
   \altI
   &=& \sol \left\{ g(\alt) + g'(\alt)(\al-\alt) - \frac{1}{2\al^{2(t)}}(\al-\alt)^2  = 0, \; \forall \; 0<\al, \alt <\hal \right\} \\[2mm]
   &=& \alt + \frac{g'(\alt) + \sqrt{[g'(\alt)]^2 + 2g(\alt)/ \al^{2(t)}}}{\al^{2(t)}}.
\end{eqnarray*}

\subsubsection{MLE of $\tha$ in Weibull distribution} 

Let $\{X_i\}_{i=1}^n \iid \mbox{Weibull}(\tha, \la)$ with pdf $(\tha/\la) (x/\la)^{\tha-1} \exp [-(x/\la)^{\tha} ]$, $x>0$, where $\{\tha, \la\}$ are two positive parameters (Hallinan 1993, Rinne 2008). The log-likelihood function is $\ell(\tha, \la) = n \log (\tha)- n\tha \log(\la) + n G(\ibx) (\tha-1) - \sum_{i=1}^n (x_i/\la)^{\tha}$, where $G(\ibx) \teq (1/n)\sum_{i=1}^n \log (x_i)$. Let $0 = \pa \ell(\tha,\la)/\pa \la = - n\tha\la^{-1} + \tha\la^{-\tha-1} \sum_{i=1}^n x_i^{\tha}$. Given $\tha$, we obtain $\la = [(1/n)\sum_{i=1}^n x_i^{\tha}]^{1/\tha}$. Replacing $\la$ in $\ell(\tha, \la)$ by $[(1/n)\sum_{i=1}^n x_i^{\tha}]^{1/\tha}$, we know that $\ell(\tha, \la)$ will reduce to $\ell(\tha) = \tha nG(\ibx) + n \log(\tha)  - n \log (\tex\sum_{i=1}^n  x_i^{\tha} ) + \mbox{ constant}$. Thus, the MLE of $\tha$ is the root of $\ell'(\tha) =0$ or the root of the equation
$$
   g(\tha) = G(\ibx) + \frac{1}{\tha}- \frac{\sum_{i=1}^n  x_i^{\tha} \log(x_i)}{\sum_{i=1}^n  x_i^{\tha}} = 0, \quad \tha >0.
$$
Since
\begin{eqnarray*}
  g'(\tha) &=& - \frac{1}{\tha^2} - \frac{[\sum_{i=1}^n x_i^{\tha} (\log x_i)^2 ]  (\sum_{i=1}^n  x_i^{\tha}) - (\sum_{i=1}^n x_i^{\tha} \log x_i )^2 } {(\sum_{i=1}^n  x_i^{\tha} )^2} \\ [2mm]
  &\ge& - \frac{1}{\tha^2} - \frac{\sum_{i=1}^n x_i^{\tha} (\log x_i)^2 } {\sum_{i=1}^n  x_i^{\tha}}  \\[2mm]
  &\ge& - \frac{1}{\tha^2} - \max_{1\le i \le n} (\log x_i)^2 \; \teq \; - \frac{1}{\tha^2} - T_{\max} \; \teq \;  b(\tha),
\end{eqnarray*}
the US iteration for calculating the MLE $\hth$ is
\begin{eqnarray*}
   \thtI 
   &=& \sol \left\{ g(\tht) + \frac{1}{\tha} - T_{\max} \tha - \frac{1}{\tht} + T_{\max} \tht = 0, \; \forall \; \tha, \tht >0  \right\} \\ [2mm]
   &=& \sol \left\{ T_{\max} \tha^2 - a_4^{(t)} \tha - 1=0, \; \forall \; \tha, \tht >0  \right\},
\end{eqnarray*}
where $a_4^{(t)} \teq g(\tht) + T_{\max} \tht -1/\tht$.

\subsubsection{MLE of $\nu$ in Student's $\ibt$-distribution}  

Let $\{X_i\}_{i=1}^n \iid t(\mu, \si^2, \nu)$ with pdf $\Ga((\nu+1)/2)/[\Ga(\nu/2)\sqrt{\nu\pi}\si] [1+(x-\mu)^2/(\nu\si^2) ]^{-(\nu+1)/2}$ for $x\in \bbR$, where $\mu\in \bbR$, $\si\,(>0), \nu\, (>0)$ are parameters. The log-likelihood function is
$$
   \ell(\mu, \si^2, \nu) = n\left[\log \Ga\left(\frac{\nu+1}{2}\right) - \log \Ga\left(\frac{\nu}{2}\right) - \frac{\log \nu +\si^2}{2}\right] - \frac{\nu+1}{2} \sum_{i=1}^n \log\left[1+\frac{(x_i-\mu)^2}{\nu\si^2} \right].
$$
Although EM/MM algorithms can be employed to calculate the MLEs $\{\hmu, \hsi^2\}$ for a given $\nu$, there is no efficient algorithm to calculate $\hnu$ given $\{\mu,\si^2\}$. The aim of this subsection is to compute $\hnu$ by fixing $\{\mu, \si^2\}$, which is the root of $0 = \ell'(\nu|\mu, \si^2)$; i.e.,
$$
   0 = g(\nu) = \frac{n}{2}\left[ \psi\left(\frac{\nu+1}{2}\right) - \psi \left(\frac{\nu}{2}\right) - \frac{1}{\nu}\right] - \frac{1}{2} \sum_{i=1}^n \log\left(1+\frac{z_i}{\nu} \right) + \frac{\nu+1}{2\nu^2} \sum_{i=1}^n \tau_i(\nu),
$$
where $z_i \teq (x_i-\mu)^2/\si^2$ and $\tau_i(\nu) \teq z_i (1+ z_i\nu^{-1} )^{-1}$. Furthermore, by defining
$$
   \de(\nu) \teq \frac{n}{4}\left[ \psi'\left(\frac{\nu+1}{2}\right) - \psi' \left(\frac{\nu}{2}\right) + \frac{2}{\nu^2}\right] \qand d_i(\nu) \teq \frac{\tau_i(\nu)}{\nu} = \frac{z_i}{\nu} \left(1+\frac{z_i}{\nu} \right)^{-1} \; \in\; (0,1),
$$
we obtain
\begin{eqnarray}
   g'(\nu) &=& \de(\nu) + \left(\frac{1}{2\nu^2} - \frac{\nu + 2}{2\nu^3} \right) \sum_{i=1}^n \tau_i(\nu) + \frac{\nu+1}{2\nu^4} \sum_{i=1}^n \tau_i^2(\nu)  \non \\ [2mm]
   &=& \de(\nu) - \frac{1}{\nu^2} \sum_{i=1}^n d_i(\nu)
   + \frac{\nu+1}{2\nu^2} \sum_{i=1}^n  d_i^2(\nu)
   = \de(\nu) + \frac{1}{2\nu^2} \sum_{i=1}^n  \left[ (\nu+1)d_i^2(\nu) - 2d_i(\nu) \right] \non \\ [2mm]
   &=& \de(\nu) + \frac{1}{2\nu^2} \left\{ \sum_{i=1}^n \frac{ [ (\nu+1)d_i(\nu) - 1]^2 }{\nu+1} -\frac{n}{\nu+1} \right\}
   \; \ge \; \de(\nu) - \frac{n}{2\nu^2 (\nu+1)}   \label{uclaone4.11} \\ [2mm]
   &=& \frac{n}{4} \left[ \sum_{m=0}^{\infty} \left(m+ \frac{\nu}{2}+\frac{1}{2}\right)^{-2} - \sum_{m=0}^{\infty} \left(m+\frac{\nu}{2}\right)^{-2}  + \frac{2}{\nu^2} \right] - \frac{n}{2\nu^2 (\nu+1)}  \non \\ [2mm]
   &=& \frac{n}{4} \left[ \sum_{m=0}^{\infty} \left(m+ \frac{1}{2}+\frac{\nu}{2}\right)^{-2} - \sum_{m=0}^{\infty} \left(m+1+\frac{\nu}{2}\right)^{-2} - \frac{4}{\nu^2} + \frac{2}{\nu^2} \right] - \frac{n}{2\nu^2} \frac{1}{(\nu+1)} \non \\
   &>& - \frac{n}{2\nu^2} - \frac{n}{2\nu^2}  \;=\; - \frac{n}{\nu^2} \; \teq \; b(\nu),   \non
\end{eqnarray}
where the inequality (\ref{uclaone4.11}) holds iff $d_i(\nu) = (\nu+1)^{-1}$.
The US iteration is
\begin{eqnarray*}
   \nu^{(t+1)} = \sol \Big\{ \left[g(\nu^{(t)}) - n/\nu^{(t)} \right]\nu + n =0, \; \forall\, \nu, \nu^{(t)} >0 \Big\} =  \frac{n \nu^{(t)}} {n  - \nu^{(t)} g(\nu^{(t)})}.
\end{eqnarray*}

\subsubsection{MLE of $\al$ in unit power--logarithmic distribution}  

The unit power--logarithmic distribution is used to model extremely left--skewed data because of a large panel of J shapes in its pdf (Chesneau 2021). Let $\{X_i\}_{i=1}^n \iid \mbox{UPL}(\al)$ with pdf $\log^{-1}(\al + 1) \cdot (1 - x^{\al}) [-\log(x)]^{-1}$ for $x\in (0, 1)$, where $\al \,(>0)$ is the parameter. The log-likelihood function is $\ell(\al) = -n \log [\log(\al+1)] + \sum_{i=1}^n \log(1-x_i^{\al}) - \sum_{i=1}^n \log [-\log(x_i)]$.  Thus, the MLE of $\al$ is the root of $0 = \ell'(\al)$ or the root of the equation
$$
   0 = g(\al) = - \frac{n}{(\al+1)\log(\al+1)} - \sum_{i=1}^n \frac{x_i^{\al} \log x_i}{1-x_i^{\al}}  = u_0(\al) - \sum_{i=1}^n \left[ \frac{ \log x_i}{x_{(1)}^{-\al}-1} \cdot u_i(\al) \right],
$$
where
\begin{eqnarray}
   u_0(\al) &\teq & - \frac{n}{(\al+1)\log(\al+1)}, \quad x_{(1)} \teq \min (x_1, \ldots, x_n),  \non \\ [2mm]
   u_i(\al) &\teq & \frac{x_{(1)}^{-\al}-1}{x_i^{-\al}-1} = \frac{\e^{-\al\log x_{(1)} }-1} {\left(\e^{-\al \log x_{(1)}}\right)^{\log x_i / \log x_{(1)}} -1} = \frac{y-1} {y^{t_i} -1}, \quad i=1, \ldots, n, \non \\ [2mm]
   y \teq y(\al) &\teq& \exp \left[ \al (-\log x_{(1)}) \right] > 1, \label{uclaone4.12} \\ [2mm]
   t_i &\teq& \frac{\log x_i} {\log x_{(1)}} \in (0,1]. \label{uclaone4.13}
\end{eqnarray}

If all $\{u_i(\al)\}_{i=0}^n$ are increasing with respect to $\al$, we can apply the fixed--block method to construct a $U$-function for $g(\al)$ at $\al = \alt$. First, since
$$
   u_0'(\al) =  \frac{n[\log(\al+1) +1]}{(\al+1)^2 \log^2 (\al+1)} >0,
$$
we know that $u_0(\al)$ is a strictly increasing function of $\al$. Next, since $y(\al)$ is an increasing function about $\al$, we only need to prove that $u_{1i}(y) \teq u_i(\al)  = (y-1)/(y ^{t_i} -1)$, $y > 1$, is increasing with respect to $y$.

Set $h_i(y) = (1-t_i)y ^{t_i} -1 + t_i y ^{t_i-1}$ for $y>1$. From (\ref{uclaone4.13}) and (\ref{uclaone4.12}), we have $h_i'(y) = t_i (1-t_i) y ^{t_i-2}(y-1) \ge  0$, so that
\begin{eqnarray}  \label{uclaone4.14}
   h_i(y) \ge h_i(1)=0.
\end{eqnarray}
Note that
$$
   u'_{1i}(y)  = \frac{y ^{t_i} -1 - t_i y ^{t_i-1}(y-1)}{(y ^{t_i} -1)^2} = \frac{h_i(y)}{(y ^{t_i} -1)^2} \sge{(\ref{uclaone4.14})}  0,
$$
implying that $u_{1i}(y)$ is increasing about $y$. By applying the fixed--block method, we obtain
\begin{eqnarray*}
   g(\al) &\sge{\sgn(\al-\alt)}&  u_0(\alt) - \sum_{i=1}^n  \left[ \frac{\log x_i}{x_{(1)}^{-\al}-1} \cdot u_i(\alt)  \right] \\ [2mm]
   &=& u_0(\alt) -  \frac{1}{x_{(1)}^{-\al}-1} \sum_{i=1}^n  u_i(\alt) \log (x_i) = U(\al|\alt),
\end{eqnarray*}
so the US iteration is given by
$$
   \altI = \sol \left\{  x_{(1)}^{-\al}= 1 - \frac{(\alt+1) \log(\alt+1)}{n} \sum_{i=1}^n  \frac{(x_{(1)}^{-\alt}-1) \log (x_i)} {x_i^{-\alt}-1}, \; \forall\, \al, \alt >0 \right\}.
$$

\subsection{MLEs of parameters in a class of discrete distributions}

\subsubsection{MLE of $\tha$ in zeta distribution} 

The zeta (or Zipf) distribution can be used to model the frequency of occurrence of a word randomly chosen from a text, or the population rank of a city randomly chosen from a country. Let $\{X_i\}_{i=1}^n\iid \Zeta(\tha)$ with pmf $Z^{-1}(\tha) x^{-\tha}$, $x=1, 2,\ldots, \infty$, where $\tha\, (>1)$ is a parameter and $Z(\tha) = \sum_{x=1}^{\infty} x^{-\tha}$ is called Riemann's zeta function (Nair \Et 2018, p.131--134). The log-likelihood is $\ell(\tha) = -  nG(\ibx) \tha  - n \log Z(\tha)$, where $G(\ibx) \teq (1/n)\sum_{i=1}^n \log (x_i)$. Let $0=\ell'(\tha) = -  nG(\ibx) - n Z'(\tha)/Z(\tha)$, then the MLE of $\tha$ is the root of the equation: $g(\tha) = -  G(\ibx) - Z'(\tha)/Z(\tha) = 0$ for $\tha >1$. Note that
$$
   g'(\tha) = -\;\frac{Z''(\tha)Z(\tha) - [Z'(\tha)]^2}{Z^2(\tha)} \ge -\; \frac{Z''(\tha) }{Z(\tha)}
  \sge{(\rm \ref{uclaoneC.3}) \; \& \; (\rm \ref{uclaoneC.4})}  a_5 - \frac{a_6} {(\tha-1)^2} \; \teq \; b(\tha),
$$
where $a_5 \teq - \log 2 \cdot (\log 2 + 2)$ and $a_6 \teq  2 \log 2+ 2$, then
the US iteration is
\begin{eqnarray*}
   \thtI 
   &=& \sol \left\{ g(\tht) + a_5(\tha-\tht) + \frac{a_6}{\tha-1} - \frac{a_6}{\tht-1} = 0, \; \forall \; \tha, \tht > 1  \right\} \\ [2mm]
   &=& \sol \left\{ a_5  \tha^2 + a_7^{(t)} \tha + a_8^{(t)}=0, \; \forall \; \tha, \tht > 1  \right\},
\end{eqnarray*}
where $a_7^{(t)} \teq g(\tht) - a_5 (\tht+1)- a_6/(\tht-1)$ and $a_8^{(t)} \teq -g(\tht) + a_5 \tht + a_6 \tht /(\tht-1)$.

\subsubsection{MLE of $\al$ in gamma--Poisson distribution} 

Let $\{X_i\}_{i=1}^n\iid \mbox{GPoisson}(\al,\be)$ with pmf $\Ga(x+\al)\Ga^{-1}(\al) \be^\al (\be + 1)^{-x-\al}/x!$, $x=0,1,2,\ldots,\infty$, where $\al$ and $\be$ are two positive parameters (Arbous \& Kerrich 1951, Bates \& Neyman 1952). Define $\Bx =(1/n)\sum_{i=1}^n x_i$, the log-likelihood function is
$$
   \ell(\al, \be)= \sum_{i=1}^n \log \Ga(x_i+\al) -n \log \Ga(\al) + n \al \log (\be)-n(\al+\Bx) \log (\be+1) + \mbox{constant}.
$$
Let $0 = \pa \ell(\al,\be)/\pa \be = n \al/\be-n(\al+\Bx)/(\be+1)$, we obtain $\be = \al/\Bx$. Replacing $\be$ in $\ell(\al, \be)$ by $\al/\Bx$, we know that $\ell(\al, \be)$ becomes
$$
   \ell(\al) = \sum_{i=1}^n  \log \Ga(x_i+\al)-n \log \Ga(\al) + n \al \log \al -\ n(\al+\Bx) \log (\al+\Bx) + \mbox{constant}.
$$
Thus, the MLE of $\al$ is the root of the following equation:
\begin{eqnarray*}
   0 &=& \ell'(\al) = \sum_{i=1}^n  \left[\psi(x_i+\al)- \psi(\al)\right] + n \log \al - n \log (\al+\Bx) \\ [2mm]
   &\sr{(\ref{uclaone4.10})}& \sum_{i=1}^n  \left[  -\sum_{m=0}^{\infty}  \frac{1}{m+x_i+\al} + \sum_{m=0}^{\infty} \frac{1}{m+\al} \right]  + n \log \al - n \log (\al+\Bx)  \\ [2mm]
   &=& \sum_{i=1}^n  \left[  - \sum_{m=x_i}^{\infty}  \frac{1}{m+\al} + \sum_{m=0}^{\infty} \frac{1}{m+\al} \right]  + n \log \al - n \log (\al+\Bx) \\ [2mm]
   &=& \sum_{i=1}^n  \left[ I(x_i \ge 1) \sum_{m=0}^{x_i-1}\frac{1}{m+\al} \right]  + n \log \al - n \log (\al+\Bx) \teq  g(\al).
\end{eqnarray*}
Since
\begin{eqnarray*}
   g'(\al) &=& \sum_{i=1}^n  \left[  - I(x_i\ge 1)\sum_{m=0}^{x_i-1} \frac{1}{(m+\al)^2} \right] + \frac{n}{\al} - \frac{n}{\al+\Bx} \\ [2mm]
   &\ge & - \sum_{i=1}^n  \left[I(x_i\ge 2) \sum_{m=1}^{x_i-1} \frac{1}{m^2} \right] - \frac{ 1}{\al^2}\sum_{i=1}^n  I(x_i \ge 1) + \frac{n\Bx}{\al(\al+\Bx)} \\ [2mm]
   & \ge& a_{10} + a_{11}/\al^2 \teq \; b(\al),
\end{eqnarray*}
where $a_{10} \teq - \sum_{i=1}^n  \left[I(x_i\ge 2) \sum_{m=1}^{x_i-1} m^{-2} \right]$ and $a_{11} \teq - \sum_{i=1}^n  I(x_i\ge 1)$, the US iteration is
\begin{eqnarray*}
   \altI 
   &=& \sol \left\{ g(\alt) + a_{10} (\al-\alt) - a_{11}/\al + a_{11}/\alt = 0, \; \forall \; \al, \alt >0  \right\} \\ [2mm]
   &=& \sol \left\{ a_{10} \al^2 + a_{12}^{(t)} \al - a_{11} =0, \; \forall \; \al, \alt >0  \right\},
\end{eqnarray*}
where $a_{12}^{(t)} \teq g(\alt) - a_{10} \alt + a_{11}/\alt$.

\subsubsection{MLE of $\be$ in Geeta distribution} 

The Geeta distribution is L-shaped like the logarithmic series distribution, Yule distribution and the discrete Pareto distribution, but is far more versatile than them as it has two parameters. Let $X \sim \mbox{Geeta}(\mu, \be)$ with pmf (Consul 1990)
$$
   \frac{\Ga(\be x -1)}{x! \Ga(\be x - x)}  (\mu-1)^{x-1} \left[\mu (\be-1) \right]^{(\be-1)x} \left(\mu \be-1 \right)^{-\be x+1}, \quad x=1, 2, \ldots, \infty,
$$
where $\mu\,(>1)$ is the mean parameter and $\be \,(>1)$ is another parameter. Let $\{X_i\}_{i=1}^n\iid \mbox{Geeta}(\mu, \be)$, the observations be $1, 2, \ldots, K$ with corresponding frequencies $f_1, f_2, \ldots, f_K$, where $n= \sum_{k=1}^K f_k$. Note that when $x=1$, $\Ga(\be x -1)/[x! \Ga(\be x - x)] = 1$, and
\begin{eqnarray*}
   \prod_{i=1}^n \frac{\Ga(\be x_i -1)}{x_i! \Ga(\be x_i - x_i)}
    = \prod_{i=1}^n \frac{(\be x_i-2) (\be x_i-3)\cdots (\be x_i-x_i) } {x_i!}  = \prod_{i=1}^n \frac{ \prod_{j=2}^{x_i} (\be x_i-j)} {x_i!},
\end{eqnarray*}
then the log-likelihood function is given by
\begin{eqnarray*}
   \ell(\mu,\be) &=& n \Bx\Big\{ \log(\mu-1) + (\be-1) \log[\mu(\be-1)] - \be \log(\mu\be-1) \Big\} \\ [2mm]
   & & + \; n\Big[ \log(\mu\be-1) -\log(\mu-1) \Big] + \sum_{k=2}^K f_k \sum_{j=2}^k \log(k \be-j) -\sum_{k=2}^{K} f_k \log(k!),
\end{eqnarray*}
where $\Bx = (1/n) \sum_{k=1}^K kf_k \ge 1$ is the sample mean.

Given $\be$, the MLE of $\mu$ is $\hmu = \Bx$. By replacing $\mu$ with $\Bx$, the log-likelihood function of $\be$ is $\ell(\be) = n \Bx (\be-1) \log[\Bx(\be-1)] - n (\Bx \be-1) \log(\Bx\be-1) + \sum_{k=2}^{K} f_k \sum_{j=2}^k \log(k \be-j) + \mbox{constant}$. Let
$$
   0=\ell'(\be) = n \Bx \Big[\log(\be-1) + \log \Bx- \log( \Bx \be -1) \Big] + \sum_{k=2}^K \sum_{j=2}^k  \frac{f_k}  {\be-j/k} \;\teq\; g(\be),
$$
and define $u_0(\be) = \log(\be-1) + \log \Bx- \log( \Bx \be -1)$,
$$
   u_{kj}(\be) = \frac{\be-1}{\be-j/k}, \quad k=2, \ldots, K; \; j=2, \ldots, k.
$$
Note that $u_0'(\be) = (\Bx-1)/[(\be-1)( \Bx \be -1)] \ge 0$ and $u_{kj}'(\be) = (1- j/k)/(\be-j/k)^2 \ge 0$, then $u_0(\be)$ and all $\{u_{kj}(\be)\}$ are   increasing with respect to $\be$. By applying the fixed--block method, we can construct a $U$-function for $g(\be)$ as
\begin{eqnarray*}
   g(\be) &=& n \Bx \cdot u_0(\be) + \sum_{k=2}^K \sum_{j=2}^k  \frac{f_k}{\be-1} u_{kj}(\be)  \\ [2mm]
   &\sge{\sgn(\be-\Bet)}& n \Bx\cdot u_0(\Bet) + \sum_{k=2}^K \sum_{j=2}^k  \frac{f_k}{\be-1} u_{kj}(\Bet) \teq U(\be|\Bet).
\end{eqnarray*}
Then, the US iteration for calculating the MLE $\hbe$ is
$$
   \betI = 1 - \frac{ \sum_{k=2}^K \sum_{j=2}^k  (\Bet-1)f_k / (\Bet-j/k) } {n \Bx [\log(\Bet -1) + \log \Bx- \log(\Bx\Bet -1)]}.
$$

\end{document}